\documentclass[11pt,a4paper]{article}
\usepackage[T1]{fontenc}
\usepackage[utf8]{inputenc}
\usepackage{lmodern}
\usepackage[a4paper,margin=26mm,headheight=14pt]{geometry}
\usepackage{amsmath,amssymb,amsthm,mathtools}
\usepackage{microtype}
\usepackage{enumitem}
\usepackage{booktabs,array}
\usepackage{tikz}
\usetikzlibrary{arrows.meta,calc,decorations.pathreplacing,positioning}
\usepackage{caption,fancyhdr,xcolor}
\usepackage{placeins}
\usepackage{aliascnt}
\usepackage{hyperref}
\usepackage[nameinlink,noabbrev]{cleveref}
\usepackage{amsmath,amssymb,graphicx,pgfplots}
\pgfplotsset{compat=1.18}
\usetikzlibrary{arrows.meta,decorations.markings}
\definecolor{ink}{HTML}{20364C}
\definecolor{teal}{HTML}{087E8B}
\definecolor{ochre}{HTML}{B46B13}
\definecolor{pale}{HTML}{EDF4F7}
\definecolor{softgray}{HTML}{6A747C}
\hypersetup{colorlinks=true,linkcolor=ink,citecolor=teal,urlcolor=teal,
 bookmarksdepth=2,pdftitle={Thurston algorithm for real entire maps},
 pdfauthor={Sebastian van Strien}}
\setlist[enumerate]{itemsep=0.3em,topsep=0.4em}
\setlist[itemize]{itemsep=0.2em,topsep=0.4em}
\numberwithin{equation}{section}
\newtheorem{theorem}{Theorem}[section]
\newaliascnt{lemma}{theorem}
\newtheorem{lemma}[lemma]{Lemma}
\aliascntresetthe{lemma}
\newaliascnt{proposition}{theorem}
\newtheorem{proposition}[proposition]{Proposition}
\aliascntresetthe{proposition}
\newaliascnt{corollary}{theorem}
\newtheorem{corollary}[corollary]{Corollary}
\aliascntresetthe{corollary}
\theoremstyle{definition}
\newaliascnt{definition}{theorem}
\newtheorem{definition}[definition]{Definition}
\aliascntresetthe{definition}
\newaliascnt{remark}{theorem}
\newtheorem{remark}[remark]{Remark}
\aliascntresetthe{remark}
\newaliascnt{example}{theorem}
\newtheorem{example}[example]{Example}
\aliascntresetthe{example}
\crefname{theorem}{Theorem}{Theorems}
\crefname{lemma}{Lemma}{Lemmas}
\crefname{proposition}{Proposition}{Propositions}
\crefname{corollary}{Corollary}{Corollaries}
\crefname{definition}{Definition}{Definitions}
\crefname{remark}{Remark}{Remarks}
\crefname{example}{Example}{Examples}
\Crefname{theorem}{Theorem}{Theorems}
\Crefname{lemma}{Lemma}{Lemmas}
\Crefname{proposition}{Proposition}{Propositions}
\Crefname{corollary}{Corollary}{Corollaries}
\Crefname{definition}{Definition}{Definitions}
\Crefname{remark}{Remark}{Remarks}
\Crefname{example}{Example}{Examples}
\newcommand{\R}{\mathbb R}
\newcommand{\C}{\mathbb C}
\newcommand{\D}{\mathbb D}
\newcommand{\HH}{\mathbb H}
\newcommand{\cC}{\mathcal C}
\newcommand{\cB}{\mathcal B}
\newcommand{\eps}{\varepsilon}
\newcommand{\dd}{\,\mathrm d}
\newcommand{\abs}[1]{\lvert #1\rvert}
\newcommand{\norm}[1]{\lVert #1\rVert}
\DeclareMathOperator{\dist}{dist}
\DeclareMathOperator{\diam}{diam}
\DeclareMathOperator{\Crit}{Crit}

\DeclareMathOperator{\sep}{sep}
\DeclareMathOperator{\ord}{ord}
\title{\color{ink}\bfseries Thurston algorithm for real entire maps:\\ existence and uniqueness}
\author{Sebastian van Strien}
\date{12 November 2026}
\begin{document}
\maketitle
\thispagestyle{empty}
\begin{abstract}
This paper deals with the real Thurston pullback algorithm for entire maps with finitely
many finite singular values, all real. Given a finite forward-invariant
marking containing the singular values and satisfying an interval
no-Levy-cycle condition, the pullback is shown to have a unique fixed point and it is shown that
every orbit of the Thurston map converges to it. The resulting entire map uniquely realises the
marked dynamics up to positive affine conjugacy. Thus this paper established  
Thurston realisation and uniqueness in the setting of real entire maps 
with a finite postsingular set. 

As usual, the realisation problem reduces to the existence of a fixed point
of the Thurston map $T$ which in the real setting acts on a simplex $W$ of dimension $s-2$
if the required postsingular set consists of $s$ points. 
The key estimate concerns a finite real preimage interval: when its
critical multiplicity is bounded and its endpoint values stay separated
from the singular values unless they are equal to one of them, its hyperbolic diameter
grows at most double-logarithmically as the separation tends to zero.
 An expansive
marking (absence of Levy cycles) then rules out configurations $x\in W$ for which 
the ratio $\|T(x)-x\|_\infty  / \dist(x,\partial W)$ is bounded and which are close 
to the boundary of $W$, i.e. with  $\dist(x,\partial W)$ close to zero. 
Using a version of the Brouwer fixed point theorem,  this immediately  implies 
that the Thurston map $T$ has a fixed point. 
The argument permits infinitely many real and non-real critical points and different
behaviours as the real variable tends to $+\infty$ and $-\infty$.

Additional estimates for fixed labelled inverse branches, which are not strictly necessary for the proof of 
the Main Theorem,  are given in
the appendices.
\end{abstract}
\clearpage
\setcounter{tocdepth}{1}
\tableofcontents
\clearpage

\section{Introduction and statement of results}\label{sec:introduction-new}

Thurston's realisation problem asks when prescribed topological dynamics
can be represented by a holomorphic map. For postcritically finite
branched coverings of the sphere, his characterisation expresses the
answer in terms of obstructions associated with invariant multicurves;
Douady and Hubbard supplied the proof in \cite{DH93}. The transcendental
setting requires additional arguments. Hubbard, Schleicher and
Shishikura established a characterisation for postsingularly finite
topological exponential maps \cite{HSS09}; results for small
postsingular sets are developed in \cite{Pro24}.

Given that the general transcendental setting has been open for so long, 
it seems to make sense to address this question under the additional assumption 
that the maps are real. 

For real interval maps the order structure supplies a different route.
Milnor and Thurston \cite{MT88} develop kneading theory for interval
maps, including the unimodal setting on compact intervals. Using this,  
one can derive the existence of unimodal maps of an admissible 
combinatorial type from an intermediate value theorem. 
The analogous result for full multimodal families cannot be derived
in this manner, but was proved in \cite[Chapter~II, \S4]{dMvS93} using a real Thurston operator and a
fixed-point criterion on an open simplex. The present paper develops
that method for real entire maps with finitely many singular values,
all real. The maps considered here can have infinitely many real critical
points, and their behaviour as $x\to-\infty$ and $x\to+\infty$ can differ.

The new difficulty is that an entire map in a fixed natural family can
develop large derivatives even when its singular values stay in a
bounded interval. A proof based on compactness of the whole normalised
family would therefore impose an unnecessarily strong requirement.
Here the needed estimates are established along configurations for
which the Thurston displacement $d(Tx,x)$ is at most a fixed multiple $K$ of the
smallest marked gap $m(x)=\min (x_{i+1}-x_i)$. 
Hyperbolic length estimates control the affine
normalisation of the pullback interval before any compactness of the
unnormalised maps is invoked.

\subsection{The real natural family}
A nonconstant entire function \(f\) is \emph{real} if
\(f(\bar z)=\overline{f(z)}\). Its finite singular set \(S(f)\) is the
closure of the critical values and finite asymptotic values with the 
equivalent covering property
\begin{equation}\label{eq:covering}
 f:\C\setminus f^{-1}(S(f))\longrightarrow\C\setminus S(f)
 \quad\text{is a holomorphic covering.}
\end{equation}
A plane homeomorphism is called \emph{real symmetric} if it commutes
with complex conjugation. All plane homeomorphisms used below preserve
orientation, and their restrictions to \(\R\) are increasing.

Fix a real entire map \(f_0\) with
\[
 S(f_0)=\{v_{1,0}<\cdots<v_{q,0}\}\subset\R,
 \qquad 1\le q<\infty.
\]
Its real natural family \(\mathcal N(f_0)\) consists of the real entire
maps admitting a representation
\begin{equation}\label{eq:natural}
 f=\phi\circ f_0\circ\psi,
\end{equation}
where \(\phi,\psi\) are real-symmetric quasiconformal plane maps
increasing on \(\R\).
The equivalences transport the order and local degree of the real
critical points, the singular-value labels, and the labels and
orientations of the real branches. A \emph{real branch} is a component
of \(\R\setminus\Crit_{\R}(f)\), where
\(\Crit_{\R}(f)=\{t\in\R:f'(t)=0\}\). An inflection critical point (an odd-local-degree critical point at
which monotonicity does not reverse) remains a critical point.

For the realisation theorem it is assumed that \(f_0|\R\) is not
monotone. This guarantees a real local extremum and is consistent with
the expansive finite diagrams obtained from the interval no-Levy-cycle
condition below. The quantitative estimates also apply to a fixed finite
marking that satisfies \eqref{eq:noflat} and \eqref{eq:expansive}, even
when $f_0|\R$ is monotone, provided the required critical labels exist. No restriction is
placed on the complex arrangement of preimages of real intervals.

\subsection{Main results}
The Main Theorem below gives existence, uniqueness, and convergence for
an expansive finite marking. Three independent statements drive its
proof: \cref{thm:boundary-fixedpoint} is the displacement-near-boundary
fixed-point criterion; \cref{ana:diameter} bounds the hyperbolic diameter
of a finite preimage interval \emph{without} assuming a pullback
displacement bound; and \cref{gap:separation} says that for each fixed
relative-displacement constant $K$, the corresponding configurations
have a positive distance from the boundary of $W$. The lower bound
depends on $K$.

\subsection{Finite marked realisation}
Let \(g\) be a topological model for \(f_0\), with singular set
\(S(g)\) as in \cref{def:model-new}. Write
\[
 P(g)=\bigcup_{n\ge0}g^n(S(g)).
\]
I consider a finite forward-invariant set \(X\supset P(g)\) with
at least two points. Its elements are the \emph{marked points} of the
model. The coordinates $x_i$ in the Main Theorem are their ordered
target counterparts; the preimages $\widehat y_i$ in
\eqref{eq:affine} lie in the source plane before affine rescaling.
The absence of real Levy cycles means that every interval joining two
distinct points of \(X\) eventually meets a real critical point of $g$; the precise convention is given in \cref{def:levy-new}.
Intersection at an endpoint is allowed. After adding the finitely many
critical points in the convex hull of \(X\), as in
\cref{lem:critical-saturation}, write
\[
 X=\{u_1<\cdots<u_s\},\qquad
 g(u_i)=u_{\pi(i)}.
\]
Thus \(\pi\) is a self-map of the finite ordered label set. Its
expansiveness, defined in \cref{def:pi-expansive}, is the finite
combinatorial form of the real Levy condition.

\begin{theorem}[Main Theorem]\label{thm:realization-new}
Let $f_0$ be a real entire map with a finite nonempty singular set
$S(f_0)\subset\R$ and with a nonmonotone real restriction. Let $g$
be a labelled real topological model for $f_0$ in the sense of
\cref{def:model-new}. Suppose $X\subset\R$ is finite, has at least
two points, satisfies $g(X)\subset X$ and $P(g)\subset X$, and has no
real Levy cycles as in \cref{def:levy-new}. After adjoining the real
critical points in the convex hull $\operatorname{conv}X$ of $X$ to the set $X$, fix the inherited inverse
branch and critical-point labels. Then the
Thurston map on
\[
 W=\{x=(x_1,\ldots,x_s):0=x_1<x_2<\cdots<x_s=1\},
\]
defined in Section~\ref{sec:thurston}, 
has a unique fixed point \(x^*\), and every orbit of the Thurston map
converges to \(x^*\). Here uniqueness holds within the fixed natural family \(\mathcal N(f_0)\) with the prescribed labels. 

There is a real entire map \(F_*\in\mathcal N(f_0)\) with the
prescribed branch data such that
\[
 F_*(x_i^*)=x_{\pi(i)}^*\quad(1\le i\le s),
 \qquad S(F_*)=\{x_{i_1}^*,\ldots,x_{i_q}^*\},
\]
where \(i_1<\cdots<i_q\) are the indices of the singular values in
\(X\).
The correspondence \(u_i\mapsto x_i^*\) is an order-preserving
conjugacy of the marked finite dynamics. The marked realisation is
unique with the endpoint normalisation \(x_1^*=0\), \(x_s^*=1\);
equivalently, before that normalisation it is unique up to positive
real affine conjugacy.
\end{theorem}

The endpoint normalisation in this theorem is a \emph{dynamical}
normalisation: the same coordinates are used in the source and target
of \(F_*\). It is distinct from the source normalisation used to select
the maps \(f_v\) in the natural family. Precisely, \eqref{eq:affine} defines $M_x(z)=\widehat y_1+
(\widehat y_s-\widehat y_1)z$ and $F_x=f_x\circ M_x$; at the fixed
point this is the map $F_*$. Keeping these two maps separate is essential.

This paper proves the finite-marked theorem. The extension from finite
to infinite marked sets for interval maps is discussed in
\cite[Chapter~II, \S4]{dMvS93}; no such extension is needed for the
arguments here.

\subsection{An example with a critical and an asymptotic value}
Consider
\[
 f(x)=(ax+b)e^{\lambda x}+A,\qquad a>0,\quad\lambda>0.
\]
Its finite singular values are the asymptotic value \(A\) and the
critical value at \(c=-1/\lambda-b/a\). This unimodal
family has {\em two} singular values whose dynamics can be prescribed.
For example, let us choose $f$ so that the asymptotic value is a fixed point 
and so that the other singular value is zero: 
\[
 f(-\infty)=A=f(A)=1,\qquad f(c)=0.
\]
These equations give
\[
 b=-a,\qquad c=1-\frac1\lambda,\qquad
 a=a(\lambda)=\lambda e^{1-\lambda},
\]
and hence
\begin{equation}\label{eq:exponential-example-new}
 f_\lambda(x)=1+\lambda e^{1-\lambda}(x-1)e^{\lambda x}.
\end{equation}
One may seek, for example, the marked order
\[
 0=f_\lambda(c)<f_\lambda^3(c)=f_\lambda^4(c)
 <f_\lambda^2(c)<1.
\]
The competing estimates already appear in this example:
\[
 1-f_\lambda^2(c)=\lambda e^{1-\lambda},
 \qquad
 \sup_{x\le1}|f_\lambda'(x)|=e\cdot \lambda.
\]
Indeed, with \(u=\lambda(x-1)\le0\), the derivative is
\(e\lambda(1+u)e^u\), whose absolute value has maximum
\(e\lambda\) on this half-ray. Here ``target'' points are the coordinates $x_i$ defined below;
``unscaled'' derivatives are derivatives of the source-normalised map
$f_\lambda$ before its pullback interval is rescaled.
The general proof compares these polynomial versus exponential effects after the affine pullback
normalisation. Later, in \cref{def:source-normalization-new} such maps will be rescaled so that they are normalised. 

\subsection{History of the algorithm}
Thurston announced his theorem about the realisability and uniqueness of 
topological branched covering maps around 1982,  but the proof was only published in 1993
by Douady and Hubbard \cite{DH93}.  The spider algorithm of Hubbard and Schleicher \cite{HS94} gives a somewhat different proof in the setting of general polynomials. For general complex entire maps only  the family $\lambda e^z$  has been treated  \cite{HSS09}, although some results were developed for small postsingular sets in \cite{Pro24}.

The real Thurston algorithm has also an extensive history in polynomial and interval dynamics.
Existence and uniqueness in the real polynomial setting of course already follow from Thurston's theorem; see \cite{DH93}. 
The existence of a fixed point of the Thurston map was later generalised to families of $\ell$-modal smooth interval maps under natural conditions in my book dating back to 1993 with Welington de Melo \cite[Section II.4]{dMvS93}, using a displacement fixed-point Theorem, which is restated as  
 \cref{thm:boundary-fixedpoint} in this paper. 
The special case of families of real polynomials with all critical points real and non-degenerate is then an immediate consequence.   The pictures in \cite{dMvS93} were produced by an implementation of this algorithm.
Bonifant, Milnor and Sutherland \cite{BMS} describe an effective numerical
procedure for constructing critically finite real polynomials with only real critical points with 
specified combinatorics, relying on the existence of the Thurston fixed point obtained in 
\cite{DH93}. The spider algorithm of Hubbard and Schleicher \cite{HS94} treats
polynomials without a reality assumption. Related realisation problems
for multimodal circle maps are considered in \cite{dFdMSV26}.
The current paper solves the Thurston realisation problem for arbitrary natural families of real entire maps. 

\subsection{Outline of the proof}

The idea of the proof goes as follows. First parameterise the source-normalised natural family by its
ordered singular values and define the finite marked Thurston map which acts on the 
set $W$ of ordered points in the postsingular set (which is a simplex in the real case). 
Next the displacement criterion \cref{thm:boundary-fixedpoint} is shown. This theorem asserts that 
if a continuous map $T$ on a bounded open convex subset $W$ of a finite-dimensional normed space satisfies
 \begin{equation} \|T(x)-x\|_\infty  / \dist(x,\partial W) \to \infty\mbox{ as  }\dist(x,\partial W)\to 0
 \label{eq:displacement-intro} 
 \end{equation}
   then it has a fixed point. In our context $ \dist(x,\partial W)\approx m(x)$
 where $m(x)$ is the smallest gap between consecutive target coordinates, i.e. $m(x)=\min_i(x_{i+1}-x_i)$. 
 So in order to show that the Thurston map has a fixed point, it is enough to 
 show that there exists no  $x\in W$ with $m(x)$ arbitrarily small
 while  $\|T(x)-x\|_\infty\le K m(x)$ for a fixed $K$.   
 In the analytic part  of this paper, the proof starts with a point $x=(x_1,\dots,x_s) \in W$, 
 taking the map $f_x$ with the appropriate singular values, 
 and choosing the vector $(\widehat y_1,\dots, \widehat y_s)$ 
 so that $f_x(\widehat y_i)=x_{\pi(i)}$, see \cref{def:thurston-new}. 
The finite interval $J=[\widehat y_1,\widehat y_s]$ spanned
by these unscaled inverse images could be very large. 
The  finite-interval estimate \cref{ana:diameter}  then does most of the work, by  giving the following upper bound for the hyperbolic 
length of $J$ in $U$:
\[
 \diam_U J \le  C_0+C_1\log\bigl(1+\log(1/m(x))\bigr),
 \]
where   $U$ is the component of $f^{-1}D(0,4)$ containing the interval $J$.
Using this,  Koebe estimates in \cref{ana:scaled},  bound the first and second derivatives of
the affine-normalised map by a power of \(1+\log(1/m(x))\). This in turn gives
a stretched exponential upper bound for $m(x)$ in terms of the ratio of the length and some hyperbolic radius, see \cref{ana:stretch}. 
Expansiveness, i.e. the absence of Levy cycles, uses the presence of critical points   in the definition of the points $\hat y_i$ to 
turn the relative displacement bound $\|T(x)-x\|_\infty\le K m(x)$ into an inequality of the form
\[
 1\le C m(x)\bigl(1+\log(1/m(x))\bigr)^B,
\]
see \cref{gap:inequality}. This gives the lower bound  $m(x)\ge m_*(K)$, see \cref{gap:separation}, 
and thus 
prevents marked gaps from collapsing. This means that there are simply no points 
$x$ so that $\|T(x)-x\|_\infty\le K m(x)$ and with $m(x)$ arbitrarily small. 
Appealing to   the fixed-point criterion \cref{thm:boundary-fixedpoint}
then gives the existence of a fixed point. 

The discussion in \cref{sec:degeneration-new} explains why
sharply concentrated members of a natural family, including the
noncompact degeneration that in previous drafts of this paper led to an additional condition (which was called A2), do not
obstruct the proof in this paper and in fact motivated it. Finally, \cref{thm:convergence-new} proves uniqueness
and convergence of the Thurston algorithm. Stronger a priori derivative
and inverse-height estimates, which are not needed for the realisation
theorem but which are of independent interest, are collected in the appendices.

% Dependency graph for the proof of the Main Theorem.
% No extra TikZ libraries are needed beyond those already loaded in the manuscript.
\begin{figure}[p]
\centering
\resizebox{0.98\textwidth}{!}{%
\begin{tikzpicture}[
 >=Stealth,
 common/.style={rounded corners=2.5pt,draw=ink,fill=pale,align=center,inner sep=5pt,text width=4.8cm,font=\small},
 analytic/.style={rounded corners=2.5pt,draw=teal!85!black,fill=teal!7,align=center,inner sep=5pt,text width=4.8cm,font=\small},
 comb/.style={rounded corners=2.5pt,draw=ochre!90!black,fill=ochre!8,align=center,inner sep=5pt,text width=4.8cm,font=\small},
 metric/.style={rounded corners=2.5pt,draw=softgray,fill=softgray!8,align=center,inner sep=5pt,text width=4.8cm,font=\small},
 result/.style={rounded corners=2.5pt,very thick,draw=ink,fill=ink!7,align=center,inner sep=6pt,text width=5.3cm,font=\small},
 note/.style={rounded corners=2.5pt,draw=softgray,dashed,fill=white,align=center,inner sep=5pt,text width=7.4cm,font=\small},
 arr/.style={->,line width=.7pt,draw=softgray},
 mainarr/.style={->,line width=1.05pt,draw=ink},
 dasharr/.style={->,dashed,line width=.65pt,draw=softgray},
 lab/.style={font=\scriptsize,text=softgray,fill=white,inner sep=1.2pt}
]

% Common setup
\node[common] (coords) at (-3.2,0)
 {\cref{thm:natural-uniqueness,thm:natural-parameters}\\
  Natural-family uniqueness and singular-value coordinates};
\node[comb] (comb) at (3.2,0)
 {\cref{lem:critical-saturation,lem:expansive-levy}\\
  Critical saturation and\\no-Levy $\Longleftrightarrow$ expansiveness};
\node[comb,dashed,text width=4.2cm] (necess) at (3.2,2.2)
 {\cref{prop:levy-necessary}\\Necessity of no-Levy};
\draw[dasharr] (necess) -- node[lab,right]{motivation} (comb);

\node[common] (T) at (0,-2.5)
 {\cref{lem:thurston-continuous}\\
  The Thurston map $T:W\to W$ is well defined and continuous};
\draw[arr] (coords) -- (T);
\draw[arr] (comb) -- (T);

% headings
\node[font=\bfseries\large,text=teal!80!black] at (-4.7,-4.3) {Existence};
\node[font=\bfseries\large,text=softgray] at (6.0,-4.3) {Uniqueness and convergence};

% Existence branch
\node[analytic] (toolkit) at (-5.8,-6.0)
 {\cref{ana:disk-topology,ana:singleton,ana:range,ana:arch}\\
  Pullback topology and double-logarithmic target arches};
\node[analytic] (diam) at (-5.8,-8.7)
 {\cref{ana:diameter}\\Finite-interval hyperbolic diameter bound};
\node[analytic] (scaled) at (-5.8,-11.4)
 {\cref{ana:scaled}\\Polylogarithmic bounds for $\|F'\|$ and $\|F''\|$\\{\footnotesize using \cref{ana:koebe,ana:radius}}};
\draw[arr] (toolkit) -- (diam);
\draw[arr] (diam) -- (scaled);

\node[comb] (gap) at (0,-8.7)
 {\cref{gap:interval,gap:inequality}\\
  Critical-multiplicity control and the expansive-gap inequality\\{\footnotesize uses expansiveness from \cref{lem:expansive-levy}}};
\draw[arr] (T) -- node[lab,right]{\cref{def:admissible}} (gap);

\node[result] (sep) at (-2.5,-14.1)
 {\cref{gap:separation}\\
  Bounded relative displacement\\$\Longrightarrow$ uniform marked-gap separation\\
  {\footnotesize (proves the estimate stated in \cref{thm:apriori})}};
\draw[arr] (scaled) -- (sep);
\draw[arr] (gap) -- (sep);

\node[metric] (fixedpt) at (0.3,-11.4)
 {\cref{lem:mesh-distance,thm:boundary-fixedpoint}\\
  Mesh--boundary comparison and the displacement fixed-point criterion};

\node[result] (exist) at (0,-17.3)
 {\cref{cor:existence-new,lem:thurston-realisation}\\
  Existence of a fixed point and realisation of the marked dynamics};
\draw[mainarr] (sep) -- (exist);
\draw[arr] (fixedpt) -- (exist);

% Uniqueness / convergence branch
\node[metric] (marking) at (6.0,-6.0)
 {\cref{lem:canonical-real-marking}\\Canonical real Teichm\"uller marking};
\node[metric] (lift) at (6.0,-8.7)
 {\cref{lem:teich-lift-new}\\Nonexpansiveness and the equality case under lifting\\{\footnotesize for the Thurston map; uses \cref{thm:natural-uniqueness}}};
\node[metric] (strict) at (5.4,-11.4)
 {\cref{prop:strict-iterate-new}\\A suitable iterate strictly decreases $d_{\mathrm T}$};
\node[metric] (proper) at (7.6,-14.1)
 {\cref{lem:teich-proper-new}\\Closed $d_{\mathrm T}$-balls in $W$ are compact};
\node[result] (conv) at (6.0,-17.3)
 {\cref{thm:convergence-new}\\Uniqueness of the fixed point and global convergence};
\draw[arr] (marking) -- (lift);
\draw[arr] (lift) -- (strict);
\draw[arr] (strict) -- (conv);
\draw[arr] (proper) -- (conv);
\draw[arr] (exist) -- (conv);

% Final theorem
\node[result,text width=6.5cm] (main) at (3.0,-20.6)
 {\cref{thm:realization-new}\\[1mm]
  Existence + realisation + uniqueness + convergence};
\draw[mainarr] (exist) -- (main);
\draw[mainarr] (conv) -- (main);

% Appendix-only note
\node[note] (appendix) at (0,-23.6)
 {\cref{gap:separation} also feeds the appendix chain
  \cref{glob:core} $\to$ \cref{glob:radius-lower,glob:all-derivatives}
  $\to$ \cref{glob:height}, summarized in \cref{thm:additional-apriori}.
  These results are not used in the proof of \cref{thm:realization-new}.};
\draw[dasharr] (sep.south) -- ++(0,0) -| (appendix.north west);

\node[font=\scriptsize,text=softgray] at (0,-25.4)
 {Solid arrows: proof dependencies. \quad Dashed arrows: motivation or appendix-only consequences.};
\end{tikzpicture}
}%
\caption{Dependency graph for the principal results. Solid arrows indicate proof dependencies. Dashed arrows indicate motivation or appendix-only consequences. The left branch gives existence; the right branch gives uniqueness and convergence.}
\label{fig:result-dependencies}
\end{figure}
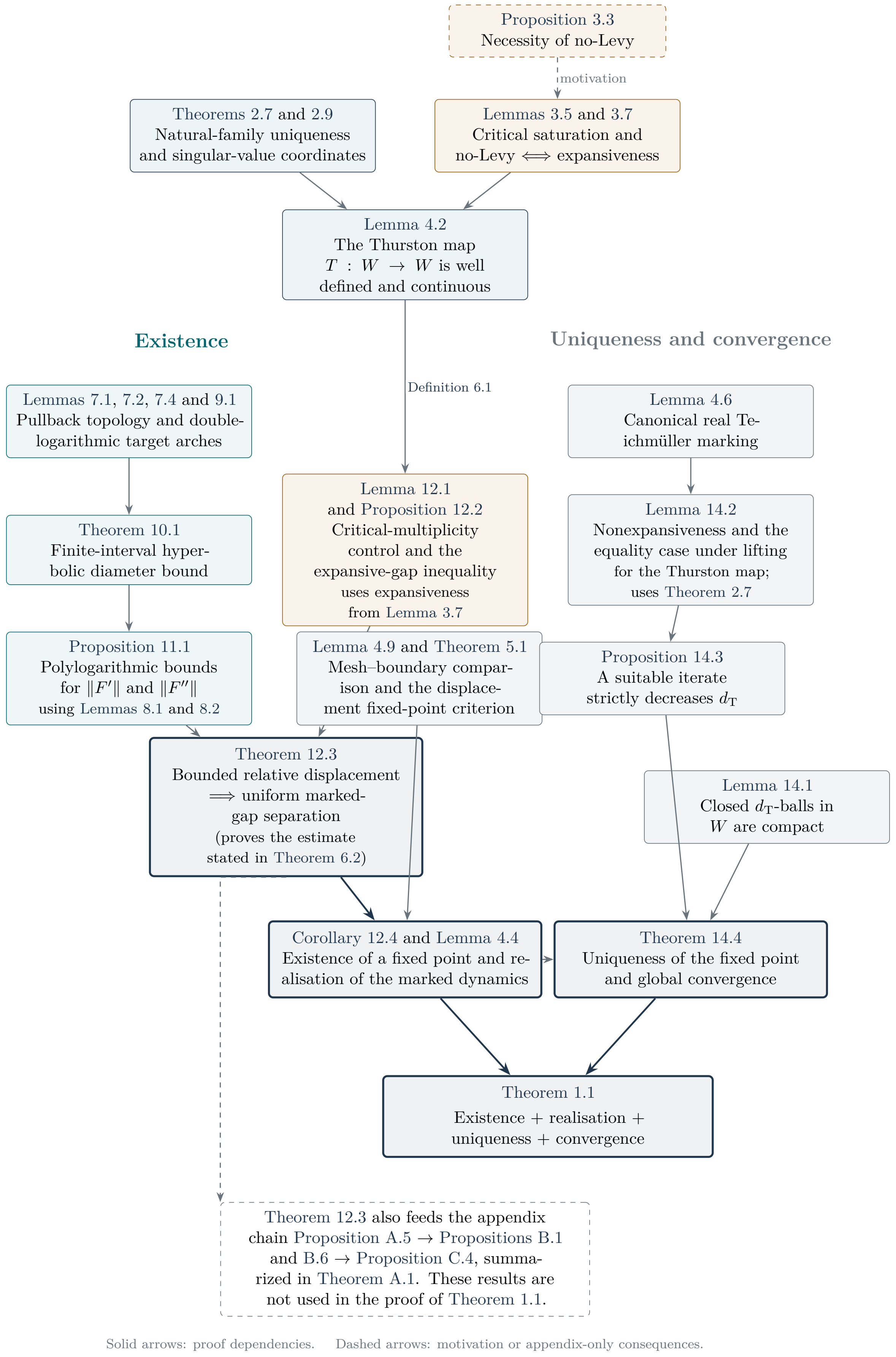

% Proposed replacement for \subsection{Outlook} (Section 1.7) of thurston_real_8nov2026.tex.
% Replace everything from "\subsection{Outlook}" up to (not including)
% "\subsection{Acknowledgements and use of artificial intelligence}".
% Add the bibitem at the end of this file, with the arXiv number once known.
% The sentence on Appendix E assumes appendix-axiomatic-existence.tex is included;
% otherwise delete that sentence (marked below).

\subsection{Outlook}

In a companion paper \cite{vS-complex} I treat the corresponding problem for entire maps of the complex plane. There it is shown that a marked entire Thurston map of stable parabolic type is Thurston equivalent to a postsingularly finite entire function if and only if it has no Levy cycle. Moreover, the realisation is unique up to affine conjugacy, and every orbit of the Thurston pullback converges to it.

The proof follows the strategy of the present paper. Here, in the absence of a fixed point, \cref{thm:boundary-fixedpoint} yields the displacement estimate \eqref{eq:displacement-intro}. Its complex counterpart is a Wolff--Denjoy type theorem for nonexpanding self-maps of Teichm\"uller space: if the pullback has no fixed point, then one measured foliation has extremal length that grows by at most a fixed factor under all iterates. The hyperbolic-diameter bound of \cref{ana:diameter} and the expansiveness argument of \cref{sec:gaps} are replaced by estimates for quadratic differentials and by the covering degrees of lifted annuli. The approach of Hubbard, Schleicher and Shishikura \cite{HSS09} has no analogue of \eqref{eq:displacement-intro}.

% History of the companion paper, following its section "Use of artificial intelligence".
% CHECK: result numbers are those of complex-v12.tex (build of 20 Nov 2026): Theorem 2.1 = l:record,
% Section 4 = sec:closed (vertical cylinders), Section 10.1 = sec:fromreal. Update if they change.
The companion paper grew out of the present one. After completing the first version of this paper, I used OpenAI's ChatGPT to explore a complex analogue of the displacement condition $\|T(x)-x\|_\infty\le Km(x)$ behind \eqref{eq:displacement-intro}, where $m(x)$ is comparable to $\dist(x,\partial W)$. Through several iterations this led to the Wolff--Denjoy type theorem above \cite[Theorem~2.1]{vS-complex}. I then asked whether that theorem gives a different proof of the real case, with geometric methods in place of the analytic results on real entire maps. For real-symmetric maps its cylinder case suffices \cite[Sections~4 and~10.1]{vS-complex}. A sharpening of the near-equality argument of \cite{HSS09}, building on the exact equality case of Douady and Hubbard \cite{DH93}, then handles the remaining case and gives the general theorem.

The theorem proved here is not a formal consequence of the complex one.
\begin{itemize}
\item It realises the marked dynamics inside the prescribed natural family $N(f_0)$, with uniqueness up to positive affine conjugacy.
\item Its hypothesis is the absence of real Levy cycles. By \cref{lem:expansive-levy} this is expansiveness of a finite combinatorial map, and it can be checked by a finite computation.
% CHECK: the next sentence uses that the real realisation is a planar Thurston
% equivalence relative to the marked set (cf. the canonical real marking, Lemma 4.6).
\item Together, the two results show that for these real maps the absence of real Levy cycles excludes Levy cycles in the plane.
\end{itemize}
The argument given here is real-variable apart from \cref{thm:natural-parameters} and \cref{ana:diameter,ana:scaled}. It is quantitative, and the Thurston map acts on the simplex $W$, where it can be implemented directly.
% Delete the next sentence if Appendix E is not included.
Its existence part extends to families of real maps that are not entire; see \cref{app:axiomatic}. Uniqueness, in contrast, uses the complex structure in an essential way.

% ---- bibitem to add ----
% \bibitem{vS-complex}
% Sebastian van Strien,
% \emph{The Thurston algorithm for general entire maps on the complex plane:
% existence and uniqueness},
% preprint, 2026, arXiv:2611.xxxxx.

%\subsection{Outlook}
%
%It seems to me that the current proof is interesting, even if once a proof 
%of the general complex Thurston Theorem for entire maps is established. 
%Moreover,  the current proof may in fact give insight into how to prove the general complex case. 
%The mechanism of the proof in this paper seems quite different from the one used in \cite{HSS09}
%because in that paper there seems to be no analogue of the inequality in \cref{eq:displacement-intro}. 
%

\subsection{Acknowledgements and use of artificial intelligence}
I thank Kostya Drach, Sergei Merenkov, Nikolai Prochorov, and Lasse Rempe
for helpful comments and discussions. I have given talks about this paper, in an increasingly 
more general reincarnation,  at Warwick in May 2025, at the Newton Institute in  Cambridge in July 2026,  Fudan University and Hebei University in China in October 2026. I would like to thank 
the audiences for their useful questions. Originally, the proof I had used some compactness 
for certain maps at the boundary of $W$, but only under an additional assumption  (which I called 
assumption A2).  The example in \cref{sec:degeneration-new} led to the realisation that 
the pivotal techniques of \cref{ana:arches} and \cref{ana:finite-interval} could be interpreted 
in a way that would enable me to prove the result for general real entire maps.

I used OpenAI's ChatGPT during the preparation of this manuscript as a
research and writing assistant. Its roles included exploring proof
strategies and counterexamples, checking intermediate calculations,
helping reorganise and edit the exposition, assisting with LaTeX and
TikZ, and helping with literature searches. I decided which suggestions
to retain and revised the resulting text. I take full responsibility for
the mathematical statements, proofs, references, and any remaining
errors. AI output was not used as a mathematical source; external results
are attributed to the cited literature.

I would like to also thank Bernat Bassols-Cornudella, who helped me to formalise the first ArXiv version
via Lean 4, conditional on 25 classical statements. Please email me if you would like to see the
Formalisation Summary. 

\section{The natural family and singular-value coordinates}
\label{sec:natural-parameters-new}
 
The ordered singular values provide \(q\) parameters, while positive
affine changes of the source provide two further parameters. This is
the real form of the finite-dimensional deformation theory of
Eremenko and Lyubich \cite[\S3]{EL92}; see also
\cite{FvS25}. The statements needed here can be proved directly.

\subsection{Source normalisation}
A \emph{source anchor point} is a specified critical point of $f_0$,
identified throughout the family by the real increasing natural
equivalence; the anchor label fixes which critical point is meant,
even if other critical points have the same critical value. Here 
one of the following normalisations is used. 

\begin{definition}[Normalised natural family]
\label{def:source-normalization-new}
\begin{enumerate}[label=\textup{(N\arabic*)}]
 \item Fix two distinct labelled real critical source anchor points
       \(A=\{c_1,c_2\}\), \(c_1<c_2\), and require
       \(\psi(c_j)=c_j\) for \(j=1,2\).
       \label{item:twoanchors}
 \item Fix one labelled real critical source anchor point \(A=\{c\}\), require
       \(\psi(c)=c\), and fix its first nonconstant Taylor term:
       \begin{equation}\label{eq:jet}
        f(c+z)-f(c)=A_\ell z^\ell+O(z^{\ell+1}),
        \qquad \ell\ge2,\qquad A_\ell\ne0.
       \end{equation}
       The coefficient \(A_\ell\) may be taken to be the corresponding
       coefficient of \(f_0\).
       \label{item:jetanchor}
\end{enumerate}
The source-normalised family is denoted by \(\mathcal N_0(f_0)\).
\end{definition}

The second normalisation is available even when \(f_0\) has several
real critical points. The word ``nonconstant'' in \eqref{eq:jet}
matters: the critical value is allowed to vary with the target
parameters. The coefficient sign is preserved by the increasing real
equivalence.

Both normalisations can be achieved by a positive affine source
change. In the first case, send the two prescribed source coordinates
to the transported critical points. In the second, translate the
critical point to \(c\), and then choose a positive scale so that its
coefficient becomes \(A_\ell\). Neither operation changes the singular
values.

For uniqueness of the holomorphic representative, two fixed labelled
points of \(f_0^{-1}(S(f_0))\cap\R\) would suffice, whether or not
both are critical. This critical-point  choice is made because
this also supplies the conformal-radius estimates in the analytic part
of the paper.

\subsection{Examples of natural families}
\begin{example}[The family of \(-ze^z\)]\label{ex:ze-natural-new}
For \(f_0(z)=-ze^z\), the real natural family consists of
\[
 f(z)=-(z-\alpha)e^{\lambda z+\beta}+A,
 \qquad \lambda>0,\quad\alpha,\beta,A\in\R.
\]
Indeed, an increasing affine target map sends the two singular values
\(-1/e,0\) to any prescribed ordered pair. The uniqueness statement
below then shows that all representatives are obtained from this map
by positive affine source changes. The displayed formula is a
reparameterisation of those source and target changes. There is
only one real critical point, so the two-critical-anchor choice (N1)
does not apply; the critical point and its quadratic coefficient may
instead be fixed using (N2).

Using our algorithm, it is possible to specify admissible dynamics of the singular values 
within this family. For example, take \( 
f_\beta(z)=Aze^z+\beta\). Our Main Theorem asserts that there exist parameters
so that the asymptotic value $\beta$ has period two \(
\beta\longmapsto0\longmapsto\beta\) and the critical point $c$ has period three
\(c \longmapsto v\longmapsto w\longmapsto c\), see Figure~\ref{fig:ze-z-cobweb}. 
% Paste this figure environment into your article.
% Required package: graphicx. Place ze_z_cobweb.pdf beside the article.
\begin{figure}[htbp]
  \centering
  \includegraphics[width=\textwidth]{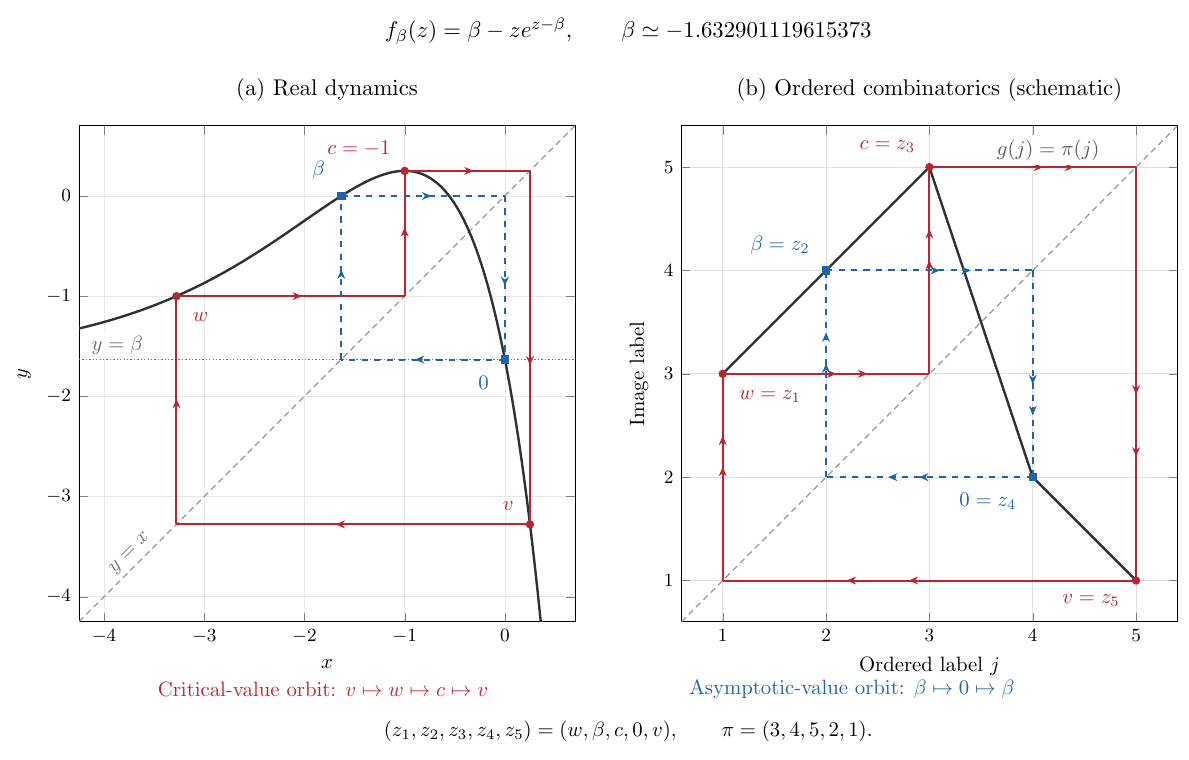}
  \caption{Real dynamics of $f_\beta(z)=\beta-ze^{z-\beta}$,
    with $\beta\simeq-1.632901119615373$.
    Here $c=-1$, $v=f_\beta(c)$, and $w=f_\beta(v)$.
    The red solid cobweb follows the critical-value cycle
    $v\mapsto w\mapsto c\mapsto v$; the blue dashed cobweb
    follows the asymptotic-value cycle $\beta\mapsto0\mapsto\beta$.
    (a) Numerical coordinates, with the graph of $f_\beta$, the diagonal,
    and the dotted horizontal asymptote $y=\beta$.
    (b) The same orbits on ordered labels
    $(z_1,\ldots,z_5)=(w,\beta,c,0,v)$.
    The black polygonal graph is the linear interpolation of
    $g(j)=\pi(j)$, where $\pi=(3,4,5,2,1)$.
    Markers lie on the graph; crossings of colored paths are not
    additional orbit points.}
  \label{fig:ze-z-cobweb}
\end{figure}
 \end{example}

\begin{example}[Exponential polynomials]
\label{ex:exp-polynomial-new}
Let \(f_0=e^Q\), where \(Q\) is a nonconstant real polynomial. Every
\(f\in\mathcal N(f_0)\) has the form
\[
 f=e^P+b,
\]
where \(P\) is a real polynomial of the same degree as \(Q\), and
\(b\) is the transported omitted value. To see this, \(f-b\) is an
entire nonvanishing function, so it has an entire logarithm \(P\).
Since \(f(x)-b>0\) on \(\R\), the logarithm can be chosen real.
The polynomial growth bounds for plane quasiconformal maps show that
\(f\) has finite order. More explicitly, the polynomial growth bounds for quasiconformal
plane maps give $\sup_{|z|\le R}\Re P(z)\le C(1+R)^N$ for some
$C,N$: indeed $\Re P(z)=\log|f(z)-b|$. Subtract $P(0)$ and apply
the Borel--Carath\'eodory inequality on the disk of radius $2R$
to bound $\sup_{|z|\le R}|P(z)-P(0)|$ by $C'(1+R)^N$.
Cauchy's coefficient estimates as $R\to\infty$ then force all
Taylor coefficients of degree greater than $N$ to vanish. The total multiplicity
of the critical points is preserved by natural equivalence. Since
\(f'=P'e^P\), this gives \(\deg P=\deg Q\).

In particular,
\begin{equation}\label{eq:quartic-family-new}
 f_\lambda(z)=\exp\bigl(-\lambda(z^2-1)^2\bigr),
 \qquad\lambda>0,
\end{equation}
belongs to the natural family of \(f_1\), has singular set
\(\{0,e^{-\lambda},1\}\), and is source-normalised by
\(A=\{-1,1\}\), its two real critical source anchor points. An explicit natural equivalence and the degeneration
of this family are examined in
\cref{sec:degeneration-new}.
\end{example}

\begin{example}[Polynomial-exponential maps]
Let $f_0(z)=P(z)e^{Q(z)}$, where $P,Q$ are real polynomials and
$Q$ is nonconstant. The asymptotic value $0$ is a singular value.
Every member of its real natural family has the form
$f(z)=b+R(z)e^{H(z)}$, where $b=\phi(0)$ is the labelled finite
asymptotic value and $R,H$ are real polynomials. Indeed, the inverse
image of $b$ consists of finitely many points, with the same
multiplicities as the zeros of $P$; factor these zeros from $f-b$
to obtain $f-b=Re^h$ with $h$ entire and real on $\R$ (up to a
constant sign absorbed in $R$).
Finite order and the Borel--Carath\'eodory--Cauchy argument of the
preceding example make $h=H$ a polynomial. In particular,
$\deg R=\deg P$. The number of logarithmic tracts over $0$,
which is preserved by the labelled natural equivalence, is
$\deg Q$ for $Pe^Q$ and $\deg H$ for $b+Re^H$; hence
$\deg H=\deg Q$. This is a statement of \emph{necessary form},
not a claim that every pair $(R,H)$ of these degrees belongs to the
same natural family.
\end{example}

% Drop-in replacement for Example 2.5.
\begin{example}[A conjugate pair of nonreal critical points]
\label{ex:nonreal-critical-pair}
Consider
\(
 F(z)=\frac{z^4}{4}+\frac{z^2}{2}.
\)
with simple critical points \( 0,i,-i\) and with \(S(F)=\left\{-\frac14,0\right\}\subset\R\). 
In particular, the two nonreal conjugate critical points have the
same real critical value.
%
%It is useful first to retain the positive affine source freedom.  For
%\(u\in\R\), \(\rho>0\), and prescribed ordered singular values
%\(v_-<v_+\), put
%\begin{equation}\label{eq:quartic-full-family}
% f_{u,\rho,v_-,v_+}(z)
% =v_++4(v_+-v_-)
% F\left(\frac{z-u}{\rho}\right).
%\end{equation}
%Equivalently,
%\[
% f_{u,\rho,v_-,v_+}(z)
% =v_+ +(v_+-v_-)
% \left(
%  \frac{(z-u)^4}{\rho^4}
%  +2\frac{(z-u)^2}{\rho^2}
% \right).
%\]
%This map belongs to \(\mathcal N(F)\), since it is obtained from
%\(F\) by positive real affine changes of source and target.  Its
%critical points are
%\[
% u,\qquad u+i\rho,\qquad u-i\rho,
%\]
%and their critical values are
%\[
% f_{u,\rho,v_-,v_+}(u)=v_+,
% \qquad
% f_{u,\rho,v_-,v_+}(u\pm i\rho)=v_-.
%\]
%Conversely, \cref{thm:natural-uniqueness} shows that every member of
%\(\mathcal N(F)\) with ordered singular values \(v_-<v_+\) is of
%this form for a unique pair \((u,\rho)\).
%
%Now impose the one-critical-point normalization (N2) from
%\cref{def:source-normalization-new}: fix the real critical point at
%\(0\) and fix the coefficient of the first nonconstant Taylor term
%there to be \(1/2\).  Writing
%\[
% \Delta=v_+-v_->0,
%\]
%these conditions give \(u=0\) and \(\rho=2\sqrt\Delta\). 
 The (N2)
source-normalized representative in $\mathcal N(F)$ with singular values $v_-,v_+\in \R$ is 
\begin{equation}\label{eq:quartic-source-normalized}
 f_{v_-,v_+}(z)
 =v_++\frac{z^2}{2}+\frac{z^4}{16(v_+-v_-)}
\end{equation}
with critical points 
\(c_0=0,
 c_\pm=\pm2i\sqrt{v_+-v_-}\). 
%\end{equation}
%and
%\[
% f_{v_-,v_+}(c_0)=v_+,
% \qquad
% f_{v_-,v_+}(c_\pm)=v_-.
%\]
For every interior parameter \(v_-<v_+\), the three critical points
are distinct and simple.  In fact,  for an arbitrary family $\mathcal N(F_0)$, distinct labelled critical
points cannot merge at an interior parameter of a natural family:
the source equivalence is a homeomorphism, hence is injective, and it
preserves local degrees. In the above specific family one can see this as follows: as \(v_+-v_-\to0\), the nonreal conjugate pair of critical points 
approach the real critical point
\(0\).  This occurs only at the boundary of the singular-value
chamber: the quartic coefficient in
\eqref{eq:quartic-source-normalized} tends to infinity, so there is no
finite locally uniform limit in this source normalization.  

Without fixing the source scale, the same geometric approach can be
hidden from the singular values.  
\begin{equation}\label{eq:quartic-full-family}
 f_{u,\rho,v_-,v_+}(z)
 =v_++4(v_+-v_-)
 F\left(\frac{z-u}{\rho}\right).
\end{equation}  Indeed, with \(v_-<v_+\) fixed,
\(u=0\), and \(\rho\to0\), the critical points
\(u,u\pm i\rho\) approach one another while the singular values remain
exactly \(v_-\) and \(v_+\).  This is an affine source degeneration,
not a merger at any member of the natural family, since \(\rho>0\) at
every interior parameter.
\end{example}

\begin{example} We have implemented the algorithm also in the case of polynomials with non-real critical points. Take for example the following quintic polynomial 
\( f(z)=b+a(z-b)(z^2+1)^2\) with \[
\begin{aligned}
a&=0.00643312448933335978126858432005919643515420305367\ldots,\\
b&=6.91871796814417406738448217377038172418956021166705\ldots.
\end{aligned}
\]
$f$ has four non-degenerate critical points are \( i,-i, c_1, c_2\) where $c_1,c_2$ are real and three singular values: 
\(f(i)=f(-i)=b\) and  \(v_1=f(c_1)\) and \(v_2=f(c_2)\). The figure below shows that this map has the following itineraries \( c_1\longmapsto v_1\longmapsto w_1=f(v_1) \longmapsto c_1\), \(
c_2\longmapsto v_2\longmapsto c_2\) and 
\( i,-i\longmapsto b\longmapsto b\). So the postsingular set consists of six points. 
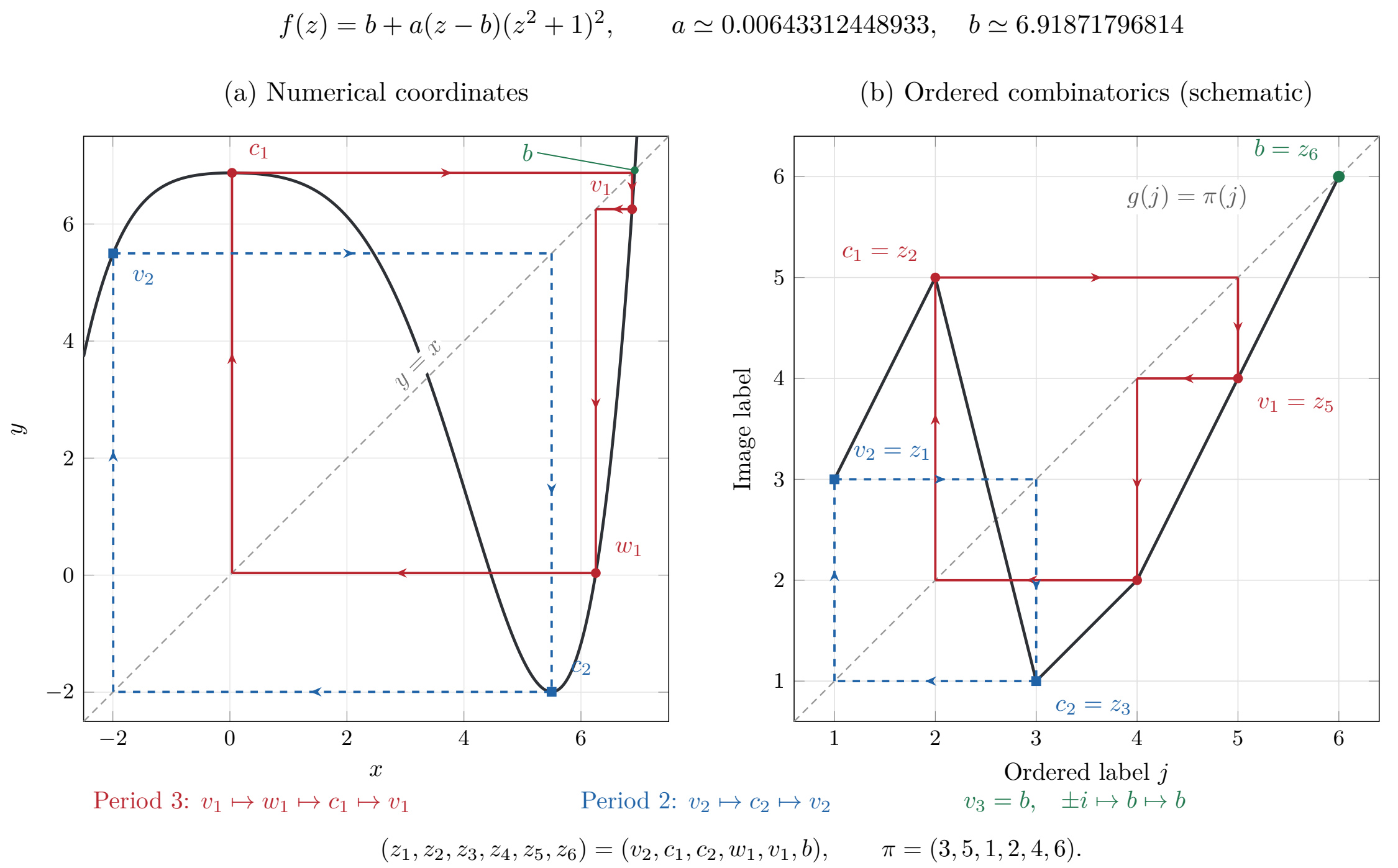
\begin{figure}[htbp]
  \centering
  \resizebox{\textwidth}{!}{%
    % ASCII-only TikZ/PGFPlots picture, suitable for \input in an article.
% Required in the preamble:
%   \usepackage{amsmath,amssymb,pgfplots}
%   \pgfplotsset{compat=1.18}
%   \usetikzlibrary{arrows.meta,decorations.markings}
%
% The left panel uses the actual numerical coordinates of the quintic.
% The right panel is the piecewise-linear interpolation of the ordered
% combinatorics; its equally spaced labels are NOT metric coordinates.
% The non-real points +/-i occur only in the annotation below the panels.
\begingroup
\definecolor{qSpiderRed}{RGB}{186,38,48}
\definecolor{qSpiderBlue}{RGB}{32,99,167}
\definecolor{qSpiderGreen}{RGB}{30,119,77}
\definecolor{qSpiderGraph}{RGB}{45,48,53}

% Parameters and orbit points (rounded only for plotting).
\def\qA{0.0064331244893333597812685843200592}
\def\qB{6.91871796814417406738448217377038}
\def\qCOne{0.03637288506609570666587584542024}
\def\qCTwo{5.49860148944924354724170989359607}
\def\qVOne{6.87432575751941686518879201982546}
\def\qVTwo{-1.99416758114323151195961764860803}
\def\qWOne{6.25369330673885611852954057577580}

\begin{tikzpicture}[
  font=\small,
  qspider/step/.style={
    line width=0.95pt,
    postaction={decorate},
    decoration={markings,
      mark=at position 0.55 with
        {\arrow{Stealth[length=4.5pt,width=3.8pt]}}}
  },
  qspider/red/.style={qspider/step,draw=qSpiderRed},
  qspider/blue/.style={qspider/step,draw=qSpiderBlue,dashed},
  qspider/label/.style={font=\small,inner sep=1.5pt,fill=white}
]

\node[font=\normalsize] at (9.3cm,10.0cm)
  {$f(z)=b+a(z-b)(z^2+1)^2,\qquad
    a\simeq0.00643312448933,\quad b\simeq6.91871796814$};

% --------------------------------------------------------------------
% (a) Actual real polynomial and the two nontrivial singular-value spiders.
% --------------------------------------------------------------------
\begin{axis}[
  name=actual,
  at={(0,0)},anchor=south west,
  width=8.4cm,height=8.4cm,scale only axis,
  xmin=-2.5,xmax=7.5,ymin=-2.5,ymax=7.5,
  axis lines=box,
  xlabel={$x$},ylabel={$y$},
  xtick={-2,0,2,4,6},ytick={-2,0,2,4,6},
  tick label style={font=\footnotesize},
  grid=major,grid style={black!9},
  title={(a) Numerical coordinates},
  title style={font=\normalsize,yshift=2pt},
  clip=true
]
  \addplot[black!40,densely dashed,line width=0.6pt,
           domain=-2.5:7.5,samples=2] {x};
  \addplot[qSpiderGraph,line width=1.2pt,domain=-2.5:7.12,samples=650]
    {\qB+\qA*(x-\qB)*(x*x+1)^2};

  % Red: v1 -> w1 -> c1 -> v1.
  % Vertical segments go to the graph; horizontal segments go to y=x.
  \coordinate (QR0) at (axis cs:\qVOne,\qVOne);
  \coordinate (QR1) at (axis cs:\qVOne,\qWOne);
  \coordinate (QR2) at (axis cs:\qWOne,\qWOne);
  \coordinate (QR3) at (axis cs:\qWOne,\qCOne);
  \coordinate (QR4) at (axis cs:\qCOne,\qCOne);
  \coordinate (QR5) at (axis cs:\qCOne,\qVOne);
  \draw[qspider/red] (QR0)--(QR1);
  \draw[qspider/red] (QR1)--(QR2);
  \draw[qspider/red] (QR2)--(QR3);
  \draw[qspider/red] (QR3)--(QR4);
  \draw[qspider/red] (QR4)--(QR5);
  \draw[qspider/red] (QR5)--(QR0);

  % Blue: v2 -> c2 -> v2.
  \coordinate (QB0) at (axis cs:\qVTwo,\qVTwo);
  \coordinate (QB1) at (axis cs:\qVTwo,\qCTwo);
  \coordinate (QB2) at (axis cs:\qCTwo,\qCTwo);
  \coordinate (QB3) at (axis cs:\qCTwo,\qVTwo);
  \draw[qspider/blue] (QB0)--(QB1);
  \draw[qspider/blue] (QB1)--(QB2);
  \draw[qspider/blue] (QB2)--(QB3);
  \draw[qspider/blue] (QB3)--(QB0);

  % Only points on the graph are marked, not crossings of the spiders.
  \addplot[only marks,mark=*,mark size=1.8pt,qSpiderRed]
    coordinates {(\qVOne,\qWOne) (\qWOne,\qCOne) (\qCOne,\qVOne)};
  \addplot[only marks,mark=square*,mark size=1.8pt,qSpiderBlue]
    coordinates {(\qVTwo,\qCTwo) (\qCTwo,\qVTwo)};

  % b is fixed: its spider is a single diagonal point, not a loop.
  \addplot[only marks,mark=*,mark size=1.4pt,qSpiderGreen]
    coordinates {(\qB,\qB)};
  \draw[qSpiderGreen,line width=0.5pt]
    (axis cs:\qB,\qB)--(axis cs:5.25,7.22)
    node[anchor=east,qspider/label,text=qSpiderGreen] {$b$};

  \node[anchor=south west,text=qSpiderRed,xshift=3pt,yshift=2pt]
    at (axis cs:\qCOne,\qVOne) {$c_1$};
  \node[anchor=south east,text=qSpiderRed,xshift=-4pt,yshift=2pt]
    at (axis cs:\qVOne,\qWOne) {$v_1$};
  \node[anchor=south west,text=qSpiderRed,xshift=4pt,yshift=3pt]
    at (axis cs:\qWOne,\qCOne) {$w_1$};
  \node[anchor=north west,text=qSpiderBlue,xshift=4pt,yshift=-3pt]
    at (axis cs:\qVTwo,\qCTwo) {$v_2$};
  \node[anchor=south west,text=qSpiderBlue,xshift=4pt,yshift=3pt]
    at (axis cs:\qCTwo,\qVTwo) {$c_2$};
  \node[rotate=45,text=black!55,fill=white,inner sep=1pt]
    at (axis cs:3.25,3.58) {$y=x$};
\end{axis}

% --------------------------------------------------------------------
% (b) Schematic: g is the linear interpolation of pi on the integer labels.
% --------------------------------------------------------------------
\begin{axis}[
  name=ordered,
  at={(10.2cm,0)},anchor=south west,
  width=8.4cm,height=8.4cm,scale only axis,
  xmin=0.6,xmax=6.4,ymin=0.6,ymax=6.4,
  axis lines=box,
  xlabel={Ordered label $j$},ylabel={Image label},
  xtick={1,2,3,4,5,6},ytick={1,2,3,4,5,6},
  tick label style={font=\footnotesize},
  grid=major,grid style={black!12},
  title={(b) Ordered combinatorics (schematic)},
  title style={font=\normalsize,yshift=2pt},
  clip=true
]
  \addplot[black!40,densely dashed,line width=0.6pt,
           domain=0.6:6.4,samples=2] {x};
  \addplot[qSpiderGraph,line width=1.2pt]
    coordinates {(1,3) (2,5) (3,1) (4,2) (5,4) (6,6)};

  % Red cycle 5 -> 4 -> 2 -> 5.
  \coordinate (SR0) at (axis cs:5,5);
  \coordinate (SR1) at (axis cs:5,4);
  \coordinate (SR2) at (axis cs:4,4);
  \coordinate (SR3) at (axis cs:4,2);
  \coordinate (SR4) at (axis cs:2,2);
  \coordinate (SR5) at (axis cs:2,5);
  \draw[qspider/red] (SR0)--(SR1);
  \draw[qspider/red] (SR1)--(SR2);
  \draw[qspider/red] (SR2)--(SR3);
  \draw[qspider/red] (SR3)--(SR4);
  \draw[qspider/red] (SR4)--(SR5);
  \draw[qspider/red] (SR5)--(SR0);

  % Blue cycle 1 -> 3 -> 1.
  \coordinate (SB0) at (axis cs:1,1);
  \coordinate (SB1) at (axis cs:1,3);
  \coordinate (SB2) at (axis cs:3,3);
  \coordinate (SB3) at (axis cs:3,1);
  \draw[qspider/blue] (SB0)--(SB1);
  \draw[qspider/blue] (SB1)--(SB2);
  \draw[qspider/blue] (SB2)--(SB3);
  \draw[qspider/blue] (SB3)--(SB0);

  \addplot[only marks,mark=*,mark size=1.8pt,qSpiderRed]
    coordinates {(5,4) (4,2) (2,5)};
  \addplot[only marks,mark=square*,mark size=1.8pt,qSpiderBlue]
    coordinates {(1,3) (3,1)};
  \addplot[only marks,mark=*,mark size=2.2pt,qSpiderGreen]
    coordinates {(6,6)};
  \node[anchor=south east,text=qSpiderGreen,xshift=-4pt,yshift=3pt]
    at (axis cs:6,6) {$b=z_6$};
  \node[anchor=south east,text=qSpiderRed,xshift=-3pt,yshift=3pt]
    at (axis cs:2,5) {$c_1=z_2$};
  \node[anchor=north west,text=qSpiderRed,xshift=4pt,yshift=-4pt]
    at (axis cs:5,4) {$v_1=z_5$};
  \node[anchor=south west,text=qSpiderBlue,xshift=4pt,yshift=4pt]
    at (axis cs:1,3) {$v_2=z_1$};
  \node[anchor=north west,text=qSpiderBlue,xshift=4pt,yshift=-4pt]
    at (axis cs:3,1) {$c_2=z_3$};
  \node[text=black!65,fill=white,inner sep=2pt]
    at (axis cs:4.5,5.8) {$g(j)=\pi(j)$};
\end{axis}

% The orbits are written starting at the singular values.
\node[anchor=west,text=qSpiderRed] at (0cm,-1.15cm)
  {Period $3$: $v_1\mapsto w_1\mapsto c_1\mapsto v_1$};
\node[anchor=west,text=qSpiderBlue] at (7.0cm,-1.15cm)
  {Period $2$: $v_2\mapsto c_2\mapsto v_2$};
\node[anchor=west,text=qSpiderGreen] at (12.5cm,-1.15cm)
  {$v_3=b$,\quad $\pm i\mapsto b\mapsto b$};
\node at (9.3cm,-1.85cm)
  {$(z_1,z_2,z_3,z_4,z_5,z_6)=(v_2,c_1,c_2,w_1,v_1,b),
    \qquad \pi=(3,5,1,2,4,6)$.};
\end{tikzpicture}
\endgroup%
  }
  \caption{The spider diagram of the three real singular values of the quintic.
  The red and blue paths have periods three and two, respectively. $v_1$ is almost a fixed point, but actually has period three. 
  The green point is the common critical value of the non-real
  critical pair and is fixed. The right-hand panel shows the
  ordered combinatorics with equally spaced labels.}
  \label{fig:mixed-quintic-spiders}
\end{figure}
\end{example} 

\subsection{Uniqueness modulo source affine maps}
\begin{theorem}[Uniqueness for prescribed singular values]
\label{thm:natural-uniqueness}
Suppose
\[
 f_j=\phi_j\circ f_0\circ\psi_j\in\mathcal N(f_0),
 \qquad j=1,2,
\]
and \(\phi_1(v_{k,0})=\phi_2(v_{k,0})\) for every \(k\).
Then
\[
 f_1=f_2\circ L
\]
for a positive real affine map \(L\). If both representatives satisfy
the same source normalisation of
\cref{def:source-normalization-new}, then \(f_1=f_2\).
\end{theorem}

\begin{proof}
The real-symmetric increasing maps \(\phi_1\) and \(\phi_2\) are
isotopic through real-symmetric plane homeomorphisms \(\phi_t\) whose
values on \(S(f_0)\) remain fixed. One way to see this is to interpolate
their increasing restrictions to \(\R\), relative to the finitely
many singular points, extend the resulting boundary isotopy to the
upper half-plane, and reflect it to the lower half-plane. The residual
homeomorphism fixing the real boundary pointwise is isotopic to the
identity by the Alexander isotopy on a closed disk.

Lift this isotopy using the covering in \eqref{eq:covering}, starting
with \(\psi_1\), to obtain
\[
 f_1=\phi_t\circ f_0\circ\widetilde\psi_t.
\]
The lifts extend across the discrete set of finite preimages of the
singular values by the local power model. The inverse lifted isotopy
shows that each extension is a plane homeomorphism. The lifts are real
symmetric and increasing on \(\R\), by uniqueness and continuity of
lifting. Thus
\[
 f_1=\phi_2\circ f_0\circ\widetilde\psi_2
     =f_2\circ L,
 \qquad L=\psi_2^{-1}\circ\widetilde\psi_2.
\]
Away from the discrete preimage of the singular set, a local inverse
of \(f_2\) shows that \(L\) is holomorphic. Its isolated
singularities are removable because \(L\) is a homeomorphism. Thus
\(L\) is an entire conformal automorphism of \(\C\). Its reality
and monotonicity give \(L(z)=az+b\), with \(a>0\) and \(b\in\R\).

A critical source anchor point over a singular value stays fixed during the lifted
isotopy: its lift is a continuous path in a discrete fibre. Under
normalisation (N1), \(L\) therefore fixes both \(c_1\) and \(c_2\),
so \(L=\mathrm{id}\). Under (N2), it fixes \(c\) and has the form
\(L(c+z)=c+az\). Comparing the first nonconstant Taylor coefficients
in \(f_1=f_2\circ L\) gives \(A_\ell=A_\ell a^\ell\). Since
\(a>0\), again \(L=\mathrm{id}\).
\end{proof}

\begin{remark}[Why arbitrary fixed source points do not normalise the family]
\label{rem:arbitrary-source-new}
The observation behind the following example is due to Nikolai
Prochorov. Let
\[
 f(z)=ze^{z+1}+1,
 \qquad f_\kappa(z)=f(\kappa z),\qquad\kappa>0.
\]
All these maps have the same singular values \(0,1\), but they can
be written
\(f_\kappa=\phi_\kappa\circ f\circ\psi_\kappa\) with
\(\phi_\kappa|\{0,1\}=\mathrm{id}\) and
\(\psi_\kappa(0)=0\), \(\psi_\kappa(1)=1\).
Indeed, choose \(\phi_\kappa\) increasing and real symmetric so that
\[
 \phi_\kappa(e^2+1)=\kappa e^{\kappa+1}+1.
\]
An isotopy fixing \(0,1\) lifts under \(f\) to a map
\(\chi_\kappa\) satisfying
\(f\circ\chi_\kappa=\phi_\kappa\circ f\), with
\(\chi_\kappa(0)=0\), \(\chi_\kappa(1)=\kappa\).
Set \(\psi_\kappa(z)=\chi_\kappa^{-1}(\kappa z)\).
These representatives are distinct and
\(f_\kappa(1)\to\infty\) as \(\kappa\to\infty\).
The point \(1\) is not a preimage of a singular value of \(f\),
so fixing it in a chosen quasiconformal natural equivalence does not eliminate
the holomorphic affine source freedom.

Other target coordinates may be used as additional normalization
data if both their values and their chosen source inverse-branch
labels are included in the deformation data. The critical source anchor points
in \cref{def:source-normalization-new} make this bookkeeping automatic.
\end{remark}

\subsection{Existence and continuity of the singular-value parameterisation}
\begin{theorem}[Singular-value coordinates]
\label{thm:natural-parameters}
Fix either source normalisation in
\cref{def:source-normalization-new}. For every
\[
 v=(v_1,\ldots,v_q)\in\mathcal V_q
 :=\{v\in\R^q:v_1<\cdots<v_q\},
\]
there is exactly one \(f_v\in\mathcal N_0(f_0)\) whose ordered
singular values are \(v_1,\ldots,v_q\). The map
\(v\mapsto f_v\) is continuous for locally uniform convergence on
\(\C\). Locally in \(\mathcal V_q\), the natural equivalences and
all fixed finite collections of labelled real critical points and
labelled monotonicity branches (identified by the natural equivalence)
can be chosen continuously.

In particular, the conclusion holds on
\(\{0\le v_1<\cdots<v_q\le1\}\).
\end{theorem}

\begin{proof}
Let \(h_v:\R\to\R\) be the increasing piecewise affine map taking
\(v_{j,0}\) to \(v_j\), affine on each interval between consecutive
reference singular values, and of slope one on the two exterior
half-rays. For \(q=1\), take the corresponding translation. The map
\[
 \phi_v(x+iy)=h_v(x)+iy
\]
is a real-symmetric quasiconformal homeomorphism. Its Beltrami
coefficient depends continuously on \(v\) in the \(L^\infty\) norm
locally in \(\mathcal V_q\): the breakpoints in the source are fixed,
and the finitely many positive slopes vary continuously.

Pull its Beltrami coefficient back by \(f_0\): off the discrete
critical set, put
\[
 \mu_v(z)=\mu_{\phi_v}(f_0(z))\,
          \frac{\overline{f_0'(z)}}{f_0'(z)}.
\]
The values on the critical set are irrelevant. The measurable Riemann
mapping theorem supplies a real-symmetric quasiconformal map \(H_v\)
with coefficient \(\mu_v\). Then
\[
 \widetilde f_v=\phi_v\circ f_0\circ H_v^{-1}
\]
has zero Beltrami coefficient and is entire. Its singular values are
exactly \(v_1,\ldots,v_q\), because the natural equivalence preserves
the covering characterisation of singular values. This is the standard
straightening construction; see \cite{LV73,EL92}.

For (N1), normalise \(H_v\) by fixing \(c_1,c_2,\infty\). The
symmetry of \(\mu_v\) and uniqueness of the normalised solution make
\(H_v\) real symmetric, and its restriction to \(\R\) is increasing.
Then \(\widetilde f_v\) already has the required source
normalisation. For (N2), first fix \(c,c+1,\infty\), and write
\[
 \widetilde f_v(c+z)-\widetilde f_v(c)
     =b_v z^\ell+O(z^{\ell+1}).
\]
The coefficient \(b_v\) is nonzero and has the sign of \(A_\ell\).
Precompose by
\[
 L_v(c+z)=c+\left(\frac{A_\ell}{b_v}\right)^{1/\ell}z
\]
to obtain \(f_v=\widetilde f_v\circ L_v\). In either case,
\(f_v\in\mathcal N_0(f_0)\).

The normalised measurable Riemann mapping theorem depends continuously
on the Beltrami coefficient when the coefficients converge in
\(L^\infty\) and remain uniformly bounded away from one. These
conditions hold locally in \(\mathcal V_q\). Thus \(H_v\),
\(H_v^{-1}\), and \(\widetilde f_v\) vary locally uniformly. Cauchy's
formula makes \(b_v\) continuous, so the final affine correction in
case (N2) is continuous as well. Applying $H_v$ and the final affine correction to a fixed finite
list of reference critical points or monotonicity branches proves
the last continuity assertion. Uniqueness follows from
\cref{thm:natural-uniqueness}.
\end{proof}

The \emph{chamber of distinct singular values} is the open region
$\mathcal V_q=\{v_1<\cdots<v_q\}$, where the singular-value order and
labels cannot change. Its boundary includes configurations where
consecutive singular values collide. No continuous extension of
$v\mapsto f_v$ across those collisions is asserted. The quantitative
bound \cref{ana:scaled}, combined with \cref{gap:separation}, handles
the relevant configurations without such an extension.

\section{Finite markings and expansive combinatorics}
\label{sec:combinatorics-new}

The singular-value coordinates from \cref{thm:natural-parameters}
describe the holomorphic maps. A \emph{finite marking} consists of a
finite ordered forward-invariant set $X$ containing the singular
values, together with the critical-point and inverse-branch labels
of its elements. Its order and the map $g|_X$ specify the desired
finite dynamics; the labels specify which real inverse is selected.
I first arrange this marking so that its combinatorics records every
real critical point that is relevant to the pullback.

\subsection{Topological models and the absence of Levy cycles}

\begin{definition}[Labelled real topological model]
\label{def:model-new}
A real topological model for $f_0$ consists of a continuous map
$g:\R\to\R$ and increasing homeomorphisms
$\phi_g,\psi_g:\R\to\R$ such that
\[
 g=\phi_g\circ f_0|_{\R}\circ\psi_g.
\]
The model includes the distinguished sets
\[
 S(g)=\phi_g(S(f_0)),\qquad
 \Crit(g)=\psi_g^{-1}(\Crit_{\R}(f_0)),\qquad
 P(g)=\bigcup_{n\ge0}g^n(S(g)).
\]
Each point of $\Crit(g)$ carries the local degree of its corresponding
critical point of $f_0$. The distinguished singular values include those
arising from nonreal critical points and nonreal asymptotic tracts.
Thus $S(g)$ is part of the model data, even when some of its points have
no special status for the real graph alone. Likewise, stationary real
critical points that are not turning points remain distinguished.
\end{definition}

I assume that a finite forward-invariant set $X\supset P(g)$ has been
specified. All points of $X$ are finite. The real graph can have
infinitely many critical points, but only finitely many lie in any
compact real interval.

\begin{definition}[Absence of real Levy cycles]
\label{def:levy-new}
The model has \emph{no Levy cycles on $X$} if, for every pair of distinct
points $a,b\in X$, there is $n\ge0$ for which
\begin{equation}\label{eq:model-levy}
 g^n([a,b])\cap\Crit(g)\ne\varnothing.
\end{equation}
Here and below $[a,b]$ denotes the interval with endpoints $a,b$,
irrespective of their order. The notation $g^n([a,b])$ means the
\emph{full image of the interval} under the $n$-th iterate, not merely
$[g^n(a),g^n(b)]$. A critical point in the boundary of the
interval counts, and $n=0$ is permitted.
\end{definition}

This is the interval formulation used in this paper. In particular, the
condition concerns the distinguished critical points of the model, and
not merely the turning points of a piecewise linear drawing.

\begin{proposition}[Necessity of the real Levy condition]
\label{prop:levy-necessary}
Let $f$ be a real entire map with finite nonempty postsingular set
$P(f)\subset\R$. Then, for every two distinct points $a,b\in P(f)$,
there is $j\ge0$ such that
\[
 f^j([a,b])\cap\Crit_{\R}(f)\ne\varnothing.
\]
In particular, the marked dynamics of a postsingularly finite real
entire map necessarily satisfies the no-Levy-cycle condition of
\cref{def:levy-new}.
\end{proposition}

\begin{proof}
Suppose that the conclusion fails for $a<b$ and put
\[
 I_j=f^j([a,b]).
\]
No interval $I_j$ contains a real critical point. Hence $f$ is strictly
monotone on every $I_j$ and maps $I_j$ diffeomorphically onto
$I_{j+1}$. The two endpoints of $I_j$ belong to the finite set $P(f)$.
Consequently two unoriented endpoint pairs eventually repeat. After
replacing the original interval by a forward image, there are a closed
interval $J=[\alpha,\beta]$, with $\alpha,\beta\in P(f)$, and an integer
$n\ge1$ such that
\[
 F:=f^n:J\longrightarrow J
\]
is a diffeomorphism. If $F$ reverses orientation, replace $F$ by $F^2$.
Thus $F$ may be assumed increasing, and then it fixes both endpoints of
$J$.

The analytic function $F(x)-x$ is not identically zero on $J$; otherwise
the entire map $F$ would be the identity. Between any two consecutive
fixed points of $F$ in $J$, the function $F(x)-x$ has constant sign.
If it is positive, its right endpoint $p$ satisfies $F'(p)\le1$; if it
is negative, its left endpoint does. Hence $F$ has a fixed point
$p\in J$ with
\[
 0<F'(p)\le1,
\]
because $F$ is a diffeomorphism on $J$. This gives an attracting or
parabolic periodic point of $f$ which is not superattracting.

I use here the standard finite-type consequence that a
postsingularly finite entire map is subhyperbolic: its Fatou set is
empty or a finite union of superattracting basins, and it has no
parabolic cycles; see, for example, \cite{EL92} or \cite{MB}. Thus a non-repelling
periodic point of a postsingularly finite entire map must be
superattracting. That is impossible here because $F'(p)>0$.
This contradiction proves the proposition.
\end{proof}

\begin{remark}[Inverse-branch interpretation]
\label{rem:levy-schwarz}
The preceding proof also has the Schwarz--Pick interpretation suggested
by the real pullback construction. After passing to an even iterate if
necessary, $F=f^n$ preserves orientation on $J$. Since
$S(F)\subset P(f)\subset\R$, the inverse branch
$h=(F|_J)^{-1}$ extends holomorphically to a self-map of the upper
half-plane and fixes the two endpoints of $J$. The Julia--Wolff form of
Schwarz--Pick says that such a branch cannot make both fixed endpoints
strictly repelling for $F$; equivalently, $F$ has a non-repelling fixed
point on $J$. The elementary real argument in the proof gives this last
conclusion without invoking the boundary form of Schwarz--Pick.
\end{remark}

\begin{lemma}[Adding the relevant critical points]
\label{lem:critical-saturation}
Let $X\supset S(g)$ be finite and forward invariant, and suppose that
$X$ contains at least two points. Put $J_g=\operatorname{conv}X$ and
\[
 X^+=X\cup\bigl(\Crit(g)\cap J_g\bigr).
\]
Then $J_g$ is forward invariant, $X^+$ is finite and forward invariant,
and $g$ is strictly monotone on every interval between consecutive
points of $X^+$. If $g$ has no Levy cycles on $X$, it has no Levy cycles
on $X^+$.
\end{lemma}

\begin{proof}
The images of the endpoints of $J_g$ belong to $X$, and every real local
extremum has its value in $S(g)\subset X$. If $g(J_g)$ left $J_g$, an
interior maximum or minimum would therefore have a value outside
$J_g$, a contradiction. This proves $g(J_g)\subset J_g$.

The set $\Crit(g)\cap J_g$ is finite: its image under $\psi_g$ lies in a
compact real interval, and the zeros of the nonzero entire function
$f_0'$ are isolated. Every newly added point maps into $S(g)\subset X$,
so $X^+$ is forward invariant. There is no critical point in the
interior of a consecutive marked interval; the real analytic model
therefore has derivative of one sign there and is strictly monotone.
The same is true of $g$ by its increasing topological equivalence.

For a pair of points in $X^+$, either both belonged to $X$, in which
case the given condition applies, or at least one is a newly added
critical point, in which case \eqref{eq:model-levy} holds already for
$n=0$.
\end{proof}

I henceforth replace $X$ by $X^+$ and enumerate it as
\[
 X=\{u_1<\cdots<u_s\}.
\]
This replacement changes neither the realisation problem nor the
prescribed dynamics on the original marked set. Define
\begin{equation}\label{eq:pi-definition}
 g(u_i)=u_{\pi(i)},\qquad
 \pi:\{1,\ldots,s\}\longrightarrow\{1,\ldots,s\},
\end{equation}
and let
\[
 \cC=\{j:u_j\in\Crit(g)\},\qquad
 S(g)=\{u_{i_1},\ldots,u_{i_q}\},\quad i_1<\cdots<i_q.
\]
Every source label records either its particular critical point or its
particular component of $\R\setminus\Crit(g)$. These labels are
transported to each holomorphic map in the natural class. Because the
marking is finite, only finitely many real branches are selected. This
is a restriction on the branches used by the Thurston algorithm defined below, not on the
total number of branches of $f_0$.

\begin{definition}[Expansive finite combinatorics]
\label{def:pi-expansive}
Let $G:[1,s]\to[1,s]$ be the piecewise linear interpolation of
$G(i)=\pi(i)$. The pair $(\pi,\cC)$ is \emph{expansive} if
\begin{equation}\label{eq:expansive}
 \text{for every }i<s\text{ there is }n\ge0
 \text{ such that }G^n([i,i+1])\cap\cC\ne\varnothing.
\end{equation}
The points of $\cC$ are actual finite critical-point labels;
formal symbols $-\infty$ and $+\infty$ used in drawings of the
real-line ends are not critical labels and do not count. Write $N=N(\pi,\cC)$ for a common upper bound on the
first hitting times in \eqref{eq:expansive}.
\end{definition}

\begin{lemma}[Equivalence of the two finite conditions]
\label{lem:expansive-levy}
For the marking provided by \cref{lem:critical-saturation},
\begin{equation}\label{eq:noflat}
 \pi(i)\ne\pi(i+1)\qquad(1\le i<s),
\end{equation}
and \eqref{eq:model-levy} is equivalent to \eqref{eq:expansive}.
\end{lemma}

\begin{proof}
Strict monotonicity on $[u_i,u_{i+1}]$ gives
\[
 g([u_i,u_{i+1}])=[u_{\pi(i)},u_{\pi(i+1)}]
\]
and proves \eqref{eq:noflat}. Let $h:[1,s]\to J_g$ be increasing and
piecewise affine with $h(i)=u_i$. If $I\subset[1,s]$ is an interval
with integer endpoints, it is a union of unit edges. Taking the union
of their images shows
\[
 g(h(I))=h(G(I)).
\]
Both resulting intervals again have marked endpoints. Induction gives
$g^n([u_i,u_{i+1}])=h(G^n([i,i+1]))$ for every $n$. All critical
points in $J_g$ occur among the marked points, so the two critical
intersection conditions are equivalent for adjacent pairs. Finally,
every nondegenerate interval with endpoints in $X$ contains an adjacent
marked interval, which proves the assertion for all pairs.
\end{proof}

The definition is a finite combinatorial test: direct an arrow from
each unit edge to every unit edge contained in its image under $G$.
Mark the edges that have an endpoint in $\cC$. Expansiveness means
that each edge reaches a marked edge. Since the number of edges is
finite, $N$ is finite. This formulation is closely related to the
critical-edge condition for the real polynomial Thurston algorithm
in \cite[Section~2]{BMS}.

\begin{remark}[Finite coordinates and real ends]
\label{rem:finite-ends-new}
By the two \emph{real ends} the directions $x\to-\infty$
and $x\to+\infty$ on the real line are meant, not points in the domain of $g$.
Formal symbols $-\infty$ and $+\infty$ may depict these directions. They are not coordinates of the finite marking, and they are
not included in its affine normalisation. A \emph{terminal monotone branch} is a maximal open interval of
monotonicity extending to $-\infty$ or $+\infty$ and having no
further real critical point in that direction. On such a branch a
finite limiting value is a singular value, but is not attained
on that open branch. It therefore cannot be assigned to a finite
regular inverse mark there. If infinitely many critical points
continue indefinitely as $x\to-\infty$ or $x\to+\infty$,
a limiting value in that direction need not
exist at all.
\end{remark}

\section{The Thurston map}
\label{sec:thurston}

Fix an expansive finite marking as in the preceding section. The Thurston map pulls the coordinates $x_i$ (the \emph{target
marks}) back along their prescribed branches, then applies a positive
affine normalisation. The corresponding unscaled source preimages
are denoted $\widehat y_i$; after normalising their extreme values
to $0$ and $1$ are denoted by $y_i$.

Set
\[
 W=\{x=(x_1,\ldots,x_s):0=x_1<x_2<\cdots<x_s=1\},
 \qquad
 v(x)=(x_{i_1},\ldots,x_{i_q}).
\]
For $x\in W$, let $f_x=f_{v(x)}$ be the representative given by
\cref{thm:natural-parameters}, with the fixed source normalisation.
In particular,
\begin{equation}\label{eq:singular-marking}
 S(f_x)=\{x_{i_1},\ldots,x_{i_q}\}\subset\{x_1,\ldots,x_s\}.
\end{equation}
The branch and critical labels of the model determine corresponding
labels for $f_x$.

\begin{definition}[Thurston operator]
\label{def:thurston-new}
Take $x\in W$. For each $i$, choose the finite real point $\widehat y_i$ on the branch,
or at the critical point, corresponding to $u_i$, subject to
\begin{equation}\label{eq:marking}
 f_x(\widehat y_i)=x_{\pi(i)}.
\end{equation}
For a critical label this includes
\begin{equation}\label{eq:critical-marking}
 f_x'(\widehat y_j)=0\qquad(j\in\cC).
\end{equation}
Write
\begin{equation}\label{eq:affine}
 \begin{split}
 L_x&=\widehat y_s-\widehat y_1,\qquad
 M_x(z)=\widehat y_1+L_xz,\\
 y_i&=M_x^{-1}(\widehat y_i)
      =\frac{\widehat y_i-\widehat y_1}{L_x},\qquad
 F_x=f_x\circ M_x.
 \end{split}
\end{equation}
The Thurston operator is
\[
 T:W\longrightarrow W,\qquad T(x)=y=(y_1,\ldots,y_s).
\]
\end{definition}

\begin{lemma}[Existence, ordering, and continuity of the pullback]
\label{lem:thurston-continuous}
The points $\widehat y_i$ exist, are uniquely determined by their
labels, and satisfy $\widehat y_1<\cdots<\widehat y_s$. The map
$T:W\to W$ is continuous.
\end{lemma}

\begin{proof}
For a regular source label, the target $g(u_i)=u_{\pi(i)}$ lies
strictly inside the image of the selected model branch. The finite
endpoints of that image are singular values. The position of the
target relative to those endpoints is therefore specified by its
index in the ordered set of target coordinates. The same strict inequalities hold for
$x_{\pi(i)}$ and the corresponding branch image of $f_x$. Strict
monotonicity supplies exactly one inverse point on that branch.
Critical source labels have their prescribed singular values and
their uniquely labelled critical points.

Within an increasing branch the target indices increase with the
source indices; within a decreasing branch they decrease. The
corresponding inverse points consequently have the required source
order. Across branches and critical points their order follows from
the fixed labels. This proves $L_x>0$ and $T(x)\in W$.

For continuity, work near a fixed parameter $x^0\in W$. The local
quasiconformal construction in \cref{thm:natural-parameters} makes
the labelled critical points and branch endpoints vary continuously,
and $f_x$ varies locally uniformly. Every regular inverse point is a
simple solution of \eqref{eq:marking}; the local inverse theorem, or
Rouch\'e's theorem on a small disk, gives its continuous dependence on
$x$. A critical source label follows its continuously transported
critical point. Thus all the finitely many $\widehat y_i$ are
continuous, as are $L_x$, $M_x$, and $T$.
\end{proof}

\begin{lemma}[Independence of the preliminary source normalisation]
\label{lem:thurston-anchor-independent}
Changing the preliminary positive affine source normalisation of
$f_x$ does not change $T(x)$ or $F_x$.
\end{lemma}

\begin{proof}
If $\widetilde f_x=f_x\circ A_x$ with $A_x$ positive affine, then
$\widetilde{\widehat y}_i=A_x^{-1}(\widehat y_i)$ and
$\widetilde M_x=A_x^{-1}\circ M_x$. Hence
\[
 \widetilde M_x^{-1}(\widetilde{\widehat y}_i)=y_i,
 \qquad
 \widetilde f_x\circ\widetilde M_x=f_x\circ M_x.
\]
The uniqueness up to positive affine precomposition in
\cref{thm:natural-uniqueness} applies to the different source
normalisations under consideration.
\end{proof}

\begin{lemma}[A fixed point realises the required dynamics on the finite set]
\label{lem:thurston-realisation}
If $T(x)=x$, then $F_x\in\mathcal N(f_0)$ satisfies
\[
 F_x(x_i)=x_{\pi(i)}\qquad(1\le i\le s).
\]
Its singular set is $\{x_{i_1},\ldots,x_{i_q}\}$, and the increasing
bijection $x_i\mapsto u_i$ conjugates the marked dynamics to that of
$g$. In particular, every singular value of $F_x$ has a finite forward
orbit with the prescribed combinatorics.
\end{lemma}

\begin{proof}
By construction $F_x(y_i)=x_{\pi(i)}$. At a fixed point $y_i=x_i$.
Precomposition by the affine automorphism $M_x$ preserves singular
values and the natural equivalence class. The remaining assertions
follow from the definition of $\pi$.
\end{proof}

\begin{remark}
The realised dynamical map is $F_x=f_x\circ M_x$. The preliminary
representative $f_x$ is normalised at its source anchor points, whereas
$F_x$ is normalised by its two extreme dynamical marks. These two
normalisations have different purposes.
\end{remark}

\subsection{The real Teichm\"uller chamber}
\label{subsec:real-teich-chamber}

The coordinates on $W$ also describe the marking information usually
retained in Teichm\"uller space. Fix a reference configuration $a\in W$
and fix it throughout this subsection. Put
\[
 P_x=\{x_1,\ldots,x_s,\infty\},\qquad
 \mathcal T_a=\mathcal T(\widehat\C,P_a).
\]
A point of $\mathcal T_a$ is represented by a quasiconformal marking
$H:\widehat\C\to\widehat\C$, modulo conformal postcomposition and
isotopy relative to $P_a$. Normalising $H$ to fix $0,1,\infty$ removes
the conformal ambiguity. In general, the images of the punctures do
not determine the isotopy class; see \cite[Section~3.1]{HSS09}.
In the real setting the following observation removes this extra
choice.

\begin{lemma}[Canonical real marking]\label{lem:canonical-real-marking}
For every $x\in W$ there is exactly one point of $\mathcal T_a$
represented by a real-symmetric quasiconformal map $H$ that is
increasing on $\R$ and satisfies $H(a_i)=x_i$ and $H(\infty)=\infty$.
\end{lemma}

\begin{proof}
Let $h_x$ be affine on each $[a_i,a_{i+1}]$, send $a_i$ to $x_i$,
and equal the identity outside $[0,1]$. Then
$H_x(u+iv)=h_x(u)+iv$ is a real-symmetric bi-Lipschitz, hence
quasiconformal, extension.

For uniqueness, the composition of the inverse of one such marking
with another is a real-symmetric homeomorphism $G$ fixing $P_a$
pointwise. Write $g=G|_{\R}$ and set
\[
 k_t(u)=(1-t)u+t g^{-1}(u),\qquad
 K_t(u+iv)=k_t(u)+iv\quad(0\le t\le1).
\]
The maps $K_t\circ G$ give an isotopy relative to $P_a$ from $G$ to
a map fixing $\widehat\R$ pointwise. Its restriction to the closed
upper hemisphere is a disk homeomorphism fixing the boundary.
The Alexander isotopy contracts it to the identity relative to that
boundary; reflection gives the same construction below. Thus $G$
is isotopic to the identity relative to $P_a$. Only a topological
isotopy is needed here.
\end{proof}

Consequently
\[
 j_a:W\longrightarrow\mathcal T_a,\qquad j_a(x)=[H_x],
\]
is canonical for this fixed reference configuration $a$; here
$[H_x]$ means the equivalence class of $H_x$ under conformal
postcomposition and isotopy fixing $P_a$ pointwise. More precisely,
let $p:\mathcal T_a\to\mathcal M_{0,s+1}$ be the usual covering that
forgets the isotopy class and remembers the labelled punctures,
normalised at $0,1,\infty$; see \cite{Gar87}.
The set $W$ is one component of the real locus of this moduli space.
Since $W$ is convex, each component of $p^{-1}(W)$ maps
homeomorphically onto $W$, and $j_a(W)$ is the component containing
$[\mathrm{id}]$. It is also the component through $[\mathrm{id}]$
of the fixed locus of
\[
 \mathfrak r[H]=[\rho\circ H\circ\rho],\qquad \rho(z)=\bar z.
\]
Indeed, the representatives $H_x$ commute with reflection. Conversely,
a connected subset of this fixed locus through $[\mathrm{id}]$
projects into the same real moduli component $W$, and so lies in
this component of $p^{-1}(W)$.

The gaps $d_i=x_{i+1}-x_i>0$, with $\sum_i d_i=1$, identify $W$
with an open real ball of dimension $s-2$, whereas $\mathcal T_a$
has complex dimension $s-2$; for $s=2$ both are a point. In this
precise sense $W$ is the real analogue of the Teichm\"uller ball
relevant to the algorithm. Here ``ball'' describes its topology,
not a ball of finite radius in the Teichm\"uller metric.
Within this chamber the ordered coordinates already determine the
marking class, so the Thurston map requires no separate isotopy
coordinate. Arbitrary complex markings over the same real
configuration can still differ by nontrivial pure mapping classes (isotopy classes of
homeomorphisms fixing each labelled puncture) when $s\ge3$. The prescribed inverse-branch and critical-point labels
remain part of the fixed dynamical data. Finally, the metric defined
in \cref{sec:convergence-new} is exactly the ambient Teichm\"uller
metric restricted to $j_a(W)$.

\subsection{Two examples}

\begin{example}[A cubic polynomial with two critical source anchor points]
\label{ex:thurston-cubic}
Take $f_0(z)=z^3/3-z$, with critical points $-1,1$, and use
\[
 \pi=(1,4,3,2,5),\qquad \cC=\{2,4\},\qquad
 S(f_x)=\{x_2,x_4\}.
\]
The representative whose critical points remain at $-1,1$ is
\[
 f_x(z)=\frac{3(x_4-x_2)}4\left(\frac{z^3}{3}-z\right)
       +\frac{x_4+x_2}{2}.
\]
Thus $f_x(-1)=x_4$ and $f_x(1)=x_2$. The first inverse point lies
to the left of $-1$, the second is $-1$, the third lies between
$-1$ and $1$, the fourth is $1$, and the fifth lies to the right of
$1$. Every unit edge has a critical endpoint, so the diagram is
expansive with $N=0$. The affine normalisation of these five inverse
points gives $T(x)$.
\end{example}

\begin{example}[A linear-exponential map and a half-ray]
\label{ex:thurston-exponential}
Take $f_0(z)=(z-1)e^z+1$, whose real critical point is $0$, critical
value is $0$, and finite asymptotic value is $1$. Use
\[
 \pi=(2,1,2,4),\qquad \cC=\{2\},\qquad
 S(f_x)=\{x_1,x_3\}=\{0,v\},\quad v=x_3.
\]
Fix the critical point at $0$ and the quadratic Taylor coefficient
at $1/2$, as for $f_0$. Then
\[
 f_x(z)=v\bigl(1+(\lambda z-1)e^{\lambda z}\bigr),
 \qquad \lambda=v^{-1/2}.
\]
The targets for the first and third inverse points are both $x_2$,
with $0<x_2<v$; they are taken on $(-\infty,0)$ and $(0,\infty)$,
respectively. The second inverse point is $0$, and the fourth is the
positive preimage of $1$. The last unit edge maps to $[2,4]$ and
therefore meets the critical label after one iterate; the other edges
already have that label as an endpoint. All pullback points are
finite, although the selected branches are half-rays and their
affine span can be large.
\end{example}

\subsection{The mesh and the distance to the boundary}

For the analytic estimates 
\[
 m(x)=\min_{1\le i<s}(x_{i+1}-x_i),\qquad
 \norm{x-y}_\infty=\max_i|x_i-y_i|
\]
is used. 
Let $\delta_W(x)=\dist_\infty(x,\partial W)$, where the distance
and the relative boundary are taken in the affine subspace
$\{x_1=0,x_s=1\}$.

\begin{lemma}[Comparison of the mesh and boundary distance]
\label{lem:mesh-distance}
If $s\ge3$, then
\begin{equation}\label{eq:mesh-distance}
 \frac12m(x)\le\delta_W(x)\le m(x).
\end{equation}
Consequently a sequence approaches $\partial W$ precisely when its
mesh tends to zero, and $\norm{T(x)-x}_\infty/m(x)\to\infty$ as $x\to\partial W$
if and only if $\norm{T(x)-x}_\infty/\delta_W(x)\to\infty$.
\end{lemma}

\begin{proof}
A change of every free coordinate by less than $m(x)/2$ preserves
all strict inequalities. Conversely, a gap of size $m(x)$ can be
closed by moving its free endpoint to the other endpoint; for an
interior gap both endpoints can instead be moved halfway towards
each other. This proves the two inequalities.
\end{proof}

When $s=2$, the space $W$ is a singleton and $T$ has a fixed point
without a boundary argument. It may therefore be assumed that $s\ge3$
when applying the next section.

\section{A fixed point criterion}
\label{sec:fixedpoint}

The following topological criterion isolates what is required from
the analytic estimates. Its simplex case occurs in
\cite[Chapter~II, Section~4, pp.~121--123]{dMvS93}; the general
convex-body statement is explicitly recorded as the de Melo--van
Strien fixed point theorem in \cite[Proposition~8.1]{Smania05}.

\begin{theorem}[The de Melo--van Strien fixed point criterion]
\label{thm:boundary-fixedpoint}
Let $V$ be a nonempty bounded open convex subset of a
finite-dimensional normed space, and let $T:V\to V$ be continuous.
If
\begin{equation}\label{eq:boundary-displacement}
 \frac{\norm{T(x)-x}}{\dist(x,\partial V)}\longrightarrow\infty
 \quad\text{whenever}\quad\dist(x,\partial V)\longrightarrow0,
\end{equation}
then $T$ has a fixed point in $V$.
\end{theorem}

I give the geometric retraction proof. The zero-dimensional case is
immediate; in the proof assume positive dimension.

%\begin{proof}[Geometric proof of \cref{thm:boundary-fixedpoint}]

% Replace Section 5.1, from \subsection{The geometric proof} up to
% (but not including) \subsection{A direct proof on the convex domain}.
% The manuscript already loads tikz, xcolor and placeins.
\subsection{The geometric proof}

\begin{figure}[tbp]
\centering
\begin{tikzpicture}[>=Stealth,scale=1.15,
  every node/.style={font=\small}]
  % The dashed lines are supporting planes at the south poles.
  \begin{scope}[shift={(0,0)}]
    \draw[softgray,thick] (0,0) circle (1.4);
    \draw[softgray,dashed] (-1.7,-1.4)--(1.7,-1.4);
    \coordinate (X) at (-0.12,-1.025);
    \coordinate (Tx) at (-0.25,-0.75);
    \coordinate (P) at (0.074,-1.398);
    \draw[teal,thick,-{Stealth}] (Tx)--(P);
    \fill[ink] (Tx) circle (1.5pt) node[above] {$T(x)$};
    \fill[ink] (X) circle (1.5pt) node[left] {$x$};
    \fill[ink] (P) circle (1.5pt) node[right] {$\Psi(x)$};
    \node[ink] at (0,1.78) {$T(x)_n\geq x_n$};
    \node[softgray] at (0,-1.78) {$z_n=-1$};
  \end{scope}
  \begin{scope}[shift={(5.4,0)}]
    \draw[softgray,thick] (0,0) circle (1.4);
    \draw[softgray,dashed] (-1.7,-1.4)--(3.1,-1.4);
    \coordinate (X) at (-0.15,-1.31);
    \coordinate (Tx) at (0.10,-1.323);
    \coordinate (Y) at (2.67,-1.4);
    \coordinate (P) at (-0.509,-1.305);
    \draw[ochre,thick] (X)--(Y);
    \draw[teal,thick,-{Stealth}] (Tx)--(P);
    \fill[ink] (X) circle (1.5pt) node[above] {$x$};
    \fill[ink] (Tx) circle (1.5pt) node[above right] {$T(x)$};
    \fill[ink] (Y) circle (1.5pt) node[below] {$y$};
    \fill[ink] (P) circle (1.5pt) node[left] {$\Psi(x)$};
    \node[ink] at (0,1.78) {$T(x)_n<x_n$};
    \node[softgray] at (0,-1.78) {$z_n=-1$};
  \end{scope}
\end{tikzpicture}
\caption{The ray-exit argument near the south pole. Left: if $T(x)$
is no lower than $x$, the exit point lies in a shrinking lower cap.
Right: if $T(x)$ is lower, the line through $x$ and $T(x)$ meets the
supporting plane at $y$ beyond $T(x)$. When displacement dominates
boundary distance, $y$ moves far away relative to the height of $x$
above that plane, and the exit point $\Psi(x)$ approaches $x$.
The diagrams are schematic.}
\label{fig:fixedpoint-ray}
\end{figure}
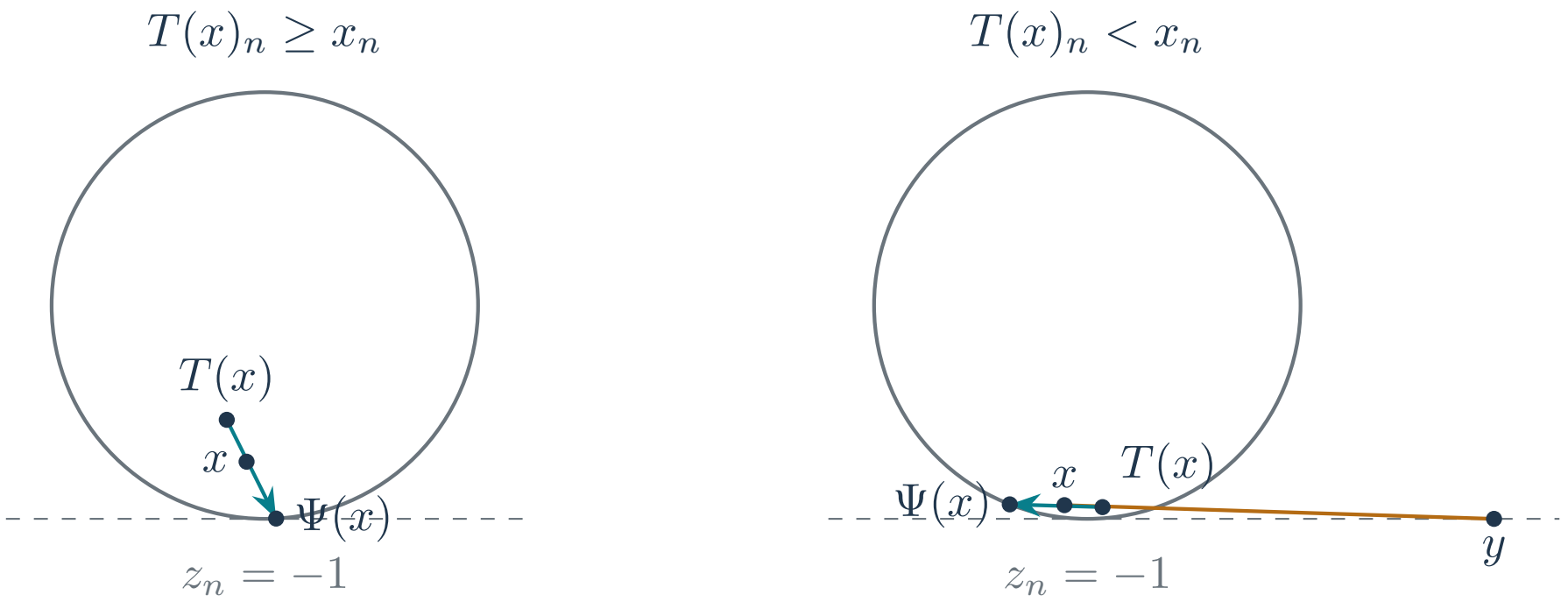

\begin{proof}[Geometric proof of \cref{thm:boundary-fixedpoint}]
First reduce to the Euclidean unit ball $B$. Choose $p\in V$ and identify
$p$ with the origin. Radial projection from $p$ identifies the boundaries
of $V$ and $B$; extending it linearly along radial segments gives a
bi-Lipschitz homeomorphism $\Phi:\overline V\to\overline B$.
For example, if $\rho B\subset V\subset RB$, the gauge
$p_V(z)=\inf\{t>0:z\in tV\}$ satisfies
$|p_V(z)-p_V(w)|\leq |z-w|/\rho$. Its restriction to the unit sphere
is bounded above and away from zero, which gives the bi-Lipschitz
bounds for the radial map and its inverse. Equivalence of norms takes
care of the original norm on $V$. Since $\Phi$ maps boundary to
boundary, the displacement hypothesis \eqref{eq:boundary-displacement}
is preserved under conjugation. So consider the situation that $T:B\to B$.

Suppose that $T$ has no fixed point. Let $\Psi(x)$ be the point where
the ray starting at $x$ and pointing away from $T(x)$ meets $\partial B$.
Then $T(x),x,\Psi(x)$ occur in this order on a line, and
$\Psi:B\to\partial B$ is continuous. {\bf Claim:} 
\begin{equation}\label{eq:ray-exit-small}
 |\Psi(x)-x|\longrightarrow0
 \qquad\text{as }\operatorname{dist}(x,\partial B)\longrightarrow0.
\end{equation}
Here is the supporting-plane picture behind this claim. Rotate
coordinates so that the nearest boundary point is the south pole, and
write $x_n=-1+\delta$, where $\delta=\operatorname{dist}(x,\partial B)$.
If $T(x)_n\geq x_n$, the ray away from $T(x)$ exits through the spherical
cap below the horizontal plane through $x$; the diameter of this cap
tends to zero with $\delta$. If $T(x)_n<x_n$, extend the line through
$x$ and $T(x)$ downwards to meet the supporting plane $z_n=-1$ at $y$.
Here $T(x)$ lies between $x$ and $y$, and $y$ is outside the ball.
In particular,
\[
 |x-y|\geq |x-T(x)|,
 \qquad \frac{|x-y|}{\delta}\geq
 \frac{|x-T(x)|}{\operatorname{dist}(x,\partial B)}\longrightarrow\infty.
\]
Thus the line from $x$ to $y$ becomes horizontal as $x$ approaches the
boundary. The ray in the opposite direction then exits the ball close
to $x$. Figure~\ref{fig:fixedpoint-ray} shows the two cases.

For a quantitative verification of this last step, set
\[
 a=|x-T(x)|>0,\quad u=\frac{x-T(x)}a,\quad
 D=1-|x|^2,\quad b=\langle x,u\rangle.
\]
Because $T(x)=x-au\in B$, expansion of $|x-au|^2<1$ gives
$b>(a^2-D)/(2a)$. If $r=|\Psi(x)-x|$, then $|x+ru|^2=1$, and hence
\[
 r=-b+\sqrt{b^2+D}.
\]
For $b\geq0$, $r\leq\sqrt D$. For $b<0$,
$-2b<D/a$, so $r\leq-2b+\sqrt D<D/a+\sqrt D$.
Since $D=(1-|x|)(1+|x|)\leq2\delta$, in both cases
\begin{equation}\label{eq:geometric-bound}
 |\Psi(x)-x|\leq \sqrt{2\delta}
       +\frac{2\delta}{|x-T(x)|}\longrightarrow0.
\end{equation}
This proves \eqref{eq:ray-exit-small} by
\eqref{eq:boundary-displacement}. Defining $\Psi(\xi)=\xi$ for
$\xi\in\partial B$ would now give a continuous retraction of
$\overline B$ onto $\partial B$, contrary to Brouwer's no-retraction
theorem. Thus $T$ has a fixed point.
\end{proof}

\FloatBarrier
\subsection{A direct proof on the convex domain}

\begin{proof}[Convex-domain proof of \cref{thm:boundary-fixedpoint}]
One can also derive a proof directly using Brouwer's fixed point theorem. Indeed, 
fix $a\in V$, and put
\[
 \rho=\dist(a,\partial V)>0,\qquad
 R=\sup_{z\in V}\norm{z-a}<\infty.
\]
For $0<t<1$, the set
\[
 K_t=a+t(\overline V-a)
\]
is a compact convex subset of $V$. The map
\[
 H_t:K_t\to K_t,\qquad H_t(x)=(1-t)a+tT(x),
\]
is continuous, so Brouwer gives a point $x_t\in K_t$ satisfying
$x_t=(1-t)a+tT(x_t)$. In particular,
\begin{equation}\label{eq:homothety-displacement}
 \norm{T(x_t)-x_t}
 =(1-t)\norm{T(x_t)-a}\le(1-t)R.
\end{equation}
Convexity and $B(a,\rho)\subset V$ give
\[
 B(x_t,(1-t)\rho)
 =(1-t)B(a,\rho)+tT(x_t)\subset V.
\]
Consequently,
\begin{equation}\label{eq:homothety-ratio}
 \frac{\norm{T(x_t)-x_t}}{\dist(x_t,\partial V)}
 \le\frac R\rho.
\end{equation}
Take a sequence $t_n\uparrow1$. If the points $x_{t_n}$ have an
interior accumulation point, \eqref{eq:homothety-displacement} and
continuity show that it is a fixed point of $T$. Otherwise their
distance to $\partial V$ tends to zero, contradicting
\eqref{eq:boundary-displacement} and \eqref{eq:homothety-ratio}.
\end{proof}

\begin{remark}[The displacement condition and the analytic task]
\label{rem:boundary-analytic-task}
The condition of \cref{thm:boundary-fixedpoint} is about relative
displacement, not the distance of successive iterates from the
boundary. For example, $T(x)=-x$ on the unit ball satisfies that
condition although it preserves the distance to the boundary.

By \cref{lem:mesh-distance}, it suffices for the Thurston operator
to prove that, for every fixed $K<\infty$, the set
\begin{equation}
 \{x\in W:\norm{T(x)-x}_\infty\le K m(x)\} \label{eq:displacement} 
\end{equation} 
stays a positive distance from $\partial W$. The estimates in the
following sections establish precisely this assertion. They are
applied before any compactness assertion about a degenerating
family of entire functions is needed.
\end{remark}

\section{Quantitative estimates for approximate fixed points}
\label{sec:setting}\label{sec:apriori}

I now prove the boundary estimate required by
\cref{thm:boundary-fixedpoint}. Stronger derivative and inverse-height
bounds are not needed for this and are given in the appendices. The estimates use only the fixed finite labelled part of the
real graph. They do not require compactness of the natural family as
singular values collide.

\subsection{Fixed branches and the relative pullback error}
Keep the finite marking, \(\pi\), and \(\cC\) from
\cref{sec:thurston}. Fix a finite collection
\(\cB_0\) of real branches of \(f_0\) containing every branch used by
the marked pullbacks. A finite number of additional fixed branches may be included when
inverse-height bounds on them are desired. Under natural equivalence,
write \(\cB_f\) for the corresponding branches of \(f\). Even if the
real critical set is infinite, these branch labels remain fixed and finite in number throughout.
When a configuration is fixed, the subscript is  suppressed and the following notation 
\(f=f_x\), \(L=L_x\), and \(F=F_x\) is used; the same notation applies
after a preliminary positive affine source change.

For a finite real interval \(J\), define its total critical multiplicity by
\[
 d_f(J)=\sum_{\substack{c\in J\\f'(c)=0}}\ord_c(f').
\]
For a finite nonempty set \(S\), let
\[
 \sep S=\min_{\substack{u,v\in S\\u\ne v}}\abs{u-v},
 \qquad \sep S=+\infty\quad\text{if }\#S=1.
\]

\begin{definition}[Admissible approximate fixed configuration]
\label{def:admissible}
Fix \(K\ge0\). An admissible configuration consists of a marked map
\(f\) in the real natural class, a target configuration \(x\in W\),
and its finite labelled unscaled preimages
\(\widehat y_1<\cdots<\widehat y_s\), satisfying
\eqref{eq:singular-marking}--\eqref{eq:critical-marking}.
Define \(L,y,F\) by \eqref{eq:affine}. The expansive
conditions \eqref{eq:noflat} and \eqref{eq:expansive} are required, and
\begin{equation}\label{eq:error}
 \norm{x-y}_\infty\le K m,
 \qquad m=\min_{1\le i<s}(x_{i+1}-x_i).
\end{equation}
No smallness condition is imposed on \(K\).
\end{definition}

The estimates are unchanged by a positive affine change of the source
coordinate: the normalized map \(F\) and the vector \(y\) are unchanged.
The source normalizations of \cref{sec:natural-parameters-new} are used only
for derivative bounds on the unscaled maps $f$, before composition
with the affine rescaling $M_x$.

\begin{theorem}[Uniform marked-gap separation]
\label{thm:apriori}
Fix the real natural class, the finite branch and critical labels, the
compatible expansive data, and \(K<\infty\). There is \(m_*>0\), depending
only on these fixed data and on \(K\), such that every configuration in
\cref{def:admissible} satisfies
\begin{equation}\label{eq:mainseparation}
 x_{i+1}-x_i\ge m_*\qquad(1\le i<s).
\end{equation}
If $s\ge3$, then \cref{lem:mesh-distance} also gives
$\delta_W(x)\ge m_*/2$. If there are at least two singular
values, every pair of distinct singular values is separated by
at least $m_*$. The constants depend on $K$.
No hypothesis is imposed on the behavior as the real variable tends
to either infinity, or on the
total number of real critical points or asymptotic tracts.
\end{theorem}

This is the a priori estimate needed for the realisation theorem. Its
proof is completed in \cref{gap:separation}; the ensuing
\cref{cor:existence-new} converts it into the boundary-displacement
condition and hence into existence of a fixed point. Stronger derivative
and inverse-height consequences after source normalization are not used
in that argument and are deferred to the appendices.

\paragraph{The analytic mechanism.}
For \(J=[\widehat y_1,\widehat y_s]\), all endpoint and critical values
lie in \([0,1]\), so \(f(J)\subset[0,1]\). The fixed labels bound
\(d_f(J)\). In \cref{ana:diameter} the following inequality is proved: 
\[
 \diam_U J\le C_0+C_1\log\bigl(1+\log(1/m)\bigr),
 \qquad
 U\text{ the component of }f^{-1}D(0,4)\text{ containing }J.
\]
The annulus and path estimates give the required double-logarithmic
upper-length mechanism. All ingredients needed here are proved below. Koebe estimates convert this diameter
bound into polynomial bounds in \(1+\log(1/m)\) for \(F'\) and \(F''\).
Expansiveness, together with \cref{eq:error} then forces
\(1\le C m(1+\log(1/m))^B\), which therefore rules out \(m\to0\).
Only after separation is established are the source anchor points used.

\section{Elementary pullback topology}\label{ana:topology}

This section records the topology needed later. Throughout, the singular
set is understood through the covering characterization: outside
$S(f)$, the map obtained by removing its full preimage is a holomorphic
covering. No bound on the number of complex critical points or asymptotic
tracts is used.

\begin{lemma}[Disk pullbacks and their real traces]\label{ana:disk-topology}
Let $f$ be a nonconstant entire function, let $D(w_0,R)$ be a finite
Euclidean disk, and let $U$ be a component of its preimage.
Then $U$ is a proper simply connected plane domain.
If $f$ is real, $w_0$ is real, and $U$ meets $\R$, then $U$ is
conjugation invariant, $U\cap\R$ is an interval, and that interval is
a hyperbolic geodesic in $U$.
\end{lemma}
\begin{proof}
Let $\gamma$ be a Jordan curve in $U$. The compact curve $\gamma$
satisfies
\[
 \max_{z\in\gamma}\abs{f(z)-w_0}<R.
\]
The maximum principle implies the same strict inequality throughout the
bounded interior of $\gamma$. The interior, together with $\gamma$,
belongs to the same component of the disk preimage. Thus every Jordan
curve in $U$ encloses only points of $U$, proving simple connectivity.
If $U=\C$, then $f$ is bounded on $\C$, contrary to Liouville's theorem
and nonconstancy.

Now assume the real hypotheses. Conjugation sends $U$ to a component
of the same preimage. A real point of $U$ belongs to both components,
so the components coincide. Choose $t_0\in U\cap\R$ and the Riemann map
$h:\D\to U$ normalized by
\[
 h(0)=t_0,\qquad h'(0)>0.
\]
Uniqueness of this normalization gives
$h(\bar z)=\overline{h(z)}$. Hence $h((-1,1))\subset U\cap\R$.
Conversely, if $h(z)$ is real, then
$h(z)=\overline{h(z)}=h(\bar z)$; injectivity gives $z=\bar z$.
Therefore
\[
 U\cap\R=h((-1,1)).
\]
It is an interval. The real diameter is a geodesic for the
hyperbolic metric of $\D$, and conformal maps preserve hyperbolic
geodesics.
\end{proof}

\begin{lemma}[A disk containing one singular value]\label{ana:singleton}
Let $f$ be nonconstant and entire, let $v\in S(f)$, and suppose
\[
 D(v,R)\cap S(f)=\{v\}.
\]
For every component $T$ of $f^{-1}D(v,R)$, precisely one of the
following holds.
\begin{enumerate}
\item The component contains a unique point $c_T$ with $f(c_T)=v$.
The map $f:T\to D(v,R)$ is proper of finite degree $e\ge1$.
Its only possible critical point is $c_T$, of multiplicity $e-1$,
and it is conformally equivalent to the power map $z\mapsto z^e$.
\item The component contains no preimage of $v$, and
\[
 f:T\longrightarrow D(v,R)\setminus\{v\}
\]
is a universal covering.
\end{enumerate}
If $f$ and $v$ are real and $T$ meets $\R$, then $T$ is symmetric.
In the first case $c_T$ is real. A symmetric Riemann map
$h:\D\to T$, normalized by $h(0)=c_T$ and $h'(0)>0$, satisfies
\begin{equation}\label{ana:power-coordinate}
 f(h(z))=v+\sigma Rz^e,\qquad \sigma\in\{-1,1\}.
\end{equation}
In the second case $T\cap\R$ is an unbounded interval.
\end{lemma}
\begin{proof}
By \cref{ana:disk-topology}, $T$ is simply connected. The set
$Z=T\cap f^{-1}(v)$ is closed and discrete in $T$. Removing a closed
discrete set from a plane domain leaves it path connected: start with
a compact path between two points and make small detours around the
finitely many removed points met in a compact neighborhood of that
path. Consequently $T\setminus Z$ is a component of the preimage of
$D(v,R)\setminus\{v\}$, and the covering characterization gives a
connected covering
\[
 f:T\setminus Z\longrightarrow D(v,R)\setminus\{v\}.
\]
The induced map on fundamental groups is injective, with image a
subgroup of $\mathbb Z$.

If $Z$ contained two distinct points, small puncture loops around those
points would not commute in $\pi_1(T\setminus Z)$. For example, a
Riemann map to the disk and then the inclusion that fills all other
punctures send these loops to the two free generators of a twice
punctured disk. Their commutator remains nontrivial. This is impossible
for a group injecting into $\mathbb Z$. Thus $Z$ has at most one point.

If $Z$ is empty, the simply connected source makes the covering
universal. If $Z=\{c_T\}$, let $e$ be the order of the zero of
$f-v$ at $c_T$. A positive small loop around $c_T$ maps to a loop of
winding number $e$ around $v$. The covering subgroup is therefore
$e\mathbb Z$, and the covering has $e$ sheets.
Filling the puncture makes $f:T\to D(v,R)$ proper: away from $v$,
a finite covering is proper, and sufficiently close to $v$ all $e$
preimages are already accounted for in a small neighborhood of $c_T$.
There are no critical points away from $c_T$, because the punctured
restriction is a covering. The multiplicity at $c_T$ is $e-1$.

The conformal power description can also be seen directly.
Since the sole zero of $f-v$ has order $e$, the function
$(f-v)/R$ has a holomorphic $e$th root on $T$. This root is a proper
map to $\D$ of degree one, because its composition with the power map
has degree $e$. It is therefore conformal.

For the real assertions, symmetry follows from
\cref{ana:disk-topology}. In the first case, uniqueness of the
preimage of $v$ gives $c_T=\bar c_T$. In the symmetric Riemann
coordinate, the proper power description has a unit complex
coefficient. Reality on the real diameter makes that coefficient
$1$ or $-1$, proving \eqref{ana:power-coordinate}.

In the universal-cover case, write $T\cap\R=(\alpha,\beta)$,
allowing infinite endpoints. On this interval, $f-v$ has a constant
nonzero sign, and $f'$ is real and never zero. Hence $f$ is strictly
monotone there. If both endpoints were finite, continuity would put
their images on $\partial D(v,R)$. The fixed sign of $f-v$ would
make both images equal to the same one of $v-R$ and $v+R$,
contradicting strict monotonicity. At least one endpoint is infinite.
\end{proof}

\begin{remark}\label{ana:covering-reference}
The covering classification in \cref{ana:singleton} is the
punctured-disk instance of the standard finite-cover or
universal-cover dichotomy for doubly connected target domains.
See \cite[Lemma~2.7 and Proposition~2.9]{BFR}. The proof above
includes the real-trace statement needed here.
\end{remark}

\begin{lemma}[A bounded target interval contains the real span]
\label{ana:range}
Let $f$ be real entire and suppose that every real critical value
belongs to $[A,B]$. If $a<b$ are finite real points with
$f(a),f(b)\in[A,B]$, then
\[
 f([a,b])\subset[A,B].
\]
\end{lemma}
\begin{proof}
If a value larger than $B$ occurred, the maximum on the compact
interval $[a,b]$ would be attained at an interior critical point
and would exceed $B$. A value below $A$ gives the corresponding
contradiction at an interior minimum.
\end{proof}

\section{Hyperbolic metrics and conformal radii}\label{ana:metric}

I use hyperbolic metrics of curvature $-1$, denoting their density,
distance, and path length by $\lambda_U$, $d_U$, and $\ell_U$.
The two elementary principles used below are
\begin{equation}\label{ana:metric-monotonicity}
 U_1\subset U_2\quad\Longrightarrow\quad
 \lambda_{U_2}\le\lambda_{U_1}\ \hbox{on }U_1,
\end{equation}
and the fact that a holomorphic covering between hyperbolic domains
is a local isometry. These follow from Schwarz--Pick, applied to
universal disk covers.

The disk and annulus formulas in this normalization are
\begin{align}
 \lambda_{D(w_0,R)}(w)
 &=\frac{2R}{R^2-\abs{w-w_0}^2},
 \label{ana:disk-metric}\\
 \lambda_{\{r<\abs{w-w_0}<R\}}(w)
 &=
 \frac{\pi}
 {\abs{w-w_0}\,\ell\,
  \sin\bigl(\pi\log(\abs{w-w_0}/r)/\ell\bigr)},
 \qquad \ell=\log(R/r).
 \label{ana:annulus-metric}
\end{align}
For the annulus formula, use the logarithmic covering from the strip
$0<\Re z<\ell$. That strip has density
$\pi/(\ell\sin(\pi\Re z/\ell))$, obtained by mapping it to a
half-plane; the derivative of the logarithm supplies the factor
$1/\abs{w-w_0}$.
I also use
\begin{equation}\label{ana:disk-distance}
 d_{\D}(0,z)=\log\frac{1+\abs z}{1-\abs z}.
\end{equation}

\begin{lemma}[Classical univalent estimates]\label{ana:koebe}
Let $h$ be univalent on $\D$, and put $R=\abs{h'(0)}$.
Then
\begin{align}
 D(h(0),R/4)&\subset h(\D),\label{ana:quarter}\\
 \abs{h(z)-h(0)}
 &\le R\frac{\abs z}{(1-\abs z)^2},\label{ana:growth}\\
 R\frac{1-\abs z}{(1+\abs z)^3}
 &\le \abs{h'(z)}
 \le R\frac{1+\abs z}{(1-\abs z)^3}.
 \label{ana:derivative-distortion}
\end{align}
If $h$ is univalent on a proper plane domain $\Omega$, then
\begin{equation}\label{ana:pre-schwarzian}
 \left|\frac{h''(z)}{h'(z)}\right|
 \le\frac{4}{\dist(z,\partial\Omega)}.
\end{equation}
\end{lemma}
\begin{proof}
The first three assertions are Koebe's quarter, growth, and distortion
theorems; see \cite{Duren}. For the last assertion, take
$0<r<\dist(z,\partial\Omega)$ and normalize
\[
 H(\zeta)=\frac{h(z+r\zeta)-h(z)}{r h'(z)}.
\]
This is a normalized univalent map on $\D$, whose second coefficient
is $rh''(z)/(2h'(z))$. The classical second-coefficient bound
$|a_2|\le2$ gives \eqref{ana:pre-schwarzian} on letting
$r\uparrow\dist(z,\partial\Omega)$.
\end{proof}

For a proper simply connected domain $U$ and $a\in U$, its
\emph{conformal radius} is
\[
 R_U(a)=\abs{h'(0)},
\]
where $h:\D\to U$ is any Riemann map with $h(0)=a$.
Equivalently, $R_U(a)=2/\lambda_U(a)$.

\begin{lemma}[Conformal-radius comparison]\label{ana:radius}
If $U_1\subset U_2$ are proper simply connected domains and $a\in U_1$,
then
\begin{equation}\label{ana:radius-monotonicity}
 R_{U_1}(a)\le R_{U_2}(a).
\end{equation}
If $a,b\in U$, where $U$ is proper and simply connected, then
\begin{equation}\label{ana:radius-distortion}
 e^{-2d_U(a,b)}
 \le\frac{R_U(b)}{R_U(a)}
 \le e^{2d_U(a,b)}.
\end{equation}
Moreover, if $d_U(a,b)\le D$, then
\begin{equation}\label{ana:euclidean-from-hyperbolic}
 \abs{b-a}\le
 \frac{R_U(a)}4(e^{2D}-1).
\end{equation}
\end{lemma}
\begin{proof}
The monotonicity assertion follows from
\eqref{ana:metric-monotonicity} and $R_U=2/\lambda_U$.
Choose $h:\D\to U$ with $h(0)=a$, and write $b=h(z)$,
$r=\abs z$. Conformal invariance gives
\[
 R_U(b)=(1-r^2)\abs{h'(z)}.
\]
The lower derivative estimate in \cref{ana:koebe} therefore yields
\[
 R_U(b)\ge
 R_U(a)\left(\frac{1-r}{1+r}\right)^2
 =R_U(a)e^{-2d_U(a,b)}.
\]
Interchanging $a,b$ proves the upper estimate.

For the last assertion, \eqref{ana:disk-distance} gives
$r\le\tanh(D/2)$. The growth estimate and the identity
\[
 \frac{\tanh(D/2)}{(1-\tanh(D/2))^2}
 =\frac{e^{2D}-1}{4}
\]
give \eqref{ana:euclidean-from-hyperbolic}.
\end{proof}

\section{Double-logarithmic length of real target intervals}\label{ana:arches}

The next lemma bounds the Poincar\'e length of a real interval that
avoids the singular set. It is the target estimate needed for the
real corridors in \cref{ana:diameter}. The proof splits the interval
at its midpoint and estimates each half directly, using disks near
the endpoint and on transition bands, and puncture-free annuli
between successive distance scales.

\begin{lemma}[A real interval in a punctured disk]\label{ana:arch}
Let $S\subset[0,1]$ consist of $q\ge1$ points, and put
\[
 Y=D(0,4)\setminus S.
\]
Let $0<r\le1/8$ and let $0\le u\le w\le1$ satisfy
\[
 [u,w]\cap S=\varnothing,
 \qquad \dist(u,S),\dist(w,S)\ge r.
\]
Put
\[
 \Lambda=1+\log(1/r).
\]
Then
\begin{equation}\label{ana:arch-bound}
 \ell_Y([u,w])\le13q+4q\log\Lambda.
\end{equation}
\end{lemma}
\begin{proof}
The case $u=w$ is immediate. Suppose $u<w$, and put
\[
 c=\frac{u+w}{2},\qquad h=\frac{w-u}{2}.
\]
First consider the left half, parameterized by $x=u+t$,
$0\le t\le h$. For $0<t\le h$, the disk $D(u+t,t)$ meets the
real axis in $(u,u+2t)\subset(u,w)$. Since all points of $S$ are
real and $[u,w]$ avoids $S$, this disk contains no puncture.
It also lies in $D(0,4)$, because
$\abs{u+t}+t=u+2t\le w\le1$.
The disk metric and domain monotonicity therefore give
\begin{equation}\label{ana:vertical-simple}
 \lambda_Y(u+t)\le\frac2t\qquad(0<t\le h).
\end{equation}

List the distinct distances from $u$ to the points of $S$ as
\[
 r_1<\cdots<r_n\le1,\qquad 1\le n\le q,
\]
and append $r_{n+1}=2$. Then $r_1\ge r$.
The initial part $0\le t\le\min\{h,r_1/2\}$ lies in
$D(u,r_1)\subset Y$. By \eqref{ana:disk-metric}, its length is at most
\begin{equation}\label{ana:initial-cost}
 \int_0^{r_1/2}\frac{2r_1}{r_1^2-t^2}\,\dd t=\log3.
\end{equation}

The transition bands
\[
 [r_j/2,2r_j]\cap[0,h],\qquad 1\le j\le n,
\]
correspond to real subintervals of total hyperbolic length at most
$2q\log4$ by \eqref{ana:vertical-simple}; overlaps are harmless.
Outside the initial part and these bands, every remaining piece lies in
a gap of the form
\[
 2r_j<t<r_{j+1}/2
\]
for some $j\le n$; here $r_{n+1}/2=1\ge h$.
For a nonempty gap, put
$\ell_j=\log(r_{j+1}/r_j)>\log4$.
The annulus
\[
 A_j=\{z:r_j<\abs{z-u}<r_{j+1}\}
\]
contains no singular value. It is contained in $D(0,4)$, since its
outer radius is at most $2$ and $0\le u\le1$. Consequently
$A_j\subset Y$.
Using \eqref{ana:annulus-metric}, then substituting
$s=\log(t/r_j)$, bounds the length of the part of the left half
in this gap by the integral over the full radial gap:
\begin{align}
 \int_{2r_j}^{r_{j+1}/2}\lambda_{A_j}(u+t)\,\dd t
 &=
 \int_{\log2}^{\ell_j-\log2}
       \frac{\pi\,\dd s}{\ell_j\sin(\pi s/\ell_j)}
 \notag\\
 &=2\log\cot\frac{\pi\log2}{2\ell_j}
 \le2\log(2\ell_j).
 \label{ana:annular-gap-cost}
\end{align}
Here $\cot t\le1/t$ for $0<t<\pi/2$, and $\pi\log2>1$.
Since
\[
 \ell_j\le\log(2/r),
\]
and there are at most $q$ gaps, the left half satisfies
\begin{equation}\label{ana:one-vertical-cost}
 \ell_Y([u,c])\le\log3+2q\log4+
 2q\log\bigl(2\log(2/r)\bigr).
\end{equation}
For the right half write $x=w-t$, $0\le t\le h$.
The disk $D(w-t,t)$ meets the real axis in
$(w-2t,w)\subset(u,w)$ and lies in $D(0,4)$, so again
$\lambda_Y(w-t)\le2/t$ for $t>0$.
Repeating the same decomposition with the distinct distances from
$w$ to $S$ gives the bound in \eqref{ana:one-vertical-cost} for
$\ell_Y([c,w])$ as well. Adding the two halves yields
\[
 \ell_Y([u,w])
 \le
 2\log3+4q\log4+4q\log2+
 4q\log\log(2/r).
\]
Because $q\ge1$,
\[
 2\log3+4q\log4+4q\log2\le13q.
\]
Also $\log(2/r)\le1+\log(1/r)=\Lambda$.
These observations prove \eqref{ana:arch-bound}.
\end{proof}

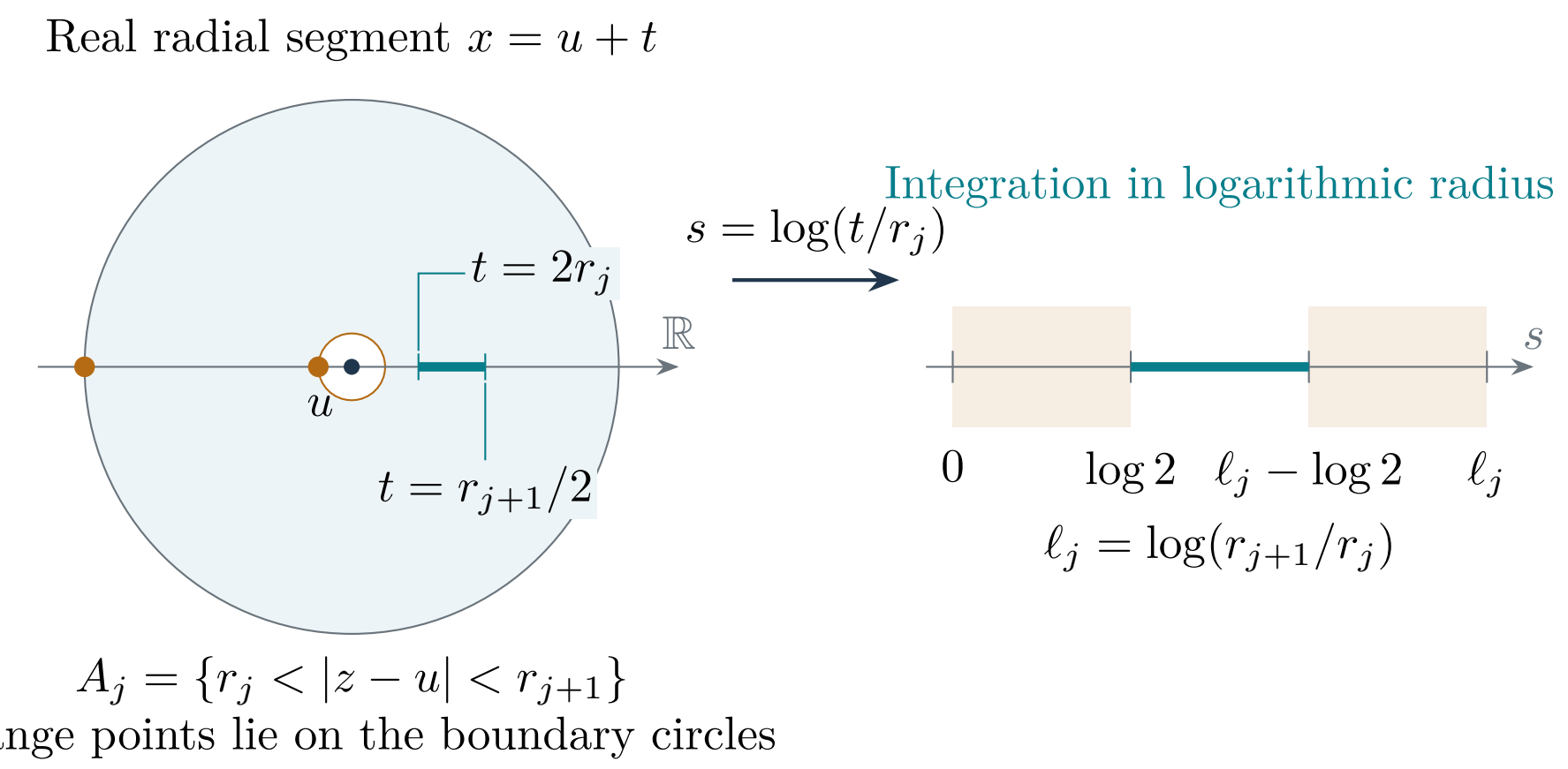
\begin{figure}[htbp]
\centering
\begin{tikzpicture}[x=1cm,y=1cm,>=Stealth,font=\small]
 \fill[pale,even odd rule] (0,0) circle (2) (0,0) circle (0.25);
 \draw[softgray] (0,0) circle (2);
 \draw[ochre] (0,0) circle (0.25);
 \draw[->,softgray] (-2.35,0)--(2.45,0) node[above] {$\R$};
 \fill[ink] (0,0) circle (1.7pt);
 \fill[ochre] (-0.25,0) circle (2.2pt);
 \fill[ochre] (-2,0) circle (2.2pt);
 \draw[line width=2pt,teal] (0.5,0)--(1,0);
 \draw[teal] (0.5,-0.1)--(0.5,0.1);
 \draw[teal] (1,-0.1)--(1,0.1);
 \draw[teal,thin] (0.5,0.12)--(0.5,0.7)--(0.85,0.7);
 \node[anchor=west,fill=pale,inner sep=1pt] at (0.85,0.7) {$t=2r_j$};
 \draw[teal,thin] (1,-0.12)--(1,-0.7);
 \node[below,fill=pale,inner sep=1pt] at (1,-0.7) {$t=r_{j+1}/2$};
 \node[below left] at (0,-0.08) {$u$};
 \node[align=center] at (0,2.45)
   {Real radial segment $x=u+t$};
 \node[align=center] at (0,-2.55)
   {$A_j=\{r_j<|z-u|<r_{j+1}\}$\\Orange points lie on the boundary circles};
 \draw[->,ink,thick] (2.85,0.65)--(4.1,0.65);
 \node[above,align=center] at (3.48,0.7) {$s=\log(t/r_j)$};
 \begin{scope}[shift={(4.5,0)}]
  \fill[ochre!12] (0,-0.45) rectangle (1.3333,0.45);
  \fill[ochre!12] (2.6667,-0.45) rectangle (4,0.45);
  \draw[->,softgray] (-0.2,0)--(4.35,0) node[above] {$s$};
  \foreach \x in {0,1.3333,2.6667,4}
    \draw[softgray] (\x,-0.12)--(\x,0.12);
  \draw[line width=2pt,teal] (1.3333,0)--(2.6667,0);
  \node[below] at (0,-0.5) {$0$};
  \node[below] at (1.3333,-0.5) {$\log2$};
  \node[below] at (2.6667,-0.5) {$\ell_j-\log2$};
  \node[below] at (4,-0.5) {$\ell_j$};
  \node[teal,align=center] at (2,1.35) {Integration in logarithmic radius};
  \node[align=center] at (2,-1.35) {$\ell_j=\log(r_{j+1}/r_j)$};
 \end{scope}
\end{tikzpicture}
\caption{The direct annular estimate along the real interval.
 The teal segment is $u+t$ with $2r_j\le t\le r_{j+1}/2$;
 only its intersection with the left half $[u,c]$ is needed.
 In logarithmic radius it becomes $[\log2,\ell_j-\log2]$, and
 integrating the annulus density gives
 $2\log\cot(\pi\log2/(2\ell_j))=O(\log\ell_j)$.
 The illustrated ratio is $r_{j+1}/r_j=8$; the argument applies
 whenever this ratio exceeds $4$. Orange bands on the right mark
 transition coordinates, not removed subsets of the target.}
\label{fig:annuli}
\end{figure}

\Cref{fig:annuli} illustrates the annular integration in the proof.

\section{Hyperbolic diameter of a finite real interval}
\label{ana:finite-interval}

The next theorem assumes only separation of singular and endpoint
values and a finite critical-multiplicity budget. It does not assume
$\|T(x)-x\|_\infty\le K m(x)$, or even require that $f$ arises from
a Thurston pullback.

\begin{theorem}[Finite-interval diameter estimate]\label{ana:diameter}
Let $f$ be a nonconstant real entire map with
\[
 S(f)\subset[0,1],\qquad \#S(f)=q\ge1.
\]
Let $d\ge0$ be an integer and $0<m\le1$. Suppose that $J=[a,b]$,
where $a<b$ are finite real points, satisfies
\begin{enumerate}
\item $f(a),f(b)\in[0,1]$;
\item $\sep(S(f))\ge m$, where the separation of a singleton is
understood as $+\infty$;
\item for $t=a,b$, either $f(t)\in S(f)$ or
$\dist(f(t),S(f))\ge m$;
\item
\[
 \sum_{\substack{c\in J\\ f'(c)=0}}\ord_c(f')\le d.
\]
\end{enumerate}
Let $U$ be the component of $f^{-1}D(0,4)$ containing $J$.
Define
\begin{align}
 M&=q(d+1),&
 B_d&=2\log\frac{1+2^{-1/(d+1)}}{1-2^{-1/(d+1)}},
 \label{ana:core-constants}\\
 \Lambda(m)&=1+\log(8/m),&
 D_{q,d}(m)&=M B_d+(M+1)\bigl(13q+4q\log\Lambda(m)\bigr).
 \label{ana:D-definition}
\end{align}
Then
\begin{equation}\label{ana:diameter-bound}
 \diam_U J\le D_{q,d}(m).
\end{equation}
In particular, for constants depending only on $q,d$,
\[
 \diam_U J\le C_0+C_1\log\bigl(1+\log(1/m)\bigr).
\]
\end{theorem}
\begin{proof}
All real critical values lie in $[0,1]$, so \cref{ana:range} gives
$f(J)\subset[0,1]$. This also shows that $J$ lies in one component $U$
of the indicated disk preimage.
By \cref{ana:disk-topology}, $U$ is a proper simply connected symmetric
domain and its real trace is a hyperbolic geodesic.

\paragraph{Step 1: small singular components meeting $J$ are finite power
components.}
Put $r=m/8$. For each $v\in S(f)$, the disk $D(v,2r)$ contains
no other singular value and lies in $D(0,4)$.
Let $T^+$ be a component of its preimage meeting $J$.
If $T^+$ were the universal-cover alternative of
\cref{ana:singleton}, then $T^+\cap\R$ would be an unbounded interval.
An unbounded interval meeting $[a,b]$ contains at least one of $a,b$.
Write that endpoint as $t$.
The endpoint hypothesis prevents $f(t)$ from lying in $D(v,2r)$:
if it is regular, its distance from $v$ is at least $m$; if it is a
different singular value, its distance from $v$ is again at least $m$.
The remaining possibility $f(t)=v$ is excluded because a universal
component omits $v$. This contradiction proves that $T^+$ is a finite
power component.

Let $c_T$ be its unique preimage of $v$ in $T^+$. Symmetry makes $c_T$ real.
It belongs to $J$. Indeed, if $c_T$ were outside $J$, the real interval
$T^+\cap\R$ joining $c_T$ to a point of $J$ would contain an endpoint
$t$ of $J$. The preceding endpoint argument would force $f(t)=v$.
Uniqueness of the center would then give $t=c_T$, a contradiction.
If $e_T$ denotes the degree of the component, its critical
multiplicity is $e_T-1$ at this center. Thus
\begin{equation}\label{ana:degree-bound}
 1\le e_T\le d+1.
\end{equation}
This argument explicitly allows a center to equal $a$ or $b$ and
allows an endpoint target to be singular.

\paragraph{Step 2: count the components and bound their inner cores.}
Split $J$ at its interior real critical points. There are at most
$d+1$ resulting closed intervals of strict monotonicity. Each has at
most one preimage of each singular value, hence at most $q$ such
preimages. Critical endpoints may be counted twice, which is harmless
for an upper bound. There are therefore at most
$M=q(d+1)$ points of $J$ mapping to $S(f)$, and at most $M$ of the
components $T^+$ just considered. When $d=0$, this is simply the
observation that $f$ is strictly monotone on $J$ and has at most $q$
such preimages.

In the symmetric coordinate for $T^+$,
\[
 f(h(\zeta))=v+\sigma\,2r\,\zeta^{e_T}.
\]
Its inner preimage of the closed disk $\overline D(v,r)$ corresponds
to
\[
 |\zeta|\le\theta_T,\qquad \theta_T=2^{-1/e_T}.
\]
The intersection of this inner preimage with $J$ is a real interval
containing $c_T$; call it an inner core. Its hyperbolic length in
$T^+$ is at most the full diameter length
\[
 \int_{-\theta_T}^{\theta_T}\frac{2\,\dd t}{1-t^2}
 =
 2\log\frac{1+\theta_T}{1-\theta_T}
 \le B_d.
\]
Since $T^+\subset U$, its length in $U$ is no larger.
Distinct inner cores are disjoint. The argument includes components
of degree one: regular preimages of singular values must also be
covered by an inner core.

\paragraph{Step 3: bound the real corridors by their target intervals.}
Between the at most $M$ inner cores there are at most $M+1$ closed
remaining corridors, ignoring empty intervals. Every critical point
of $J$ is a center already included in a core. Moreover, every
point of $J$ whose image lies in any $D(v,r)$ belongs to that inner
core, by the power description of $T^+$. It follows that on each
corridor $Q=[\alpha,\beta]$,
\[
 f'\ne0,\qquad \dist(f(t),S(f))\ge r\quad(t\in Q).
\]
Thus $f$ maps $Q$ strictly monotonically onto a real interval
$[u,w]$ avoiding $S(f)$, where
\[
 u=\min\{f(\alpha),f(\beta)\},\qquad
 w=\max\{f(\alpha),f(\beta)\}.
\]
Both endpoints belong to $[0,1]$ and have distance at least $r$
from $S(f)$, so \cref{ana:arch} applies directly to $[u,w]$.

Put $Y=D(0,4)\setminus S(f)$ and $V=U\setminus f^{-1}(S(f))$.
The restriction
\[
 f:V\longrightarrow Y
\]
is a covering: $V$ is connected, since it is obtained by removing a
closed discrete set from the plane domain $U$, and $V$ is a component
of the full preimage of $Y$.
The covering is a hyperbolic local isometry, and $f|_Q$ traverses
$[u,w]$ exactly once, in either orientation. Since $Q\subset V\subset U$,
domain monotonicity and \cref{ana:arch} give
\[
 \ell_U(Q)\le\ell_V(Q)
 =\ell_Y([u,w])\le13q+4q\log\Lambda(m).
\]

\paragraph{Step 4: add the lengths along $J$.}
The inner cores and the remaining corridors cover $J$ and have
disjoint interiors. Adding their lengths therefore gives
\[
 \ell_U(J)\le M B_d+(M+1)\bigl(13q+4q\log\Lambda(m)\bigr).
\]
Since the real trace of $U$ is a hyperbolic geodesic in the simply
connected domain $U$, its subinterval $J$ satisfies
\[
 \diam_U J=d_U(a,b)=\ell_U(J).
\]
The preceding length bound proves
\eqref{ana:diameter-bound}.
Finally,
\[
 \Lambda(m)=1+\log8+\log(1/m)
 \le(1+\log8)\bigl(1+\log(1/m)\bigr),
\]
which gives the final form of the bound.
\end{proof}

\begin{figure}[htbp]
\centering
\begin{tikzpicture}[x=1.15cm,y=1cm,>=Stealth,font=\small]
 \node[anchor=west,font=\small\bfseries] at (0,5.5)
   {Source: the real interval $J=[a,b]\subset U$};
 \foreach \lo/\hi in {1.25/2.75,4.35/5.65,7.15/8.85}{
   \fill[pale,rounded corners=3pt] (\lo,3.9) rectangle (\hi,4.55);
 }
 \draw[softgray] (-0.15,4.15)--(10.15,4.15);
 \draw[very thick,teal] (0,4.15)--(1.7,4.15);
 \draw[very thick,teal] (2.3,4.15)--(4.7,4.15);
 \draw[very thick,teal] (5.3,4.15)--(7.6,4.15);
 \draw[very thick,teal] (8.4,4.15)--(10,4.15);
 \draw[line width=4pt,ochre!70] (1.7,4.15)--(2.3,4.15);
 \draw[line width=4pt,ochre!70] (4.7,4.15)--(5.3,4.15);
 \draw[line width=4pt,ochre!70] (7.6,4.15)--(8.4,4.15);
 \foreach \x in {0,10}{\fill[ink] (\x,4.15) circle (2pt);}
 \foreach \x in {2,8}{\fill[ink] (\x,4.15) circle (2.3pt);}
 \draw[fill=white,ink,thick] (5,4.15) circle (2.4pt);
 \node[below] at (0,4.03) {$a$};
 \node[below] at (10,4.03) {$b$};
 \node[above] at (2,4.55) {power core};
 \node[above] at (5,4.55) {degree-one core};
 \node[above] at (8,4.55) {power core};
 \draw[decorate,decoration={brace,mirror,amplitude=4pt},teal]
   (2.3,3.82)--(4.7,3.82);
 \node[teal,below] at (3.5,3.65) {$Q=[\alpha,\beta]$};
 \node[teal,below] at (6.45,4.03) {corridor};
 \node[ochre] at (8,3.35) {$\ell_U(\text{core})\le B_d$};
 \draw[->,ink,thick] (3.5,3.12)--(3.5,2.55);
 \node[anchor=west] at (3.7,2.83) {$f|_Q$ strictly monotone};
 \node[anchor=west,font=\small\bfseries] at (0,2.15)
   {Target: the real interval $f(Q)=[u,w]\subset Y$};
 \draw[->,softgray] (-0.15,1.35)--(10.2,1.35) node[right] {$\R$};
 \foreach \x/\lab in {0.5/v_1,1.3/v_2,8.5/v_3,9.4/v_4}{
   \fill[ochre] (\x,1.35) circle (2.5pt);
   \node[below] at (\x,1.25) {$\lab$};
 }
 \draw[line width=2pt,teal] (2.4,1.35)--(7.3,1.35);
 \fill[ink] (2.4,1.35) circle (2pt);
 \fill[ink] (7.3,1.35) circle (2pt);
 \node[below] at (2.4,1.25) {$u$};
 \node[below] at (7.3,1.25) {$w$};
 \node[teal,below] at (4.85,1.15) {$\dist([u,w],S(f))\ge r$};
 \node[align=center] at (5,0.2)
   {$\ell_U(Q)\le\ell_{U\setminus f^{-1}(S(f))}(Q)
     =\ell_Y([u,w])\le13q+4q\log\Lambda(m)$};
 \node[align=center] at (5,-0.5)
   {$\ell_U(J)=\displaystyle\sum_{\text{cores }C}\ell_U(C)
                 +\displaystyle\sum_{\text{corridors }Q}\ell_U(Q)$};
\end{tikzpicture}
\caption{Direct length estimates on the real core/corridor decomposition.
 Every real preimage of a singular value belongs to an inner core
 (orange), including regular preimages of degree one. The shaded
 source neighborhoods represent components over the larger disks.
 Each remaining corridor (teal) maps monotonically onto a real
 interval avoiding $S(f)$. The covering local isometry and the
 target-interval estimate bound the corridor itself. Adding these
 lengths and the core lengths bounds the length of $J$. The diagram
 is schematic.}
\label{fig:cores}
\end{figure}
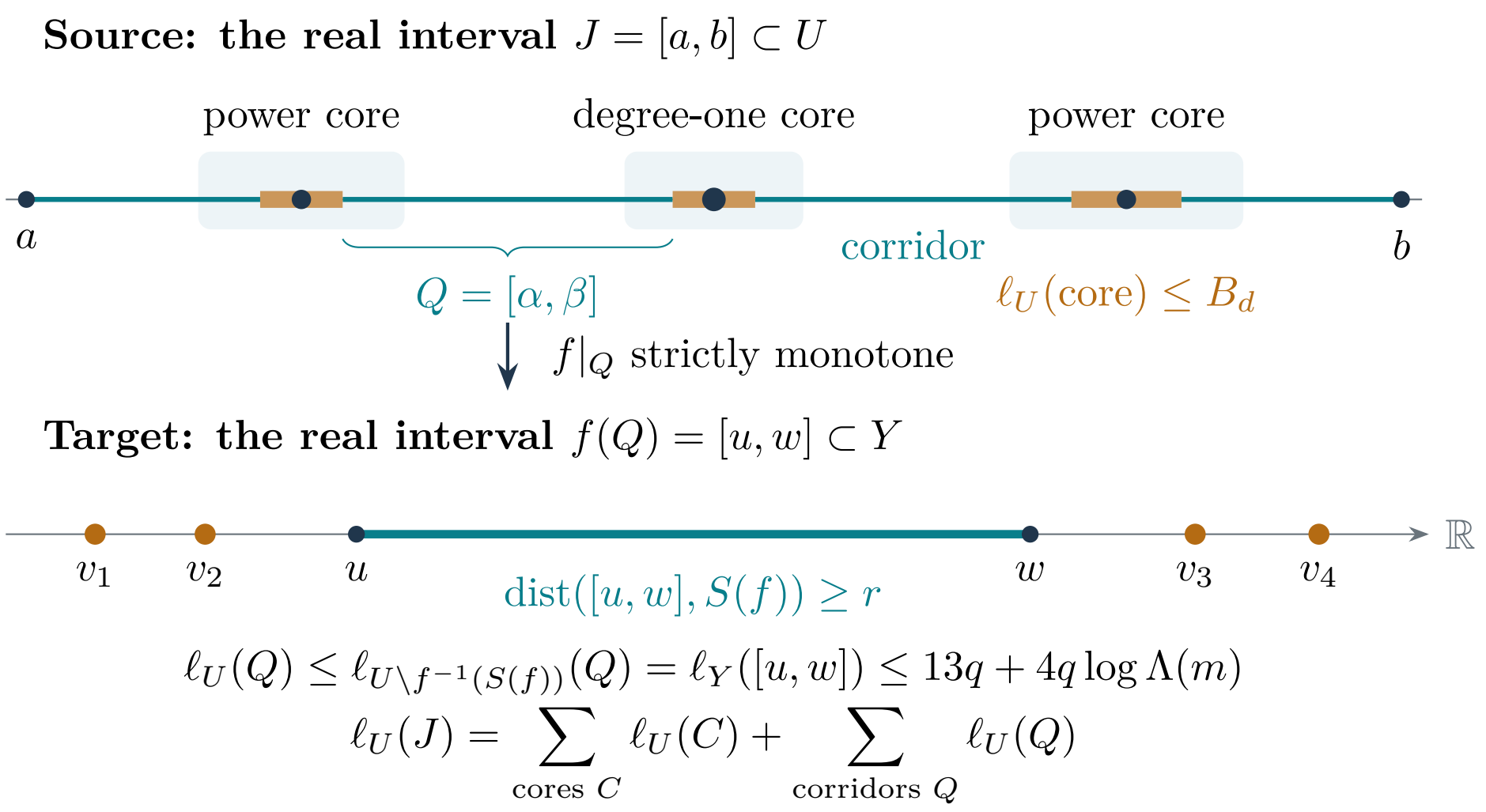

\Cref{fig:cores} shows the decomposition of $J$ into real cores and
corridors, together with the direct target-interval estimate.

\begin{remark}[The role of the endpoint condition]\label{ana:no-tail}
The proof uses no condition on either real tail of $f$ outside $J$.
It also uses no bound on the number of critical points outside $J$
or on the number of complex asymptotic tracts. 
%Endpoint condition~(3)
%in \cref{ana:diameter} is what excludes universal singleton components meeting this finite
%interval. The small radius $r=m/8$ is chosen below the distances of
%the regular endpoint targets from the singular set. In particular, I make no assumptions about successive extrema
%(``record levels'') of the real graph outside $J$, nor about the
%behaviour along either infinite real tail. 
\end{remark}

\FloatBarrier
\section{Affine rescaling estimates}\label{ana:affine}

\begin{proposition}[Uniform control after affine normalization]
\label{ana:scaled}
Under the hypotheses of \cref{ana:diameter}, put
\[
 L=b-a,\qquad F(z)=f(a+Lz),\qquad D=D_{q,d}(m).
\]
For every $t\in J=[a,b]$,
\begin{equation}\label{ana:L-over-radius}
 \frac{L}{R_U(t)}\le\frac{e^{2D}-1}{2}.
\end{equation}
Furthermore,
\begin{align}
 \norm{F'}_{[0,1]}&\le8(e^{2D}-1)\le8e^{2D},
 \label{ana:F-first}\\
 \norm{F''}_{[0,1]}&\le32(e^{2D}-1)^2\le32e^{4D}.
 \label{ana:F-second}
\end{align}
More explicitly, with $M,B_d,\Lambda(m)$ as in
\cref{ana:diameter}, set
\[
 C_D=M B_d+13q(M+1),\qquad A=8q(M+1).
\]
Then
\begin{align}
 \norm{F'}_{[0,1]}
 &\le8e^{2C_D}\Lambda(m)^A,\label{ana:polylog-first}\\
 \norm{F''}_{[0,1]}
 &\le32e^{4C_D}\Lambda(m)^{2A}.
 \label{ana:polylog-second}
\end{align}
In particular, constants $C,B$ depending only on $q,d$ satisfy
\begin{equation}\label{ana:polylog-C2}
 \norm{F'}_{[0,1]}+\norm{F''}_{[0,1]}
 \le C\bigl(1+\log(1/m)\bigr)^B.
\end{equation}
\end{proposition}
\begin{proof}
Both endpoints of $J$ are within hyperbolic distance $D$ of any
$t\in J$. Apply \eqref{ana:euclidean-from-hyperbolic} with center
$t$ to both endpoints. Since $t$ lies between them,
\[
 L=\abs{a-t}+\abs{b-t}
 \le\frac{R_U(t)}2(e^{2D}-1),
\]
proving \eqref{ana:L-over-radius}.

Koebe's quarter theorem gives
$D(t,R_U(t)/4)\subset U$. Since $\abs f<4$ on $U$,
Cauchy's estimates, first on smaller concentric disks and then by
passing to their limiting radius, yield
\begin{equation}\label{ana:raw-cauchy}
 \abs{f'(t)}\le\frac{16}{R_U(t)},\qquad
 \abs{f''(t)}\le\frac{128}{R_U(t)^2}.
\end{equation}
Now $F'(z)=Lf'(a+Lz)$ and
$F''(z)=L^2f''(a+Lz)$. Combining \eqref{ana:raw-cauchy}
with \eqref{ana:L-over-radius} proves
\eqref{ana:F-first} and \eqref{ana:F-second}.

From the definition,
\[
 D=C_D+4q(M+1)\log\Lambda(m).
\]
Exponentiating gives \eqref{ana:polylog-first} and
\eqref{ana:polylog-second}. Since $\Lambda(m)\ge1$ and
$\Lambda(m)\le(1+\log8)(1+\log(1/m))$, the last assertion follows
after adjusting the constant.
\end{proof}

\begin{corollary}[A stretched-exponential scale consequence]
\label{ana:stretch}
Under the same hypotheses, put
\[
 C_{\mathrm{sc}}=\tfrac12e^{2C_D},\qquad A=8q(M+1).
\]
For every $t\in J$,
\begin{equation}\label{ana:stretch-bound}
 m\le8e\,
 \exp\left[-\left(
       \frac{L}{C_{\mathrm{sc}}R_U(t)}
       \right)^{1/A}\right].
\end{equation}
\end{corollary}
\begin{proof}
Equation \eqref{ana:L-over-radius} and the expression for $D$ give
\[
 \frac{L}{R_U(t)}\le
 C_{\mathrm{sc}}\Lambda(m)^A.
\]
Thus $\Lambda(m)\ge
(L/(C_{\mathrm{sc}}R_U(t)))^{1/A}$.
Using the identity $m=8e\exp(-\Lambda(m))$ from  \cref{ana:D-definition} gives the result. 
\end{proof}

\begin{remark}\label{ana:rescaling-meaning}
The estimate controls $L/R_U(t)$, not the absolute length $L$.
This is exactly the ratio that occurs when Cauchy's estimates are
multiplied by the affine rescaling factors $L$ and $L^2$.
No source anchor has been used in this section. Anchors enter only
after marked-gap separation, when a lower bound for an unscaled source
scale is required.
\end{remark}

\section{Expansiveness forces uniform marked gaps}\label{sec:gaps}

In this section the expansiveness will be used, together with the fact 
that each critical point is only of bounded order and the displacement inequality \( \norm{x-y}_\infty\le Km\), to get a lower bound on $m(x)$. 

\begin{lemma}[A bounded target interval and bounded critical multiplicity]
\label{gap:interval}
Let \(J_f=[\widehat y_1,\widehat y_s]\) be the interval spanned by the unscaled inverse images in
\cref{def:admissible}. Then \(f(J_f)\subset[0,1]\). There is a finite
integer \(d\), depending only on the fixed branch and critical labels
in the natural class, such that \(d_f(J_f)\le d\) for all admissible
maps.
\end{lemma}
\begin{proof}
The endpoint images belong to \([0,1]\). If \(f\) took a value larger
than \(1\) on \(J_f\), its maximum on this compact interval would be
attained at an interior critical point and would be a singular value
larger than \(1\). This contradicts \eqref{eq:singular-marking}.
An interior minimum below \(0\) is excluded in the same way.

For the multiplicity assertion, work first in the base map. Between
any two fixed finite labelled critical points there are only finitely
many real critical points, because the zeros of \(f_0'\) are discrete.
Moving a selected point within its fixed open monotonicity branch
does not cross any new critical point. A selected unbounded branch
contains no critical point beyond its finite endpoint. Thus the
interval spanned by all selected inverse points contains a fixed
finite list of possible critical points, including any critical
endpoints. The increasing real natural equivalence preserves their
order and local degrees. The sum of their multiplicities is a
uniform bound \(d\).
\end{proof}

\begin{proposition}[The expansive-gap inequality]\label{gap:inequality}
Let \(F\in C^2([0,1])\), let \(x,y\in W\) satisfy
\(F(y_i)=x_{\pi(i)}\) and \(\norm{x-y}_\infty\le Km\), where
\(m=\min_i(x_{i+1}-x_i)\). Assume
\(F'(y_j)=0\) for \(j\in\cC\), and assume
\eqref{eq:noflat}--\eqref{eq:expansive}. If
\[
 M_1^F=\norm{F'}_{[0,1]},\qquad M_2^F=\norm{F''}_{[0,1]},
\]
then
\begin{equation}\label{eq:gapinequality}
 1\le\frac{(1+2K)^2}{2}\,
      m M_2^F(1+M_1^F)^{2N}.
\end{equation}
\end{proposition}
\begin{proof}
For an interval \(Q=[a,b]\subset[1,s]\) with integer endpoints, define
\(\ell_x(Q)=x_b-x_a\). The extrema of \(G\) on \(Q\) occur at integer
vertices, so \(G(Q)\) also has integer endpoints. The target coordinates
at those extrema are values of \(F\) at marked points between \(y_a\)
and \(y_b\). Hence
\begin{equation}\label{eq:gapstep}
 \ell_x(GQ)\le M_1^F(y_b-y_a)
 \le M_1^F\bigl(\ell_x(Q)+2Km\bigr).
\end{equation}
This is an oscillation estimate and does not require monotonicity
between the marked points.

Choose a shortest edge \(Q_0=[i,i+1]\), so \(\ell_x(Q_0)=m\), and put
\(Q_j=G^j(Q_0)\). Every \(Q_j\) is a nondegenerate integer-endpoint
interval: it contains a unit edge, whose image is nondegenerate by
\eqref{eq:noflat}. Setting \(v_j=\ell_x(Q_j)+2Km\) gives
\[
 v_{j+1}\le M_1^F v_j+2Km\le(1+M_1^F)v_j.
\]
Since \(v_0=(1+2K)m\), it follows that
\begin{equation}\label{eq:gapiteration}
 \ell_x(Q_j)+2Km\le(1+2K)(1+M_1^F)^j m.
\end{equation}

By expansiveness, some \(Q_j\), \(j\le N\), contains a critical label.
Because \(Q_j\) is a nondegenerate integer interval, it contains a
unit edge \([k,k+1]\) with that label as an endpoint. The \(y\)-length
of this edge is at most the right-hand side of
\eqref{eq:gapiteration}. On the other hand, \eqref{eq:noflat} gives
\[
 \abs{F(y_{k+1})-F(y_k)}
 =\abs{x_{\pi(k+1)}-x_{\pi(k)}}\ge m.
\]
Taylor's theorem at the critical endpoint gives
\[
 m\le\frac{M_2^F}{2}(y_{k+1}-y_k)^2
 \le\frac{M_2^F}{2}(1+2K)^2(1+M_1^F)^{2N}m^2.
\]
Dividing by \(m>0\) proves \eqref{eq:gapinequality}.
\end{proof}

\begin{theorem}[Separation under bounded relative displacement]
\label{gap:separation}
For every fixed $K<\infty$ there is $m_*(K)>0$, depending only
on $q$, the uniform critical-multiplicity bound $d$, the fixed
combinatorics, and $K$, such that every admissible configuration
with $\|x-T(x)\|_\infty\le K m(x)$ has
$m(x)\ge m_*(K)$. In particular, if $s\ge3$,
$\operatorname{dist}_\infty(x,\partial W)\ge m_*(K)/2$.
The lower bound is allowed to depend on $K$.
\end{theorem}
\begin{proof}
Every singular value is marked. Hence \(\sep S(f)\ge m\), with the
singleton convention. Each endpoint image of \(J_f\) is a marked
target, so it is either singular or at distance at least \(m\) from
the singular set. Also \(0<m\le1\). By \cref{gap:interval}, all the
hypotheses of the finite-interval estimate \cref{ana:diameter} hold.
The affine rescaling estimate \cref{ana:scaled} consequently gives
fixed constants \(C_1,C_2,A>0\) such that
\[
 M_1^F\le C_1\Lambda^A,\qquad
 M_2^F\le C_2\Lambda^{2A},\qquad
 \Lambda=1+\log(8/m).
\]
Since \(\Lambda\ge1\),
\[
 1+M_1^F\le(1+C_1)\Lambda^A.
\]
Substituting these inequalities into \eqref{eq:gapinequality} gives
\begin{equation}\label{eq:gapcontradiction}
 1\le C_3 m\Lambda^{2A(N+1)}
\end{equation}
for a constant independent of \(f,x,y\). For any fixed \(B\),
\[
 m(1+\log(8/m))^B\longrightarrow0\quad\hbox{as }m\downarrow0.
\]
For example, with \(t=\log(1/m)\), this expression is
\(e^{-t}(1+\log8+t)^B\). Choose \(m_*\in(0,1]\) small enough that
the right-hand side of \eqref{eq:gapcontradiction} is less than \(1\)
for every \(0<m<m_*\). Such \(m\) cannot occur.
\end{proof}

%\begin{remark}[What controls a large affine rescaling]\label{gap:compensation}
%The proof does not assume that the unscaled source length $L$
%from \eqref{eq:affine} is bounded.
%The analytic estimate controls \(L/R_U(t)\), rather than \(L\)
%alone. Its polynomial dependence on \(\log(1/m)\) is equivalent to
%a stretched-exponential implication when this dimensionless ratio
%is large. The expansive-gap inequality excludes the resulting small
%marked gap. Both near-singular target values and pairs of close
%singular values are covered, since all of these targets are marked.
%An endpoint that is singular is handled by a finite power
%core, not by a positive-distance assumption.
%\end{remark}

\subsection{Proof of the Main  \cref{thm:realization-new}}

\begin{corollary}[Boundary displacement and existence]
\label{cor:existence-new}
For the Thurston map associated with the fixed compatible expansive
marking, and every \(K<\infty\), the set
\[
 \{x\in W:\norm{T(x)-x}_\infty\le K m(x)\}
\]
is a compact subset of \(W\). If \(s\ge3\), then
\begin{equation}\label{eq:boundary-displacement-new}
 \frac{\norm{T(x)-x}_\infty}{m(x)}\longrightarrow\infty
 \qquad\text{as }x\longrightarrow\partial W.
\end{equation}
Consequently \(T\) has a fixed point and the prescribed marked
dynamics has a real entire realisation.
\end{corollary}
\begin{proof}
By \cref{gap:separation}, the displayed set is contained in
\(\{x\in W:m(x)\ge m_*\}\), a compact subset of \(W\).
It is closed there by continuity of \(T\) and \(m\), so it is
compact. A bounded subsequence of the quotient in
\eqref{eq:boundary-displacement-new} along a sequence approaching the
boundary would contradict this compactness. By
\cref{lem:mesh-distance}, the displacement divided by the actual
boundary distance also tends to infinity. Apply
\cref{thm:boundary-fixedpoint}. When \(s=2\), the space \(W\)
is a single point and the assertion is immediate. The realisation is
given by \cref{lem:thurston-realisation}.
\end{proof}

The proof for the existence of a fixed point used neither unscaled-map derivative bounds on the
whole family nor any extension of the singular-value parameterisation
through a collision. The stronger unscaled-map derivative and inverse-height
estimates are therefore postponed to the appendices.

\section{The two-peak degeneration and the boundary estimate}
\label{sec:twopeak-degeneration}\label{sec:degeneration-new}

The distinction between compactness of the whole normalized natural
family and the estimates needed for the Thurston map is illustrated by
the family
\begin{equation}\label{eq:twopeak-family}
 f_\lambda(z)=\exp\bigl(-\lambda(z^2-1)^2\bigr),\qquad \lambda\ge1.
\end{equation}
The source normalization fixes the two labelled critical points
\(-1\) and \(1\). The family does degenerate as \(\lambda\to\infty\),
and the present proof does not assert compactness across this
degeneration. This section shows why such a degeneration nevertheless cannot occur along
configurations whose Thurston displacement is bounded by a fixed
multiple of their smallest gap: some terms are exponentially small compared
with other terms that become unbounded; see \cref{deg:exclusion}.

\subsection{The natural class and the actual degeneration}

\begin{lemma}[A normalized family with colliding singular values]
\label{deg:family}
All maps in \eqref{eq:twopeak-family} belong to the same real natural
class, with source equivalences fixing the real line pointwise. Their
real critical points are \(-1,0,1\), and
\begin{equation}\label{eq:twopeak-singular}
 S(f_\lambda)=\{0,e^{-\lambda},1\},\qquad
 f_\lambda(\pm1)=1,\qquad f_\lambda(0)=e^{-\lambda}.
\end{equation}
On the real line,
\begin{equation}\label{eq:twopeak-limit}
 \lim_{\lambda\to\infty}f_\lambda(x)=
 \begin{cases}
  1,&x\in\{-1,1\},\\
  0,&x\notin\{-1,1\}.
 \end{cases}
\end{equation}
The convergence to zero is uniform on every compact subset of
\(\R\setminus\{-1,1\}\). There is no locally uniformly convergent
subsequence with an entire limit on a neighborhood of either peak.
\end{lemma}

\begin{proof}
Put \(P(z)=(z^2-1)^2\). Its critical values are \(0\) and \(1\).
The real-symmetric quasiconformal map
\[
 H_\lambda(w)=\Re w+i\lambda\Im w
\]
is isotopic to the identity through maps fixing \(0\) and \(1\).
Lift this isotopy through the covering associated with \(P\), starting
with the identity lift. The endpoint is a real-symmetric
quasiconformal homeomorphism \(\psi_\lambda\) satisfying
\[
 P\circ\psi_\lambda=H_\lambda\circ P.
\]
The lift extends across the finitely many preimages of \(0\) and
\(1\) by the local power models. It fixes every real point: along
the isotopy, the image of a real point stays in the same discrete
fiber of \(P\), since the target isotopy fixes the real line.

Define the radial quasiconformal homeomorphism
\[
 \phi_\lambda(w)=w\abs{w}^{\lambda-1},\qquad
 \phi_\lambda(0)=0.
\]
For \(u\in\C\),
\(\phi_\lambda(e^{-u})=e^{-\lambda\Re u-i\Im u}\).
Consequently,
\[
 \phi_\lambda\circ f_1\circ\psi_\lambda
 =\exp(-\lambda P)=f_\lambda.
\]
This proves the natural-class assertion and the claimed normalization.

Differentiation gives
\[
 f_\lambda'(z)=-4\lambda z(z^2-1)f_\lambda(z),
\]
so the finite critical values are precisely \(e^{-\lambda}\) and
\(1\). The value \(0\) is an asymptotic value along either real
end. There is no other finite asymptotic value. Indeed, if
\(f_\lambda(\gamma(t))\to a\ne0\) along a curve tending to infinity,
choose a logarithm near \(a\). On a final connected part of the
curve, \(-\lambda P(\gamma(t))\) differs from this logarithm by a
constant integer multiple of \(2\pi i\). It is therefore bounded,
contrary to properness of the polynomial \(P\). This proves
\eqref{eq:twopeak-singular}.

The real limits follow directly from \(P(x)\ge0\), with equality
only at \(x=\pm1\). On a compact real set avoiding these two points,
\(P\) has a positive minimum. Finally, a locally uniform entire
limit near either peak would be continuous and would have the
discontinuous real trace in \eqref{eq:twopeak-limit}, which is
impossible.
\end{proof}

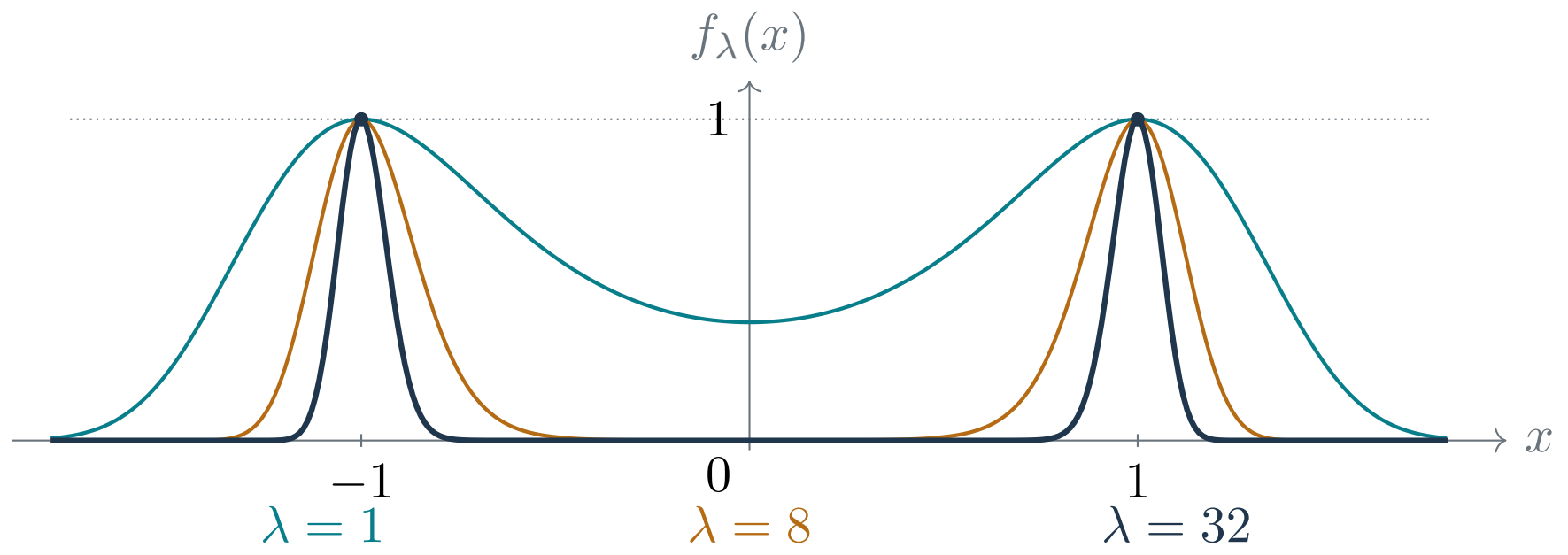
\begin{figure}[tbp]
\centering
\begin{tikzpicture}[x=2.9cm,y=2.4cm]
 \draw[->,softgray] (-1.9,0)--(1.95,0) node[right] {$x$};
 \draw[->,softgray] (0,-0.03)--(0,1.12) node[above] {$f_\lambda(x)$};
 \draw[softgray,densely dotted] (-1.75,1)--(1.75,1);
 \draw[teal,thick] plot[domain=-1.8:1.8,samples=260]
   (\x,{exp(max(-20,-(\x*\x-1)^2))});
 \draw[ochre,thick] plot[domain=-1.8:1.8,samples=300]
   (\x,{exp(max(-20,-8*(\x*\x-1)^2))});
 \draw[ink,very thick] plot[domain=-1.8:1.8,samples=360]
   (\x,{exp(max(-20,-32*(\x*\x-1)^2))});
 \fill[ink] (-1,1) circle (1.5pt) (1,1) circle (1.5pt);
 \foreach \x in {-1,1}{
   \draw[softgray] (\x,0.02)--(\x,-0.02);
   \node[below] at (\x,-0.02) {$\x$};
 }
 \node[left] at (0,1) {$1$};
 \node[below left] at (0,0) {$0$};
 \node[teal] at (-1.1,-0.26) {$\lambda=1$};
 \node[ochre] at (0,-0.26) {$\lambda=8$};
 \node[ink] at (1.1,-0.26) {$\lambda=32$};
\end{tikzpicture}
\caption{The source anchor points stay at \(-1\) and \(1\), while the
critical value at \(0\) collides with the asymptotic value \(0\).
The peaks become narrower. This is a genuine degeneration of the
normalized natural family, but the relative displacement estimate
prevents such configurations from approaching a Thurston fixed point.}
\label{fig:twopeak}
\end{figure}

%The obstruction to using a single complementary domain is visible
%without describing every inverse-image component. For
%\(Z=[0,e^{-1}]\) and \(Y=f_1^{-1}(Z)\),
%\[
% f_1(it)=e^{-(1+t^2)^2}\in Z\qquad(t\in\R),
%\]
%so \(i\R\subset Y\). Hence the anchors \(-1\) and \(1\), which
%lie outside \(Y\), belong to different components of \(\C\setminus Y\).
%A limiting argument applied separately on those components need not
%give one nondegenerate limiting coordinate on both of them.

The convergence in \eqref{eq:twopeak-limit} must not be interpreted
as locally uniform complex convergence on
\(\C\setminus\{-1,1\}\). For example,
\[
 P(1+i/2)=-\frac{15}{16}-\frac i2,
 \qquad
 \abs{f_\lambda(1+i/2)}=e^{15\lambda/16}\longrightarrow\infty.
\]
Thus the real picture describes a discontinuous pointwise limit on
\(\R\), not an entire limit on the punctured plane.

\subsection{Why the degeneration does not obstruct the pullback argument}

There are two distinct assertions in the proof. The finite-interval
estimate \cref{ana:diameter} and its rescaled consequence
\cref{ana:scaled} hold while singular values almost collide. Their constants
are allowed to deteriorate polynomially in \(\log(1/m)\). Only after
combining these estimates with expansiveness and a bound on the
relative displacement does \cref{gap:separation} give a positive
lower bound for \(m\). The uniform bounds for the source-normalized
maps are then proved on that separated set of configurations.

In particular, the proof never needs a compact extension of
\(v\mapsto f_v\) across collisions of singular values. The two-peak
family discussed in this subsection supplies a concrete reason why such an extension would be an
unnecessary and generally false intermediate requirement.

For this example, the compensation can also be seen by direct
calculation, including the potentially large affine rescaling.

\begin{lemma}[Derivative sizes in the two-peak family]
\label{deg:derivatives}
There are absolute constants \(C_1,C_2<\infty\) such that, for
\(\lambda\ge1\),
\begin{equation}\label{eq:twopeak-raw-derivatives}
 \sup_{x\in\R}\abs{f_\lambda'(x)}\le C_1\sqrt\lambda,
 \qquad
 \sup_{x\in\R}\abs{f_\lambda''(x)}\le C_2\lambda.
\end{equation}
Both orders of growth occur: \(f_\lambda''(\pm1)=-8\lambda\), and
\(\sup_\R\abs{f_\lambda'}\ge2\sqrt{2/e}\sqrt\lambda\).
\end{lemma}

\begin{proof}
Besides the first derivative above,
\[
 f_\lambda''(x)=
 \left[-4\lambda(3x^2-1)
       +16\lambda^2x^2(x^2-1)^2\right]f_\lambda(x).
\]
Set \(w=\sqrt\lambda\abs{x^2-1}\). Since \(\lambda\ge1\),
\(x^2\le1+w\). Therefore
\begin{align*}
 \abs{f_\lambda'(x)}
 &\le4\sqrt\lambda\,\sqrt{1+w}\,w e^{-w^2},\\
 \abs{f_\lambda''(x)}
 &\le\lambda\left[4(4+3w)+16(1+w)w^2\right]e^{-w^2}.
\end{align*}
The bracketed functions are bounded for \(w\ge0\), proving the
upper bounds. Substitution of \(x=\pm1\) gives the second-derivative
identity. For the first derivative take
\(x=\sqrt{1+(2\lambda)^{-1/2}}\).
\end{proof}

\begin{proposition}[Exclusion of the two-peak degeneration]
\label{deg:exclusion}
Fix compatible expansive finite marked combinatorics. For any
sequence of configurations with associated source-normalized maps
\(f_{\lambda_n}\), where \(\lambda_n\to\infty\), let
\(y_n=T(x_n)\) and
\(m_n=\min_i(x_{i+1,n}-x_{i,n})\). Then
\begin{equation}\label{eq:twopeak-displacement}
 \frac{\norm{T(x_n)-x_n}_\infty}{m_n}\longrightarrow\infty.
\end{equation}
\end{proposition}

\begin{proof}
Suppress the sequence index. Since \(0,e^{-\lambda},1\) are marked,
\[
 m\le e^{-\lambda},\qquad
 \tau=\log(1/m)\ge\lambda\ge1.
\]
Every finite real pullback target \(v=f_\lambda(\widehat y_i)\)
is positive. As a positive marked coordinate, it satisfies
\(v\ge m\). The explicit inverse equation gives
\[
 (\widehat y_i^2-1)^2
 =\frac{\log(1/v)}{\lambda}
 \le\frac{\tau}{\lambda}.
\]
Consequently, with \(L=\widehat y_s-\widehat y_1\),
\begin{equation}\label{eq:twopeak-span}
 \abs{\widehat y_i}\le
       \sqrt{1+\sqrt{\tau/\lambda}},\qquad
 L\le2\sqrt2\,(\tau/\lambda)^{1/4}.
\end{equation}
The second estimate uses \(\tau/\lambda\ge1\). It allows \(L\)
to grow; no absolute bound on the unscaled source span is being assumed.

For \(F(z)=f_\lambda(\widehat y_1+Lz)\),
\cref{deg:derivatives} and \eqref{eq:twopeak-span} yield constants
\(C_3,C_4\), independent of the configuration, such that
\begin{align}
 \norm{F'}_{[0,1]}
 &\le C_3\lambda^{1/4}\tau^{1/4}
 \le C_3\tau^{1/2},\label{eq:twopeak-scaled-first}\\
 \norm{F''}_{[0,1]}
 &\le C_4\lambda^{1/2}\tau^{1/2}
 \le C_4\tau.\label{eq:twopeak-scaled-second}
\end{align}
If \eqref{eq:twopeak-displacement} failed, some subsequence would
satisfy \(\norm{x-y}_\infty\le Km\) with one finite constant \(K\).
The expansive-gap inequality \cref{gap:inequality} would then give
\[
 1\le C_K\,m\tau^{N+1}
   =C_K\,m\bigl(\log(1/m)\bigr)^{N+1},
\]
where \(N\) depends only on the fixed expansive combinatorics.
Since \(m\le e^{-\lambda}\to0\), the right-hand side tends to
zero, a contradiction.
\end{proof}

The direct computation is a concrete instance of the general
mechanism: derivative growth and the affine source span contribute
only powers of \(\log(1/m)\), while the critical-endpoint estimate
contributes a factor \(m\). This factor makes the contradiction
possible even when the normalized maps themselves have no entire
limit. Arbitrarily small additional regular marked gaps cause no
difficulty; the proof uses \(\tau=\log(1/m)\), not the unjustified
identification \(\tau=\lambda\).

\begin{remark}[The alternative critical-jet normalization]
\label{deg:jet}
The local normalization at a single peak gives a different limit.
Set
\[
 \widetilde f_\lambda(z)
   =f_\lambda\left(1+\frac{z-1}{\sqrt\lambda}\right).
\]
Then the critical point \(1\) and its first nonconstant Taylor
coefficient \(-4\) are fixed, and
\[
 \widetilde f_\lambda(z)=
 \exp\left(-4(z-1)^2
           -\frac{4(z-1)^3}{\sqrt\lambda}
           -\frac{(z-1)^4}{\lambda}\right)
 \longrightarrow e^{-4(z-1)^2}
\]
locally uniformly on \(\C\). The other two critical points move to
\(1-\sqrt\lambda\) and \(1-2\sqrt\lambda\). Thus the behavior of
the source-normalized family depends on the normalization, as it
should. The affine-normalized pullback map and the exclusion in
\cref{deg:exclusion} are unchanged by this positive affine change
of source variable.
\end{remark}

\section{Uniqueness and convergence of the Thurston algorithm}
\label{sec:convergence-new}

The boundary estimate has established existence. To finish the
realisation theorem,  the Teichm\"uller metric on the finite marked
configurations is used. The proof below distinguishes strict decrease of
distance from a uniform Lipschitz constant smaller than one; the latter
is not needed. The standard facts about extremal maps and integrable
quadratic differentials used here are recalled explicitly; see
\cite{Gar87,LV73}.

If \(s=2\), then \(W\) is a point and there is nothing to prove.
Throughout the rest of this section suppose \(s\ge3\). Use the
punctured spheres \((\widehat\C,P_x)\) and the canonical real marking
classes from \cref{subsec:real-teich-chamber}. For \(x,x'\in W\), the
comparison class is represented by a real-symmetric plane
homeomorphism taking \(x_i\) to \(x_i'\), increasing on the real line,
and fixing \(0,1,\infty\); its uniqueness is
\cref{lem:canonical-real-marking}. Define
\begin{equation}\label{eq:teichmetric-new}
 d_{\mathrm T}(x,x')=\frac12\inf_\phi\log K(\phi),
\end{equation}
where the infimum is over quasiconformal maps in that marked class.

\subsection{Extremal maps and compact metric balls}
The extremal-map theorem for finite-type punctured spheres gives a
unique normalised extremal representative. If its dilatation is
\(K>1\), its Beltrami coefficient has the form
\begin{equation}\label{eq:teichcoefficient-new}
 \mu_\phi=k\,\frac{\overline{q}}{|q|},
 \qquad k=\frac{K-1}{K+1},
\end{equation}
where \(q(z)\,\mathrm dz^2\) is a nonzero integrable holomorphic
quadratic differential on \(\widehat\C\setminus P_x\).
Equivalently, it is meromorphic on the sphere with at most simple
poles at the punctures. Reflection preserves the marked class and
dilatation, so uniqueness makes the extremal map real symmetric.
Since it fixes \(0\) and \(1\), it is increasing on \(\R\).
These are applications of the finite-type extremal-map theorem, not
claims about extremal maps on an infinite-type surface.

\begin{lemma}[Topology and properness of the real marked metric]
\label{lem:teich-proper-new}
The metric \(d_{\mathrm T}\) induces the ordinary topology of \(W\),
and its closed balls are compact.
\end{lemma}
\begin{proof}
Fix \(p\in W\). For configurations \(x\) close to \(p\), choose
disjoint small disks around the interior marked points \(p_i\),
disjoint from \(0,1\), and smooth real-symmetric cutoff functions
\(\chi_i\) supported in those disks and equal to one near \(p_i\).
The map
\[
 z\longmapsto z+\sum_{i=2}^{s-1}(x_i-p_i)\chi_i(z)
\]
is a real-symmetric quasiconformal homeomorphism when \(x\) is close
enough to \(p\): its derivative is uniformly close to the identity,
and it is the identity outside a compact set. It sends the marked
points to their prescribed positions, and its dilatation tends to
one as \(x\to p\). Thus Euclidean convergence implies metric
convergence.

Conversely, if \(d_{\mathrm T}(x_n,p)\to0\), take normalised
extremal maps from \(p\) to \(x_n\). Their dilatations tend to one.
Compactness of normalised quasiconformal maps fixing \(0,1,\infty\)
implies that every subsequential limit is a conformal automorphism
fixing these three points, hence the identity. Consequently
\(x_n\to p\) in ordinary coordinates.

For a closed metric ball of radius \(R\) centred at \(p\), the
normalised extremal maps have dilatation at most \(e^{2R}\).
The compactness theorem, applied also to the inverses, gives a
subsequence converging uniformly on the sphere to a quasiconformal
homeomorphism fixing \(0,1,\infty\). Its images of the marked
points are distinct and form a point \(x\in W\). The limit map
is real symmetric and has dilatation at most \(e^{2R}\), so
\(d_{\mathrm T}(p,x)\le R\). Ordinary convergence to \(x\)
implies metric convergence by the first paragraph. The ball is
therefore compact.
\end{proof}

\subsection{Lifting extremal maps}
\begin{lemma}[Nonexpansiveness and its equality case]
\label{lem:teich-lift-new}
For \(x,x'\in W\),
\begin{equation}\label{eq:nonexpansive-new}
 d_{\mathrm T}(T(x),T(x'))\le d_{\mathrm T}(x,x').
\end{equation}
More generally, put \(x_j=T^j(x)\), and define for any integer $n\ge 1$, 
\[
 P_n=F_{x_0}\circ F_{x_1}\circ\cdots\circ F_{x_{n-1}},
 \qquad F_{x_j}=f_{v(x_j)}\circ M_{x_j}.
\]
If \(x\ne x'\) and
\(d_{\mathrm T}(T^n(x),T^n(x'))=d_{\mathrm T}(x,x')\),
then for an extremal target differential \(q\) as in
\eqref{eq:teichcoefficient-new},
\begin{equation}\label{eq:pullback-equality-new}
 P_n^*q=cQ,\qquad c>0,
\end{equation}
where \(Q\) is an integrable quadratic differential with at most
simple poles at \(P_{x_n}\).
\end{lemma}
\begin{proof}
Let \(\phi\) be the extremal map taking the configuration \(x\)
to \(x'\). The singular values of \(F_x\) are mapped to the
corresponding singular values of \(F_{x'}\). Because these entire
maps lie in the same labelled real natural class, the covering
lifting property gives a real-symmetric quasiconformal map
\(\widetilde\phi\) such that
\begin{equation}\label{eq:lift-square-new}
 F_{x'}\circ\widetilde\phi=\phi\circ F_x.
\end{equation}
For completeness, start with a natural equivalence between the two
entire maps. Its target map and \(\phi\) agree on the singular
labels and are isotopic relative to them. Lift that isotopy to adjust
the source map, exactly as in \cref{thm:natural-uniqueness}.
The fixed branch and critical labels select the lift, which is
increasing on \(\R\). It sends each selected preimage \(y_i\)
to \(y_i'\). In particular, it fixes \(0,1,\infty\) and
represents the required marked class from \(T(x)\) to \(T(x')\).

Away from the discrete critical set the two entire maps in
\eqref{eq:lift-square-new} are locally conformal, so
\(K(\widetilde\phi)\le K(\phi)\). This proves
\eqref{eq:nonexpansive-new}. Iterating the construction gives a
normalised lift \(\phi_n\) satisfying
\[
 \widetilde P_n\circ\phi_n=\phi\circ P_n,
 \qquad K(\phi_n)\le K(\phi),
\]
where \(\widetilde P_n\) denotes the composition associated with \(x'\).

If equality of distances holds at the \(n\)-th iterate, this lift
is extremal. Since \(x\ne x'\), its dilatation is greater than
one. The Beltrami pullback formula gives
\[
 \mu_{\phi_n}
 =k\frac{\overline{P_n^*q}}{|P_n^*q|}
 =k\frac{\overline Q}{|Q|}
 \quad\text{almost everywhere},
\]
for an integrable extremal source differential \(Q\). The quotient
\((P_n^*q)/Q\) is meromorphic on \(\C\) and positive real almost
everywhere away from a discrete exceptional set. The open mapping
theorem makes it a positive constant. This proves
\eqref{eq:pullback-equality-new}.
\end{proof}

\subsection{Strict decrease of distance}
\begin{proposition}[A strictly decreasing iterate]
\label{prop:strict-iterate-new}
If \(f_0\) is transcendental, then
\[
 d_{\mathrm T}(T(x),T(x'))<d_{\mathrm T}(x,x')
 \qquad(x\ne x').
\]
If \(f_0\) is a polynomial of degree \(D\ge2\), choose an integer
\(n_0\ge1\) with \(D^{n_0}+2>s\). Then
\[
 d_{\mathrm T}(T^{n_0}(x),T^{n_0}(x'))<d_{\mathrm T}(x,x')
 \qquad(x\ne x').
\]
\end{proposition}
\begin{proof}
First suppose \(f_0\) is transcendental. The same is true of every
\(F_x\), since natural equivalence preserves the infinite covering
degree. If equality held for distinct configurations,
\cref{lem:teich-lift-new} would give \(F_x^*q=cQ\), so the
pullback differential would have finite area.

Choose a small disk \(D_0\) avoiding the finite singular set, all
zeros and poles of \(q\), and the possible exceptional value in
Picard's theorem. A transcendental entire map has infinitely many
preimages of every value in this disk. The covering over \(D_0\)
therefore has infinitely many disjoint univalent sheets. The area of
\(F_x^*q\) on each sheet equals the same positive number
\(\int_{D_0}|q|\). Its total area is infinite, a contradiction.

Now suppose that \(f_0\) is a polynomial of degree \(D\). The composition \(P_{n_0}\)
in \cref{lem:teich-lift-new} has degree \(d=D^{n_0}\).
A nonzero meromorphic quadratic differential on the sphere has divisor
degree \(-4\). Since every allowed pole of \(q\) is simple,
including the possible pole at infinity, \(q\) has at least three
distinct finite simple poles. Let their number be \(r\ge3\).

There are \(rd\) preimages of these poles counted with local
degree. If \(e_z\ge2\) is a local degree at a ramified finite
preimage, then \(e_z\le2(e_z-1)\). The total finite critical
multiplicity of a degree-\(d\) polynomial is \(d-1\).
Consequently the number of unramified preimages of the finite poles
is at least
\[
 rd-2(d-1)\ge d+2>s.
\]
Each is a simple pole of \(P_{n_0}^*q\). But an integrable source
differential \(Q\) may have finite poles only at the \(s\)
marked source points. Thus equality in
\eqref{eq:pullback-equality-new} is impossible. This proves the
second assertion.
\end{proof}

\begin{theorem}[Uniqueness and global convergence]
\label{thm:convergence-new}
The fixed point of the Thurston map is unique, and \(T^n(x)\)
converges to it for every \(x\in W\). The dynamically normalised
entire map realising the prescribed marked data is unique within the fixed natural family \(\mathcal N(f_0)\) with the prescribed labels. 
\end{theorem}
\begin{proof}
Existence was established in \cref{cor:existence-new}; denote a
fixed point by \(p\). Take \(n_0=1\) in the transcendental case
and take \(n_0\) as in \cref{prop:strict-iterate-new} in the
polynomial case. Two distinct fixed points of \(T\) would have
their distance strictly decreased by \(T^{n_0}\), which is
impossible. Thus \(p\) is unique.

For any \(x\), nonexpansiveness makes
\(r_n=d_{\mathrm T}(T^n(x),p)\) nonincreasing, with limit
\(r_\infty\ge0\). The orbit lies in a compact metric ball by
\cref{lem:teich-proper-new}. Choose a subsequence
\(T^{n_j}(x)\to a\). Continuity gives
\[
 d_{\mathrm T}(a,p)=r_\infty,
 \qquad
 d_{\mathrm T}(T^{n_0}(a),p)=r_\infty.
\]
Strict decrease under \(T^{n_0}\) forces \(a=p\), and therefore
\(r_\infty=0\). The whole orbit converges to \(p\).

Finally, suppose an entire map \(F\) with the prescribed labelled
dynamics has a marked configuration \(x\) with extreme points
\(0,1\). Source-normalise it by a positive affine map \(A\),
so that \(f_{v(x)}=F\circ A\), using
\cref{thm:natural-uniqueness}. Its unscaled inverse images are
\(\widehat y_i=A^{-1}(x_i)\), hence its affine normalising map
is \(M_x=A^{-1}\). It follows that \(T(x)=x\) and
\(F=F_x\). Thus \(x=p\) and \(F=F_p\). Removing the endpoint
normalisation gives uniqueness up to positive real affine conjugacy.
This completes the proof of \cref{thm:realization-new}.
\end{proof}

\appendix

\part*{Appendices: additional analytic estimates}
\addcontentsline{toc}{part}{Appendices: additional analytic estimates}

\section{Additional a priori estimates: real ends and a finite anchored core}\label{glob:core-section}

The realisation theorem is already complete before this appendix: the
uniform marked-gap separation proved in \cref{thm:apriori,gap:separation}
gives the boundary displacement, while uniqueness and convergence are
proved in \cref{sec:convergence-new}. The material in
\cref{glob:core-section,glob:derivative-section,glob:height-section}
records stronger analytic consequences which may be useful independently,
but none of them is used in the proof of \cref{thm:realization-new}.

\begin{theorem}[Additional derivative and inverse-height bounds]
\label{thm:additional-apriori}
Assume the hypotheses of \cref{thm:apriori} and one of the two source
normalizations from \cref{def:source-normalization-new}: either two fixed
labelled real critical points, or one fixed labelled real critical point
together with the first nonconstant Taylor coefficient in \eqref{eq:jet}.
Then there are constants \(M_1,M_2<\infty\) such that
\[
 |f'(t)|\le M_1,\qquad |f''(t)|\le M_2
\]
for every \(t\) in the real closures of the selected branches with
\(f(t)\in[-1,2]\). Every selected inverse branch \(g_{I,f}\) continues
univalently to the upper half-plane and satisfies
\[
 \left|\Im g_{I,f}\left(\frac{1+i}{2}\right)\right|
 \ge h_*>0,
 \qquad
 h_* =\frac{\eta^2}{12M_1},\quad
 \eta=\min\{m_*/2,1/2\}.
\]
The constants may depend on the source normalization and the fixed
labelled branches. No assumption is made on the total number of real
critical points or asymptotic tracts.
\end{theorem}

The separation result has been proved without any source normalization.
Write $m_*>0$ for the lower bound obtained in \cref{gap:separation}.
In particular,
\begin{equation}\label{glob:singular-separation}
 \sep S(f)\ge m_*,
\end{equation}
with $\sep S(f)=+\infty$ when the singular set is a singleton. I
retain the fixed finite collection of labelled branches from the
theorem and impose one of its two source normalizations.

\begin{lemma}[Classification of the real ends]\label{glob:ends}
At either end of the real line, one of the following descriptions applies:
\begin{enumerate}
 \item There are infinitely many real critical points in that direction.
 After the first such critical point, the image of that tail is contained in the convex
 hull of $S(f)$.
 \item There is a terminal monotone branch whose limit is a finite
 asymptotic value in $S(f)$.
 \item There is a terminal monotone branch whose limit is $+\infty$
 or $-\infty$.
\end{enumerate}
The alternatives at the two ends need not agree. Every finite endpoint
of the image of a real branch is a singular value.
\end{lemma}
\begin{proof}
The zeros of $f'$ are discrete, since $f$ is nonconstant. If there are
infinitely many real critical points towards $+\infty$, enumerate
consecutive ones as $c_0<c_1<\cdots$, with $c_j\to+\infty$.
The derivative has constant nonzero sign on $(c_j,c_{j+1})$, so
$f$ is strictly monotone there. Both endpoint values are singular,
and hence
\[
 f([c_j,c_{j+1}])\subset[\min S(f),\max S(f)].
\]
These intervals cover $[c_0,\infty)$.
If there are only finitely many critical points towards the end, the
function is strictly monotone on a terminal interval. Its limit
exists in the extended real line. A finite limit is an asymptotic
value by definition. The argument at $-\infty$ is identical.

A finite source endpoint of a maximal monotonicity branch is critical.
An infinite source endpoint with a finite image limit contributes an
asymptotic value. These exhaust the finite endpoints of the branch image.
\end{proof}

\begin{remark}\label{glob:no-whole-line}
Expansiveness in this theorem involves an actual real critical label.
Thus $\Crit_{\R}(f)\ne\varnothing$, and no selected branch is all of
$\R$. Every selected unbounded branch is a half-ray with a finite
critical endpoint. Stationary critical points, of odd local degree
at least three, are allowed throughout the argument.
\end{remark}

\begin{lemma}[A finite critical budget]\label{glob:budget}
Fix finitely many labelled real branches and labelled real critical
points. There is an integer $d_*\ge0$ such that the convex hull of
those critical points and arbitrary finite points on the selected
branches contains total real critical multiplicity at most $d_*$.
Finitely many points may be chosen on each branch.
\end{lemma}
\begin{proof}
In the base map choose a compact interval containing every finite
endpoint of a selected branch and every additional labelled critical
point. Every selected bounded branch is contained in this interval.
A selected unbounded branch extends beyond it only along a
critical-free half-ray. Thus moving chosen points along these branches
encounters no additional critical points beyond this fixed finite list.
The increasing source homeomorphism preserves the order, labels,
and local degrees of critical points, giving the same multiplicity
bound throughout the family.
\end{proof}

Set
\begin{equation}\label{glob:rho}
 \rho=\min\{m_*/16,\,1/16\}.
\end{equation}
For each selected branch make these finite choices:
\begin{enumerate}
 \item If both source endpoints are finite, choose both endpoints.
 \item If the branch tends to $+\infty$ in value, choose its finite
 critical endpoint and the unique point where $f=2$.
 \item If it tends to $-\infty$ in value, choose its finite critical
 endpoint and the unique point where $f=-1$.
 \item If it tends to a finite asymptotic value $v$, choose its finite
 critical endpoint and the unique entry point $t_{\rm ent}$ with
 \begin{equation}\label{glob:entry}
  f(t_{\rm ent})=v+\sigma\rho/e,
 \end{equation}
 where $\sigma\in\{-1,1\}$ is the sign of $f-v$ on that branch.
\end{enumerate}
The last choice exists because the two finite endpoints of its image
are distinct singular values, separated by at least $m_*$. The
cutoffs at $-1$ and $2$ exist because the finite endpoint value
belongs to $[0,1]$. Let $E_f$ be the convex hull of these points
and the source anchor points. The source anchor points are labelled critical points.

\begin{proposition}[Uniform diameter of the finite core]\label{glob:core}
Let $V_f$ be the component of $f^{-1}D(-1,12)$ containing $E_f$.
There is a uniform constant $D_*<\infty$ such that
\begin{equation}\label{glob:core-diameter}
 \diam_{V_f} E_f\le D_*.
\end{equation}
Moreover,
\begin{equation}\label{glob:core-image}
 f(E_f)\subset[-1,2].
\end{equation}
\end{proposition}
\begin{proof}
Every point defining the convex hull has image in $[-1,2]$.
Every real critical value belongs to $[0,1]$. An excursion beyond
$[-1,2]$ in the convex hull would give an interior extremum with
a real critical value outside $[0,1]$. This proves
\eqref{glob:core-image}.

Each endpoint image of $E_f$ is singular, is $-1$ or $2$, or is
an entry value. At an entry value its distance from $v$ is $\rho/e$,
and its distance from every other singular value is at least
$m_*-\rho/e>\rho/e$. Thus every regular endpoint image is at
distance at least $\rho/e$ from $S(f)$.

Apply \cref{ana:diameter} to
\[
 H_f=\frac{f+1}{3},\qquad \mu=\frac{\rho}{3e}.
\]
Then $H_f(E_f)\subset[0,1]$,
$S(H_f)=(S(f)+1)/3\subset[1/3,2/3]$, and the singular separation
is at least $m_*/3\ge\mu$. Each endpoint image is singular or is
at distance at least $\mu$ from $S(H_f)$. Lemma~\ref{glob:budget}
bounds the critical multiplicity in $E_f$ by $d_*$.
The relevant component of $H_f^{-1}D(0,4)$ is exactly $V_f$.
If $E_f$ is a singleton, its diameter is zero. Otherwise the
diameter lemma with this fixed positive $\mu$ proves the result.
\end{proof}

\begin{figure}[htbp]
\centering
\begin{tikzpicture}[x=1.02cm,y=0.75cm,>=Stealth,font=\small]
 \draw[<->,softgray] (-0.25,0.2)--(10.9,0.2);
 \draw[line width=4pt,teal!45] (2.8,0.2)--(9,0.2);
 \draw[densely dashed,ochre] (-0.05,3.2)--(3.6,3.2);
 \node[ochre,above left] at (1.3,3.22) {asymptotic value \(v\)};
 \draw[very thick,teal] (0,3.1)
   .. controls (1.4,3.05) and (2.1,2.95) .. (2.8,2.7);
 \draw[very thick,ink] (2.8,2.7)
   .. controls (4.0,2.0) and (4.2,0.7) .. (5,0.7)
   .. controls (6,0.7) and (7.4,2.3) .. (9,4.2);
 \draw[very thick,softgray,->] (9,4.2)
   .. controls (9.5,4.9) and (10,5.1) .. (10.4,5.6);
 \draw[densely dashed,softgray] (2.8,0.2)--(2.8,2.7);
 \draw[densely dashed,softgray] (5,0.2)--(5,0.7);
 \draw[densely dashed,softgray] (9,0.2)--(9,4.2);
 \draw[densely dashed,ochre] (8.1,4.2)--(10.2,4.2);
 \node[ochre,right] at (10.2,4.2) {\(2\)};
 \foreach \p in {(2.8,2.7),(5,0.7),(9,4.2)}
   \fill[ink] \p circle (2.4pt);
 \foreach \x in {2.8,5,9}\fill[ink] (\x,0.2) circle (2pt);
 \node[below] at (2.8,0.03) {$t_{\mathrm{ent}}$};
 \node[below] at (5,0.03) {$c$};
 \node[below] at (9,0.03) {$t_{\mathrm{cut}}$};
 \node[anchor=west,fill=white,inner sep=2pt] at (3.1,2.85)
   {\(f(t_{\mathrm{ent}})=v-\rho/e\)};
 \node[above] at (5,1.55) {\(f(c)\in S(f)\)};
 \node[teal,align=center] at (1.2,-0.85) {logarithmic tail};
 \node[teal] at (6.1,-0.85) {finite core \(E_f\)};
 \node[softgray] at (10,-0.85) {\(f>2\)};
\end{tikzpicture}
\caption{A mixed configuration, drawn schematically. The finite core
 joins the source anchor points to an entry point on a finite-asymptotic
 half-ray and to a fixed target cutoff on a branch escaping in value.
 Bounded hyperbolic diameter and the anchors control derivatives on
 the core. A logarithmic coordinate controls the remaining left tail.
 The part beyond the right cutoff is outside the target window for
 the derivative assertion.}
\label{fig:tails}
\end{figure}
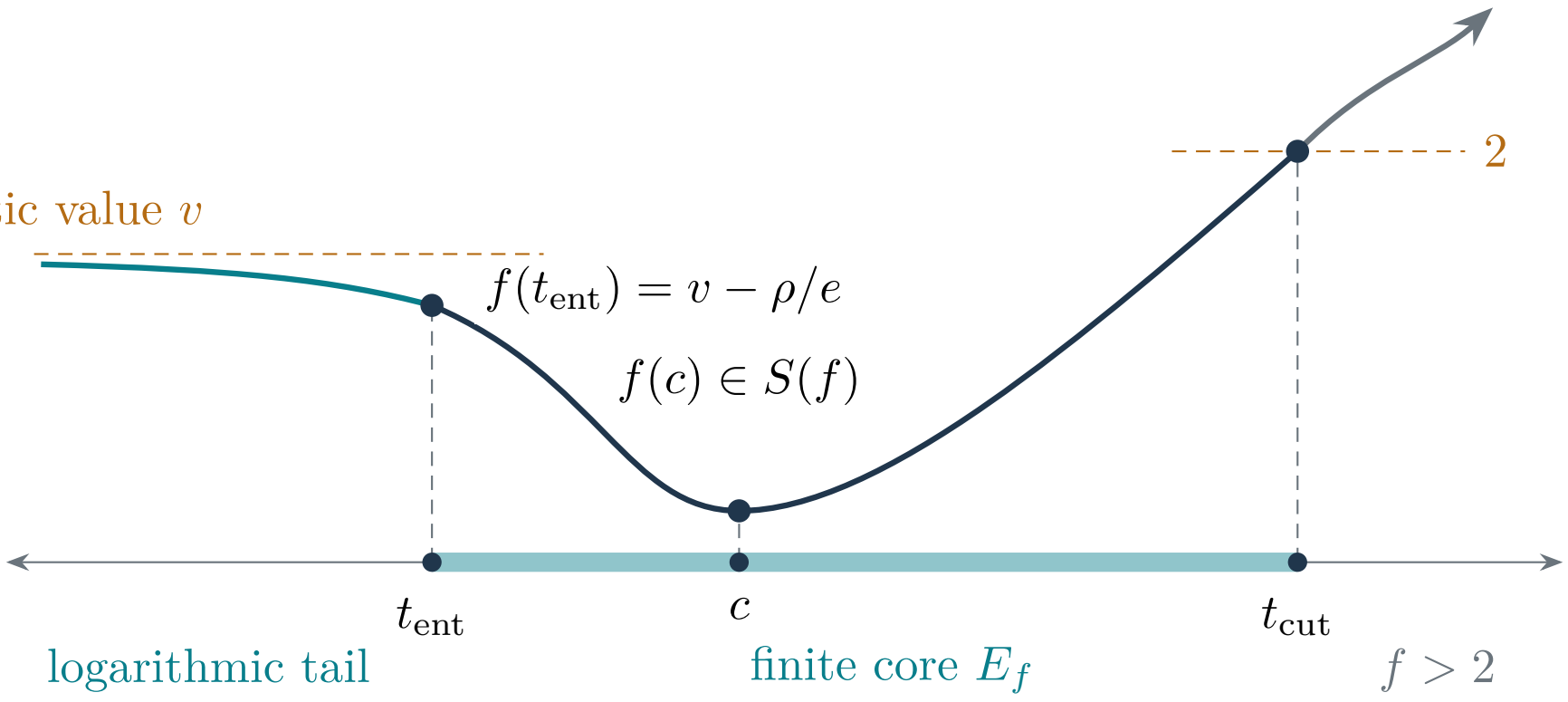

\Cref{fig:tails} illustrates a finite core with two different end types.

\section{Source normalization and real derivative bounds}
\label{glob:derivative-section}

The source normalizations in \cref{def:source-normalization-new} now provide a lower conformal-radius bound on the separated finite core.

Write $R_U(t)$ for the conformal radius of $U$ at $t$, as in
\cref{ana:radius}.

\begin{proposition}[The anchors give a lower source scale]
\label{glob:radius-lower}
Under either source normalization in Theorem~\ref{thm:additional-apriori}, there is
$r_*>0$ such that
\begin{equation}\label{glob:radius-lower-equation}
 R_{V_f}(t)\ge r_*\qquad(t\in E_f).
\end{equation}
One may take
\begin{equation}\label{glob:rstar}
 r_*=
 \begin{cases}
  4(c_2-c_1)e^{-4D_*},
  &\text{for two fixed labelled critical points},\\[1mm]
  (\rho/|A_\ell|)^{1/\ell}e^{-2D_*},
  &\text{for a fixed critical point and coefficient $A_\ell$}.
 \end{cases}
\end{equation}
\end{proposition}
\begin{proof}
For two critical anchors, Koebe growth for a Riemann map centered
at $c_1$ gives
\[
 c_2-c_1\le
 \frac{R_{V_f}(c_1)}4(e^{2D_*}-1)
 \le\frac{R_{V_f}(c_1)}4e^{2D_*}.
\]
Hence $R_{V_f}(c_1)\ge4(c_2-c_1)e^{-2D_*}$.

For a fixed critical jet, let $T_f$ be the component of
$f^{-1}D(f(c),\rho)$ containing $c$. This target disk contains
just the singular value $f(c)$. Because the component contains a
preimage of the center, Lemma~\ref{ana:singleton} makes it a finite
power component of degree $\ell$. A symmetric Riemann map
$h:\D\to T_f$, with $h(0)=c$ and $h'(0)>0$, satisfies
\[
 f(h(z))=f(c)+\sigma\rho z^\ell,\qquad \sigma\in\{-1,1\}.
\]
Indeed, the proper disk map has its only zero at the origin, with
multiplicity $\ell$. Comparing its first nonconstant coefficient
with the prescribed critical expansion gives
\begin{equation}\label{glob:power-radius}
 R_{T_f}(c)=h'(0)=\left(\frac{\rho}{|A_\ell|}\right)^{1/\ell}.
\end{equation}
The disk $D(f(c),\rho)$ lies in $D(-1,12)$, and $T_f$ meets $V_f$
at $c$, so $T_f\subset V_f$. Monotonicity of conformal radius
gives the same lower bound for $R_{V_f}(c)$.

To propagate either bound, take a Riemann map $h:\D\to V_f$
centered at the chosen anchor $c$, and write $t=h(z)$.
Koebe's derivative estimate yields
\[
 R_{V_f}(t)=(1-|z|^2)|h'(z)|
 \ge R_{V_f}(c)
       \left(\frac{1-|z|}{1+|z|}\right)^2
 =R_{V_f}(c)e^{-2d_{V_f}(c,t)}.
\]
Use \eqref{glob:core-diameter}. In the two-anchor case use $c=c_1$.
\end{proof}

\begin{remark}\label{glob:jet-direction}
The finite power component is essential for the direction of
\eqref{glob:power-radius}. A fixed nonzero Taylor coefficient and
Cauchy's inequality alone give an upper conformal-radius bound,
whereas the derivative estimates require a lower bound.
\end{remark}

\begin{corollary}[Derivatives on the finite core]\label{glob:core-derivatives}
For every $t\in E_f$,
\begin{equation}\label{glob:core-derivative-constants}
 |f'(t)|\le M_{1,\mathrm{core}}:=\frac{48}{r_*},\qquad
 |f''(t)|\le M_{2,\mathrm{core}}:=\frac{384}{r_*^2}.
\end{equation}
\end{corollary}
\begin{proof}
Koebe's quarter theorem puts $D(t,R_{V_f}(t)/4)$ inside $V_f$.
On this domain $|f+1|<12$. Cauchy's estimates on slightly smaller
disks, followed by a limit, give
\[
 |f'(t)|\le48/R_{V_f}(t),\qquad
 |f''(t)|\le384/R_{V_f}(t)^2.
\]
Apply Proposition~\ref{glob:radius-lower}.
\end{proof}

\subsection{The complete finite-asymptotic half-rays}
\label{glob:log-tail-section}

The core contains every selected bounded branch and every part of a
value-escaping branch over $[-1,2]$. It remains to treat a half-ray
with a finite asymptotic limit.

\begin{lemma}[A normalized logarithmic inverse]\label{glob:log-inverse}
Let a selected half-ray tend to $v\in S(f)$, and let $t_{\rm ent}$
be its entry point. There is a univalent map
$h:\{z:\Re z>0\}\to\C$ satisfying
\begin{equation}\label{glob:log-identity}
 f(h(z))=v+\sigma\rho e^{-z},\qquad h(1)=t_{\rm ent}.
\end{equation}
The map $h$ is real on the positive real axis, and $h([1,\infty))$
is exactly the remainder of the half-ray beyond $t_{\rm ent}$.
\end{lemma}
\begin{proof}
The remainder lies in one component $T$ of $f^{-1}D(v,\rho)$.
If $f:T\to D(v,\rho)$ were finite and proper, the preimage of a
small closed disk about $v$ would be compact in $T$. It could not
contain an unbounded real ray whose image tends to $v$.
Lemma~\ref{ana:singleton} therefore makes
$f:T\to D(v,\rho)\setminus\{v\}$ a universal covering.

The map $z\mapsto v+\sigma\rho e^{-z}$ is another universal covering
of the same punctured disk, from the right half-plane. Equivalence
of universal coverings gives the conformal map $h$ with its
specified value at $1$. The conjugate map
$z\mapsto\overline{h(\bar z)}$ is another lift with the same value
there. Uniqueness of lifts makes it equal to $h$, so $h(t)$ is
real for $t>0$.

For $t\ge1$, the image in \eqref{glob:log-identity} runs from the
entry value to $v$. Each value has a unique inverse on the selected
branch, and the lift starts on that branch. Thus it traces exactly
the remaining half-ray. A finite limit of $h(t)$ as $t\to\infty$
would be a preimage of $v$ in the same open component $T$,
contrary to the omitted-value alternative.
\end{proof}

\begin{lemma}[Distortion in the logarithmic coordinate]
\label{glob:log-distortion}
Put $a=|h'(1)|$. For $t\ge1$,
\begin{equation}\label{glob:log-distortion-equations}
 a t^{-3}\le |h'(t)|\le at,\qquad
 \left|\frac{h''(t)}{h'(t)}\right|\le\frac4t.
\end{equation}
Moreover,
\begin{equation}\label{glob:a-lower}
 a\ge a_*:=\frac{\rho}{eM_{1,\mathrm{core}}}>0.
\end{equation}
\end{lemma}
\begin{proof}
Set $\chi(\zeta)=(1+\zeta)/(1-\zeta)$ and $H=h\circ\chi$.
Then $H$ is univalent on $\D$ and $|H'(0)|=2a$.
For $r=(t-1)/(t+1)\in[0,1)$, Koebe's derivative bounds give
\[
 2a\frac{1-r}{(1+r)^3}\le |H'(r)|
 \le 2a\frac{1+r}{(1-r)^3}.
\]
Divide by $\chi'(r)=2/(1-r)^2$ to obtain
\[
 a\left(\frac{1-r}{1+r}\right)^3\le |h'(t)|
 \le a\frac{1+r}{1-r}.
\]
These are the claimed bounds for $h'$.

The disk $D(t,t)$ is contained in the right half-plane. The function
\[
 k(\zeta)=\frac{h(t+t\zeta)-h(t)}{t h'(t)}
\]
is normalized and univalent on $\D$. The classical second-coefficient
estimate $|k''(0)|\le4$ gives
$t|h''(t)/h'(t)|\le4$; see \cite{Duren}.

Finally, differentiation of \eqref{glob:log-identity} at $1$ gives
$|f'(t_{\rm ent})|a=\rho/e$. The entry value is regular, so this
derivative is nonzero. Since $t_{\rm ent}\in E_f$,
Corollary~\ref{glob:core-derivatives} proves \eqref{glob:a-lower}.
\end{proof}

\begin{proposition}[Derivatives on all selected branches]
\label{glob:all-derivatives}
There are uniform constants $M_1,M_2<\infty$ such that
$|f'(t)|\le M_1$ and $|f''(t)|\le M_2$ on every selected branch
wherever $f(t)\in[-1,2]$. If the branch image is bounded, these
estimates hold on its entire length, including an infinite source
tail when present.
\end{proposition}
\begin{proof}
Differentiate \eqref{glob:log-identity} along $t\ge1$:
\[
 f'(h(t))h'(t)=-\sigma\rho e^{-t},\qquad
 f''(h(t))h'(t)^2+f'(h(t))h''(t)=\sigma\rho e^{-t}.
\]
It follows that
\begin{align}
 |f'(h(t))|&\le\frac{\rho}{a}t^3e^{-t},
 \label{glob:tail-first}\\
 |f''(h(t))|&\le\frac{\rho}{a^2}(t^6+4t^5)e^{-t}.
 \label{glob:tail-second}
\end{align}
For the second inequality, use the exact identity
\[
 f''(h(t))=\frac{\sigma\rho e^{-t}}{h'(t)^2}
          \left(1+\frac{h''(t)}{h'(t)}\right).
\]
The lower bound $a\ge a_*$ is uniform. For example, set
\[
 A_3=27e^{-3},\qquad A_6=6^6e^{-6}+4\cdot5^5e^{-5}.
\]
Then one can take
\begin{equation}\label{glob:final-derivative-constants}
 \begin{split}
 M_1&=\max\{M_{1,\mathrm{core}},\,\rho A_3/a_*\},\\
 M_2&=\max\{M_{2,\mathrm{core}},\,\rho A_6/a_*^2\}.
 \end{split}
\end{equation}
These bounds cover the core and every remaining finite-asymptotic
tail. The bounded branches lie in the core. On a value-escaping
branch, the part over $[-1,2]$ lies between its finite endpoint
and its cutoff at $-1$ or $2$, hence also lies in the core.
Lemma~\ref{glob:ends} exhausts the required points.
\end{proof}

\begin{remark}\label{glob:target-window}
A bounded target window is necessary on a value-escaping branch:
even the fixed map $f(t)=t^2$ has unbounded derivative on its whole
escaping half-ray. Any fixed compact target interval may replace
$[-1,2]$, with adjusted constants.
\end{remark}

\section{Uniform inverse heights}\label{glob:height-section}

\subsection{Inverse branches on the upper half-plane}

\begin{lemma}\label{glob:inverse-extension}
For a real monotonicity branch $I$ and $J=f(I)$, the real inverse
extends to a univalent map $g:\HH\to\C$ with
$f\circ g=\operatorname{id}$. Its image lies in one open half-plane.
It is real and holomorphic across every point of $J$.
\end{lemma}
\begin{proof}
All singular values are real. Therefore $f$ restricts over $\HH$
to a covering, and each component of $f^{-1}(\HH)$ maps conformally
onto the simply connected domain $\HH$. A local inverse at a
regular real point selects one such component and its univalent
inverse $g$. No point of $g(\HH)$ is real, since $f(\R)\subset\R$.
Connectedness puts its image entirely in one half-plane.

The inverse-function theorem supplies real holomorphic inverse
germs along all of $J$. Along any compact subinterval, finitely
many overlapping germs connect to the initial one. Their
restrictions to $\HH$ agree with $g$ by uniqueness of analytic
continuation. Hence the real extension holds across all of $J$.
Singular values caused by other preimages or tracts may lie inside
$J$; they do not obstruct this regular inverse branch.
\end{proof}

\subsection{A boundary derivative controls an interior height}

\begin{lemma}[Positive harmonic representation]\label{glob:poisson}
Every positive harmonic function $u$ on $\HH$ has a representation
\begin{equation}\label{glob:poisson-equation}
 u(x+iy)=\beta y+
 \int_{\R}\frac{y}{(x-t)^2+y^2}\,\dd\mu(t),\qquad \beta\ge0,
\end{equation}
where $\mu$ is positive and
$\int_{\R}(1+t^2)^{-1}\dd\mu(t)<\infty$.
If $u$ extends continuously with value zero across an open real
interval $J$, then $\mu(J)=0$.
\end{lemma}
\begin{proof}
Map $\HH$ to $\D$ by $\chi(w)=(w-i)/(w+i)$ and put
$v=u\circ\chi^{-1}$. For $0<r<1$, let $\nu_r$ on the unit circle
have density $v(re^{i\theta})\,\dd\theta/(2\pi)$.
The mean-value property gives total mass $v(0)$ independently
of $r$. Choose $r_n\uparrow1$ along which these positive measures
converge weakly to a finite measure $\nu$.

For fixed $z\in\D$ and $r_n>|z|$, Poisson's formula gives
\[
 v(z)=\int_{\partial\D}
 \frac{1-|z/r_n|^2}{|\zeta-z/r_n|^2}\,\dd\nu_{r_n}(\zeta).
\]
The kernels converge uniformly, so
\[
 v(z)=\int_{\partial\D}
 \frac{1-|z|^2}{|\zeta-z|^2}\,\dd\nu(\zeta).
\]
The point $\zeta=1$ corresponds to infinity; its Poisson kernel
at $\chi(x+iy)$ is $y$. For a finite real $t$,
\[
 \frac{1-|\chi(x+iy)|^2}{|\chi(t)-\chi(x+iy)|^2}
 =\frac{y(1+t^2)}{(x-t)^2+y^2}.
\]
Take $\beta=\nu(\{1\})$, push the rest of $\nu$ to $\R$, and
multiply it by $1+t^2$. This proves the representation and its
integrability condition.

If $u$ has zero continuous boundary values on $J$, choose a
nonnegative continuous function $\varphi$ compactly supported in
$J$. Then $\int\varphi(x)u(x+iy)\,\dd x\to0$ as $y\downarrow0$.
Tonelli's theorem and Fatou's lemma give the lower bound
\[
 \int_{\R}\liminf_{y\downarrow0}
 \left(\int_{\R}\frac{y\varphi(x)}{(x-t)^2+y^2}\,\dd x\right)
 \dd\mu(t)
 =\pi\int_{\R}\varphi(t)\,\dd\mu(t)
\]
for this limit. The inner kernels are the unnormalized Poisson
approximate identity. Thus $\int\varphi\,\dd\mu=0$ for every such
test function, proving $\mu(J)=0$.
\end{proof}

\begin{lemma}[Pick height estimate]\label{glob:pick}
Let $G:\HH\to\HH$ be holomorphic with a real-valued holomorphic
extension across an open real interval $J$. If
\[
 u\in[-1,2],\qquad 0<\eta\le1/2,\qquad
 (u-\eta,u+\eta)\subset J,
\]
then
\begin{equation}\label{glob:pick-equation}
 \operatorname{Im}G\left(\frac{1+i}{2}\right)
 \ge\frac{\eta^2}{12}G'(u).
\end{equation}
\end{lemma}
\begin{proof}
Apply Lemma~\ref{glob:poisson} to $\operatorname{Im}G$.
The representing measure vanishes on $J$. Divide the formula
at $u+iy$ by $y$, then use monotone convergence and the real
holomorphic extension to get
\begin{equation}\label{glob:pick-boundary-derivative}
 G'(u)=\beta+\int_{\R}\frac{\dd\mu(t)}{(t-u)^2}.
\end{equation}
On the support of $\mu$, $|t-u|\ge\eta$ and $|u-1/2|\le3/2$.
Consequently
\begin{align*}
 (t-1/2)^2+1/4
 &\le2(t-u)^2+2(u-1/2)^2+1/4\\
 &\le2(t-u)^2+19/4\\
 &\le\frac6{\eta^2}(t-u)^2.
\end{align*}
The last inequality uses $\eta\le1/2$. Also
$\beta/2\ge(\eta^2/12)\beta$. Evaluating the representation at
$(1+i)/2$ proves
\[
 \operatorname{Im}G\left(\frac{1+i}{2}\right)
 =\frac{\beta}{2}
 +\int_{\R}\frac{(1/2)\,\dd\mu(t)}{(t-1/2)^2+1/4}
 \ge\frac{\eta^2}{12}G'(u).
\]
\end{proof}

\begin{proposition}[Uniform inverse heights]\label{glob:height}
For every selected branch $I$,
\begin{equation}\label{glob:height-bound}
 \left|\operatorname{Im}g_{I,f}\left(\frac{1+i}{2}\right)\right|
 \ge h_*:=\frac{\eta^2}{12M_1}>0,\qquad
 \eta=\min\{m_*/2,\,1/2\}.
\end{equation}
\end{proposition}
\begin{proof}
Let $J=f(I)$. Its finite endpoints are singular values by
Lemma~\ref{glob:ends}. If bounded, they are distinct by strict
monotonicity, so its length is at least $m_*$. Choose
\[
 u=
 \begin{cases}
  \text{the midpoint of }J,&J\text{ is bounded},\\
  2,&J=(v,\infty),\\
  -1,&J=(-\infty,v).
 \end{cases}
\]
These cases are exhaustive by Remark~\ref{glob:no-whole-line},
and give $u\in[-1,2]$ and $(u-\eta,u+\eta)\subset J$.
Choose $\sigma\in\{-1,1\}$ so that $G=\sigma g_{I,f}$ maps
$\HH$ to $\HH$. Then
\[
 G'(u)=|g_{I,f}'(u)|
 =\frac1{|f'(g_{I,f}(u))|}
 \ge\frac1{M_1},
\]
since the real point $g_{I,f}(u)$ belongs to the selected branch
and maps into $[-1,2]$. Lemma~\ref{glob:pick} proves the result.
This establishes Theorem~\ref{thm:additional-apriori}.
\end{proof}

\section{Why the normalization and fixed branch labels matter}
\label{glob:scope-section}

\subsection{A critical point and a fixed point alone are insufficient}

\begin{proposition}[A Gaussian degeneration]\label{glob:gaussian}
In the general real entire setting, fixing $c<p$ and requiring
$f'(c)=0$ and $f(p)=p$ does not imply uniform derivative or
inverse-height bounds, even with separated singular values and
an expansive finite marking with a fixed relative error constant.
\end{proposition}
\begin{proof}
Let
\[
 p=\frac12,\quad b=\frac14,\quad
 a_\varepsilon=p-b e^{-(p/\varepsilon)^2},\quad
 f_\varepsilon(z)=a_\varepsilon+b e^{-(z/\varepsilon)^2},
 \qquad\varepsilon>0.
\]
These maps lie in one real natural class by positive affine
source and target changes. Their unique real critical point is
$c=0$, and $f_\varepsilon'(0)=0$, $f_\varepsilon(p)=p$.

The only critical value of $e^{-z^2}$ is $1$, and its only finite
asymptotic value is $0$. To verify the latter, if it tended to a
finite nonzero value along a curve to infinity, choose a logarithm
near that value. The continuous function
$-z^2-\log(e^{-z^2})$ would then be a constant in
$2\pi i\mathbb Z$ on the terminal curve. This would bound $z^2$,
a contradiction. The value $0$ occurs along the real axis.
Therefore
\[
 S(f_\varepsilon)=\{a_\varepsilon,a_\varepsilon+b\},\qquad
 \sep S(f_\varepsilon)=b.
\]

Take
\[
 x=\left(0,a_\varepsilon,a_\varepsilon+\frac b4,
 a_\varepsilon+\frac b2,a_\varepsilon+\frac{3b}4,
 a_\varepsilon+b,1\right),
 \qquad
 \pi=(3,4,5,6,5,4,3),\quad \mathcal C=\{4\}.
\]
The ordered unscaled inverse images are
\begin{align*}
 (\widehat y_1,\ldots,\widehat y_7)
 =\varepsilon\bigl(&-\sqrt{\log4},-\sqrt{\log2},
 -\sqrt{\log(4/3)},0,\\
 &\sqrt{\log(4/3)},\sqrt{\log2},\sqrt{\log4}\bigr).
\end{align*}
Substitution verifies all marked equations. The two edges incident
to label $4$ hit it immediately, and every other edge maps in
one step to $[3,4]$ or $[4,5]$. Thus the marking is expansive
with hitting time at most one.

Since $1/4<a_\varepsilon<1/2$, the shortest target gap is
$m=b/4=1/16$. The normalized vector $y$ is independent of
$\varepsilon$, and $x,y\in W$ give
\[
 \|x-y\|_\infty\le1=16m.
\]
Hence the relative error condition holds with $K=16$.
Nevertheless,
\[
 \sup_{t\in\R}|f_\varepsilon'(t)|
 =\frac{\sqrt2\,b}{\varepsilon}e^{-1/2}\longrightarrow\infty,
 \qquad |f_\varepsilon''(0)|=2b/\varepsilon^2.
\]
The entire real image is contained in $[0,1]$.

The right-hand inverse is
\[
 g_{+,\varepsilon}(w)=
 \varepsilon\sqrt{-\log((w-a_\varepsilon)/b)},
\]
with branches giving the positive real inverse on
$(a_\varepsilon,a_\varepsilon+b)$; the left inverse is its
negative. At $w_*=(1+i)/2$,
\[
 g_{+,\varepsilon}(w_*)=
 \varepsilon\sqrt{-\log(e^{-(p/\varepsilon)^2}+2i)}.
\]
The square-root factor tends to a finite nonreal limit, so both
inverse heights are of order $\varepsilon$ and tend to zero.
Finally,
\[
 f_\varepsilon(z)-f_\varepsilon(0)
 =-\frac b{\varepsilon^2}z^2+O(z^4).
\]
The fixed critical-jet normalization excludes this source dilation.
\end{proof}

\begin{remark}[A sufficient extra marking for a fixed-point anchor]
\label{glob:marked-fixed-point}
The critical-point/fixed-point normalization can be used in this
general argument if $p$ is a fixed transported source label,
$f(p)=p$, and $p$ is also a finite target mark. After gap separation,
its image is singular or at distance at least $m_*$ from the
singular set. Include $p$ in the finite core, and use the fixed
distance $p-c$ in the two-point conformal-radius argument.
The critical budget remains uniform because the source label of
$p$ is fixed.

The source-label condition matters when there are infinitely many
real critical points. Merely requiring $f(p)=p$ and that $p$ be
a target mark does not bound the number of critical points between
$p$ and the other selected source labels.
\end{remark}

\subsection{Uniformity cannot extend to arbitrarily changing branch labels}

\begin{example}\label{glob:changing-branches}
Consider the single fixed map
\[
 f(z)=\frac{1+\sin(z^2)}2,\qquad S(f)=\{0,1/2,1\}.
\]
The critical values follow from $f'(z)=z\cos(z^2)$: the critical
point $0$ has value $1/2$, and the others have value $0$ or $1$.
There are no finite asymptotic values. Indeed, if $\sin(z^2)$
tended to a finite value $L$ along a curve $z\to\infty$, then
$u=e^{iz^2}$ would satisfy
\[
 u^2-2i\sin(z^2)u-1=0.
\]
The limiting roots are finite and nonzero. If distinct, continuity
forces $u$ eventually to remain near one of them; if repeated,
both roots tend to the same nonzero point. Thus $u$ tends to a
finite nonzero limit. A logarithm near that limit makes
$iz^2-\log(e^{iz^2})$ a constant in $2\pi i\mathbb Z$ on the
terminal curve, contradicting $|z^2|\to\infty$.

Let $\arcsin$ be the inverse associated with the real sine branch
$(-\pi/2,\pi/2)$, continued to $\HH$. For $j\ge1$,
\[
 g_j(w)=\sqrt{2\pi j+\arcsin(2w-1)}
\]
is the inverse of a positive real branch of $f$. The square root
is positive on the associated real target interval. Put
$\alpha=\log(1+\sqrt2)$, so $\arcsin(i)=i\alpha$. At
$w_*=(1+i)/2$,
\[
 \operatorname{Im}g_j(w_*)
 =\operatorname{Im}\sqrt{2\pi j+i\alpha}
 \sim\frac{\alpha}{2\sqrt{2\pi j}}\longrightarrow0.
\]
Also, at $t_j=\sqrt{2\pi j}$,
\[
 f(t_j)=\frac12,\qquad f'(t_j)=\sqrt{2\pi j}\longrightarrow\infty.
\]
The function, its singular values, and any chosen critical anchors
are fixed (the critical jet at $0$ has coefficient $1/2$).
Thus a fixed finite collection of branch labels is indispensable
when the real critical set is infinite.
\end{example}

% Appendix E (existence beyond entire maps), version 3 of 21 Nov 2026; v2 is appendix-axiomatic-existence-v2.tex
% Appendix E, version 3 (21 Nov 2026): existence beyond entire maps.
% Used by thurston_real_8nov2026-v3.tex via \input{appendix-axiomatic-existence-v3.tex}.
% Version 2 (appendix-axiomatic-existence-v2.tex) is unchanged.
% Reorganised after a report on version 2:
%   - one organising question and three subsections: the derivative criterion; finite-type maps
%     on strips; real estimates and their limitation;
%   - Lemma E.1 of v2 (N <= s-2) is now item (1) of the remark after the theorem;
%   - the theorem is sharpened: one-sided Taylor gives ||F'|| <= (9/4) sqrt(Phi), so
%     m Phi(m)^{N+1} -> 0 suffices (v2: exponent 2N+1); polynomial tails now need p > 2N+2
%     (v2: 4N+2);
%   - Proposition E.4 (finite type) has a proof in prose, and the modular bump follows it;
%   - Corollary E.6: the end-limit gap estimate is applied only on the lap with that end limit;
%   - one C^2 example (the Matern kernel), and the family h_r as Proposition E.8, with the
%     common strip |Im z| < 1/sqrt 2;
%   - cut: the other one-hump and several-singular-value examples, polynomials and Markov's
%     inequalities, the Epstein class, the population models, the x^2/(1+x^2) discussion
%     (period counts, the Douady-Hubbard/Tan Lei route), the catalogue of tail exponents, and the
%     rational realisation of the marking (E.2).
% Bibitems no longer cited by the appendix: Ama77, BE95, Has75, MSS73, Tan92, WC72.

\section{Existence beyond entire maps}\label{app:axiomatic}

How much of the existence proof survives when the maps are no longer entire? The proof of \cref{thm:realization-new} uses that the maps are entire at only two points:
\begin{itemize}
\item \cref{thm:natural-parameters}, which supplies a map with every prescribed vector of singular values;
\item \cref{ana:diameter,ana:scaled}, which bound the rescaled maps $F_x$ of \eqref{eq:affine}.
\end{itemize}
Everything else is real or combinatorial. This appendix gives three answers.
\begin{itemize}
\item The abstract argument needs only a bound on the rescaled second derivative (\cref{app:criterion}).
\item Finite-type maps on strips and other simply connected domains inherit that bound from the complex proof (\cref{app:strips}).
\item Some families of $C^2$ maps satisfy it by real estimates. Without such a bound, realisation can fail, even for a family of maps holomorphic on a common strip (\cref{app:real}).
\end{itemize}

Fix a labelled real topological model and a finite marking $X$ as in \cref{thm:realization-new}, saturated as in \cref{lem:critical-saturation}. Throughout this appendix
\[
 s=\#X
\]
is the number of marked points: the prescribed postsingular set, any further prescribed points, and the real critical points in $\operatorname{conv}X$ adjoined in \cref{lem:critical-saturation}. Then $(\pi,\cC)$ satisfies \eqref{eq:noflat} and is expansive. Let $N$ be its hitting-time bound as in \cref{def:pi-expansive}, and let $W$ and $v(x)$ be as in \cref{sec:thurston}.

\subsection{A derivative criterion for realisation}\label{app:criterion}

\begin{definition}[Full labelled family]\label{def:full-family}
A \emph{full labelled family} is a family $\{f_v\}$ of real $C^2$ maps on open intervals, defined for every singular-value vector $v$ that occurs as $v(x)$, $x\in W$, with the following properties.
\begin{enumerate}[label=(\alph*)]
\item $f_v=\varphi_v\circ f_0\circ\psi_v$ for increasing homeomorphisms $\varphi_v,\psi_v$, where $f_0$ realises the model. Every turning value, finite end limit and finite boundary value of $f_0$ lies in a finite labelled set $V^0$, and $\varphi_v(V^0)=v$, label by label.
\item The critical points of $f_0$ form a closed discrete set. $f_v'$ vanishes at the labelled critical points $\psi_v^{-1}(c)$, $c\in\Crit f_0$, and $f_v$ is strictly monotone between consecutive ones.
\item $v\mapsto f_v$ is continuous for locally uniform convergence, and the labelled critical points depend continuously on $v$.
\end{enumerate}
\end{definition}

\begin{theorem}[A derivative criterion]\label{thm:axiomatic}
Let $\{f_v\}$ be a full labelled family. Define $T$ and $F_x=f_{v(x)}\circ M_x$ as in \cref{def:thurston-new}. Suppose that, for a nonincreasing function $\Phi\ge1$,
\begin{equation}\label{eq:axiom-bound}
 \norm{F_x''}_{[0,1]}\le\Phi\bigl(m(x)\bigr)\quad(x\in W),
 \qquad
 m\,\Phi(m)^{N+1}\longrightarrow0\quad(m\downarrow0).
\end{equation}
Then $T$ has a fixed point $x^*$. The map $F=f_{v(x^*)}\circ M_{x^*}$ realises the marking: $F(x^*_i)=x^*_{\pi(i)}$ and $F'(x^*_j)=0$ for $j\in\cC$, with the prescribed branch labels.
\end{theorem}

\begin{proof}
\begin{itemize}
\item The proofs of \cref{lem:thurston-continuous,lem:thurston-realisation} use only (a)--(c). So $T$ is continuous, and a fixed point realises the marking.
\item The proof of \cref{ana:range} is real, so $F_x([0,1])\subset[0,1]$.
\item Write $\Phi=\Phi(m(x))$ and $h=1/(2\sqrt\Phi)\le\tfrac12$. From each $t\in[0,1]$, a step of length $h$ stays in $[0,1]$ in at least one direction. Taylor's theorem along that step gives
\[
 \abs{F_x'(t)}\le\frac1h+\frac h2\,\Phi=\frac94\sqrt\Phi .
\]
\item Substituting $M_2^F\le\Phi$ and $1+M_1^F\le\tfrac{13}4\sqrt\Phi$ into \eqref{eq:gapinequality} gives
\[
 1\le C_{K,N}\,m\,\Phi(m)^{N+1},\qquad C_{K,N}=\tfrac12(1+2K)^2\bigl(\tfrac{13}4\bigr)^{2N},
\]
along configurations with $\norm{T(x)-x}_\infty\le K\,m(x)$. By \eqref{eq:axiom-bound} this forces $m(x)\ge m_*(K)>0$.
\item The proof of \cref{cor:existence-new} now applies verbatim.
\end{itemize}
\end{proof}

\begin{remark}\label{rem:axiomatic}
\begin{enumerate}[label=(\arabic*)]
\item \emph{Hitting time.} One may take $N\le s-2$, so $m\Phi(m)^{s-1}\to0$ suffices. Indeed, by \eqref{eq:noflat} the map $G$ of \cref{def:pi-expansive} sends each unit edge onto a union of unit edges. Hence $G^n(e)$ is the union of the unit edges reached from $e$ in $n$ steps of the directed graph with $e\to e'$ when $e'\subset G(e)$. So the first hitting time in \eqref{eq:expansive} is the length of a shortest path from $e$ to an edge with an endpoint in $\cC$. Such a path visits distinct edges, and all of them except the last avoid $\cC$. At least one of the $s-1$ unit edges has an endpoint in $\cC$, so the length is at most $s-2$.
\item \emph{Admissible bounds.} Any $\Phi(m)=Cm^{-\alpha}$ with $(N+1)\alpha<1$ suffices. Condition \eqref{eq:axiom-bound} holds for every marking exactly when $\Phi(m)=m^{-o(1)}$. Indeed, if $\Phi\le C_\eps m^{-\eps}$ for every $\eps>0$, take $\eps<1/(N+1)$. Conversely, if $m\Phi(m)^k\to0$ for every $k$, then $\Phi(m)\le m^{-1/k}$ for small $m$. Compact families, with uniformly bounded spans and $C^2$ norms, have $\Phi$ constant, even with flat critical points; this is the setting of \cite[Section~II.4]{dMvS93}. For real entire natural families, \cref{ana:scaled} gives $\Phi$ polylogarithmic.
\item \emph{Only the spans matter.} The bound is needed only on the spans $[\widehat y_1,\widehat y_s]$. Over arbitrary intervals it fails even for entire maps: for $\tfrac12(1+\sin x)$ on $[-\tfrac\pi2,-\tfrac\pi2+2\pi k]$ the rescaled derivative is $\pi k$.
\item \emph{Existence only.} The theorem uses neither non-flatness of the critical points, nor the multiplicity bound of \cref{gap:interval}, nor a source normalisation, nor a sign of the Schwarzian derivative. Its conclusion is existence only.
  \begin{itemize}
  \item Let $f_v=M_t^{-1}\circ Q_{a(v)}\circ M_t$, where $Q_a(x)=4ax(1-x)$, $M_t(x)=x/(x+t(1-x))$ and $t=v(1-a(v))/\bigl(a(v)(1-v)\bigr)$. These are real-analytic unimodal maps with negative Schwarzian and critical value $v$. If $a(v)=v$ near $0$ and $1$ but $a$ is not monotone, a superstable combinatorics occurs at several $v$, and these give distinct fixed points of $T$. So \cref{sec:convergence-new} needs complex structure.
  \item A $C^\infty$ map equal to the identity on $[0,1]$, with turning values $0$ and $1$, has $[0,1]$ as a real Levy cycle. So \cref{prop:levy-necessary} needs analyticity.
  \end{itemize}
\end{enumerate}
\end{remark}

\subsection{Finite-type maps on strips and simply connected domains}\label{app:strips}

\begin{proposition}[Finite-type maps on simply connected domains]\label{prop:finite-type}
Let $\Omega_0\supset\R$ be a simply connected domain invariant under complex conjugation, for instance a strip $\{\abs{\Im z}<\delta\}$. Let $f_0:\Omega_0\to\C$ be a real-symmetric holomorphic map with a finite set $S\subset\R$ such that $f_0:\Omega_0\setminus f_0^{-1}(S)\to\C\setminus S$ is a covering. Then the real natural family \eqref{eq:natural} of $f_0$, with $f=\phi\circ f_0\circ\psi$ on $\psi^{-1}(\Omega_0)$, contains a full labelled family satisfying \eqref{eq:axiom-bound} with $\Phi(m)=C(1+\log(1/m))^B$.
\end{proposition}

\begin{proof}
For the parameters, run the proof of \cref{thm:natural-parameters} with one change: pull the Beltrami coefficient back by $f_0$ on $\Omega_0$ and extend it by $0$ to $\C\setminus\Omega_0$. The resulting maps $f_v=\phi\circ f_0\circ\psi$ are holomorphic on the simply connected symmetric domains $\Omega_v=\psi^{-1}(\Omega_0)\supset\R$. They are coverings over $\C\setminus\phi(S)$ off the preimage of $\phi(S)$, and they depend continuously on $v$. So (a)--(c) hold. The multiplicity bound of \cref{gap:interval} holds because the zeros of $f_0'$ are discrete in $\Omega_0$.

For the bound, note that \cref{ana:diameter,ana:scaled} use entireness only through \cref{ana:disk-topology,ana:singleton}. Let $U$ be the component of $f_v^{-1}D(0,4)$ containing the span.
\begin{itemize}
\item In \cref{ana:disk-topology}, the interior of a Jordan curve in $U$ lies in $\Omega_v$, because $\Omega_v$ is simply connected. So the maximum principle still applies, and $U$ is simply connected. Moreover $U\ne\C$, since $U\subset\Omega_v$, or by Liouville's theorem if $\Omega_v=\C$.
\item In \cref{ana:singleton}, the covering property is now a hypothesis. The universal-cover component $T$ still has unbounded real trace. Indeed, a finite endpoint of $T\cap\R$ is real, so it lies in $\Omega_v$, where $f_v$ is continuous. Its image therefore lies on $\partial D(v,R)$, and the sign argument in the proof of \cref{ana:singleton} excludes two finite endpoints.
\end{itemize}
Everything after these two lemmas takes place inside $U$. Sections~\ref{ana:metric}--\ref{ana:finite-interval} compare the hyperbolic metrics of $U$ and of preimage components contained in $U$, using that $f_v$ restricted to $U\setminus f_v^{-1}(S(f_v))$ is a covering of $D(0,4)\setminus S(f_v)$. The Koebe--Cauchy argument of \cref{ana:scaled} uses only the disk $D(t,R_U(t)/4)\subset U$, on which $\abs{f_v}<4$. So \cref{ana:diameter,ana:scaled} hold with the same constants, and \cref{ana:scaled} gives the bound.
\end{proof}

\begin{example}[A modular bump]\label{ex:modular}
Let $\lambda$ be the elliptic modular function on the upper half-plane \cite{Ahl79}, and put $u(z)=\lambda(ie^z)$ on the strip $\Omega_0=\{\abs{\Im z}<\pi/2\}$.
\begin{itemize}
\item The map $z\mapsto ie^z$ sends $\Omega_0$ conformally onto the upper half-plane, and $\lambda$ is a universal covering of $\C\setminus\{0,1\}$. So $u:\Omega_0\to\C\setminus\{0,1\}$ is a real-symmetric universal covering. The boundary of the strip is a natural boundary for $u$, so $u$ is not the restriction of an entire function.
\item On $\R$, $u$ decreases from $1$ to $0$, with $u(0)=\tfrac12$, $u(-x)=1-u(x)$ and $u(x)\sim16\exp(-\pi e^{x})$ as $x\to+\infty$.
\end{itemize}
Hence $f_0=4u(1-u)$ is an even bump with critical value $1$ at $0$ and end limit $0$ at both ends, a non-entire analogue of $e^{-x^2}$. The polynomial $P(t)=4t(1-t)$ is a covering over $\C\setminus\{0,1\}$ off $P^{-1}\{0,1\}$. So $f_0=P\circ u$ satisfies \cref{prop:finite-type} with $S=\{0,1\}$. Thus every expansive marking of its model is realised in its natural family, although $f_0$ has a natural boundary. The same holds for $P\circ u$ with any real polynomial $P$ with real critical values, with $S=P(\Crit P)\cup\{P(0),P(1)\}$.
\end{example}

\subsection{Real estimates and their limitation}\label{app:real}

The next result checks \eqref{eq:axiom-bound} by real means. The idea is to keep the map fixed and vary only the target.

\begin{corollary}[A fixed map with varying target]\label{cor:fixed-source}
Let $f_0$ be a $C^2$ map of $\R$ which realises the model, with closed discrete critical set and finitely many turning values and finite end limits, all in $V^0$.
\begin{itemize}
\item Suppose that $f_0'$ and $f_0''$ are bounded on every lap with a finite end limit.
\item Suppose that $\abs{f_0(x)-a}\le\rho(\abs{x})$ on each such lap, where $a$ is its limit and $\rho$ is decreasing with $\rho(t)\to0$.
\item Let $\varphi_v$ be increasing $C^2$ diffeomorphisms of $\R$, continuous in $v$, with $\varphi_v(V^0)=v$, slope $1$ outside a fixed compact set, and $\sup_v(\norm{\varphi_v'}_\infty+\norm{\varphi_v''}_\infty)<\infty$.\footnote{For instance, $\varphi_v'=\eps_v+\sum_k(\Delta v_k-\eps_v\Delta V^0_k)\beta_k$ on $[V^0_1,V^0_q]$. Here the $\beta_k\ge0$ are fixed $C^1$ bumps supported in $(V^0_k,V^0_{k+1})$ with $\int\beta_k=1$, $\Delta$ denotes consecutive differences, and $\eps_v=\min_k\Delta v_k/(2\Delta V^0_k)$; outside $[V^0_1,V^0_q]$, $\varphi_v'$ rises from $\eps_v$ to $1$ over fixed intervals. Since $\Delta v_k\le1$, the bounds are uniform in $v$. For $q=2$ one may take $\varphi_v$ affine.}
\end{itemize}
Then $f_v=\varphi_v\circ f_0$ is a full labelled family satisfying \eqref{eq:axiom-bound} with $\Phi(m)=C\bigl(1+\rho^{-1}(m/C)\bigr)^2$. In particular \cref{thm:axiomatic} applies when $\rho$ decays exponentially, with $\Phi(m)\asymp\log^2(1/m)$, and when $\rho(t)=Ct^{-p}$ with $p>2N+2$, in particular when $p>2s-2$.
\end{corollary}

\begin{proof}
\begin{itemize}
\item Properties (a)--(c) are clear, with $\psi_v=\mathrm{id}$. The targets $t_i=\varphi_v^{-1}(x_{\pi(i)})$ lie in a fixed compact set.
\item Suppose that $\widehat y_i$ lies on an open lap with finite end limit $a$. That lap is mapped strictly monotonically onto an open interval with endpoint $a$, so $t_i\ne a$. Hence $x_{\pi(i)}$ and $\varphi_v(a)$ are distinct marks, and $\abs{t_i-a}\ge m/\norm{\varphi_v'}_\infty$. Since $\abs{t_i-a}=\abs{f_0(\widehat y_i)-a}\le\rho(\abs{\widehat y_i})$, this gives $\abs{\widehat y_i}\le\rho^{-1}(m/C)$.
\item Every other labelled preimage lies on a bounded lap or on a lap tending to an infinite limit, hence in a fixed compact set. So the span has length $L\le C\bigl(1+\rho^{-1}(m/C)\bigr)$.
\item Finally $F_x''=L^2\bigl(\varphi_v''(f_0')^2+\varphi_v'f_0''\bigr)\circ M_x$, which is bounded by $CL^2$.
\end{itemize}
\end{proof}

\begin{example}[A kernel that is $C^2$ but not $C^3$]\label{ex:matern}
The Mat\'ern kernel $f_0(x)=(1+\abs{x})e^{-\abs{x}}$ of Gaussian-process regression \cite{RW06} has a single critical point, a maximum at $0$ with $f_0''(0)=-1$. It is $C^2$ but not $C^3$ there, since $f_0'''(0^\pm)=\pm2$. Its end limit is $0$ at both ends, approached exponentially, and $f_0'$ and $f_0''$ are bounded. Here $q=2$, so $\varphi_v$ may be taken affine. By \cref{cor:fixed-source}, every expansive marking of this model is realised by a map $x\mapsto A+Bf_0(cx+d)$ with $B,c>0$.
\end{example}

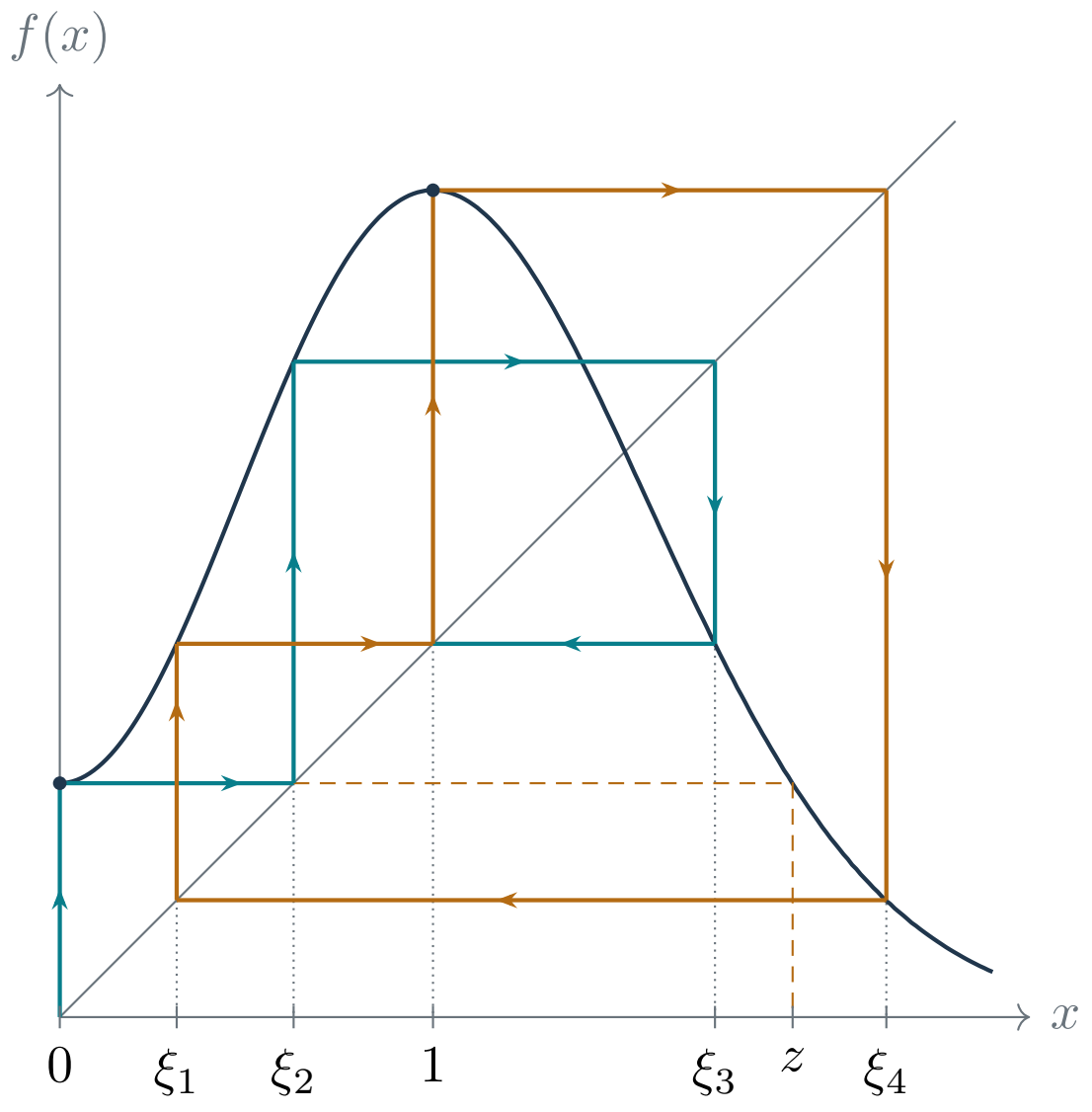
\begin{figure}[htbp]
\centering
\begin{tikzpicture}[x=2.5cm,y=2.5cm,font=\small,
  mid/.style={postaction={decorate},decoration={markings,
    mark=at position #1 with {\arrow{Stealth[length=4pt]}}}},mid/.default=0.55]
 \def\xa{0.31323}\def\xb{0.62649}\def\xc{1.75555}\def\xd{2.21481}\def\zz{1.96362}
 \draw[->,softgray] (0,0)--(2.6,0) node[right] {$x$};
 \draw[->,softgray] (0,0)--(0,2.5) node[above] {$f(x)$};
 \draw[softgray] (0,0)--(2.4,2.4);
 \foreach \a in {\xa,\xb,1,\xc,\xd}{\draw[softgray,densely dotted] (\a,\a)--(\a,0);}
 \draw[ink,thick] plot[domain=0:2.5,samples=220]
   (\x,{5.355887*(0.116972+0.883028*\x*\x)*exp(-0.883028*\x*\x)});
 \draw[ochre,densely dashed] (\xb,\xb)--(\zz,\xb)--(\zz,0);
 % \p, \q: arrow positions on the vertical and horizontal segment (moved off crossings)
 \foreach \a/\b/\p/\q in {0/\xb/0.55/0.78,\xb/\xc/0.55/0.55,\xc/1/0.55/0.55}{
   \draw[teal,thick,mid=\p] (\a,\a)--(\a,\b);
   \draw[teal,thick,mid=\q] (\a,\b)--(\b,\b);}
 \foreach \a/\b/\p/\q in {1/\xd/0.55/0.55,\xd/\xa/0.55/0.55,\xa/1/0.78/0.8}{
   \draw[ochre,thick,mid=\p] (\a,\a)--(\a,\b);
   \draw[ochre,thick,mid=\q] (\a,\b)--(\b,\b);}
 \foreach \a/\l in {0/0,\xa/\xi_1,\xb/\xi_2,1/1,\xc/\xi_3,\zz/z,\xd/\xi_4}{
   \draw[softgray] (\a,0.03)--(\a,-0.03);
   \node[below] at (\a,-0.03) {$\l$};}
 \fill[ink] (0,\xb) circle (1.3pt) (1,\xd) circle (1.3pt);
\end{tikzpicture}
\caption{Cobweb of the marking \eqref{eq:fat-marking}, realised by the entire map $f(x)=\lambda(\eps+(1-\eps)x^2)e^{-(1-\eps)x^2}$ with $\eps\approx0.11697$ and $\lambda\approx5.3559$. Only $x\ge0$ is shown; $f$ is even. The orbit of the critical point $0$ (teal) lands on the cycle $1\mapsto\xi_4\mapsto\xi_1\mapsto1$ of the critical point $1$ (ochre). The dashed line marks the point $z>1$ with $f(z)=f(0)$. Realisation requires $z<f(1)=\xi_4$, while for the maps $\lambda h_r$ of \cref{prop:fat-tail} with $f(0)<1$ one always has $f(1)<1/r<z$.}
\label{fig:cobweb-fat}
\end{figure}

The tail hypothesis cannot simply be dropped, even for maps that are holomorphic on a common strip. Consider the even model with a minimum at $0$, maxima at $\pm1$ and end limit $0$ at both ends, and the marking
\begin{equation}\label{eq:fat-marking}
 0\mapsto\xi_2\mapsto\xi_3\mapsto1\mapsto\xi_4\mapsto\xi_1\mapsto1,
 \qquad 0<\xi_1<\xi_2<1<\xi_3<\xi_4 .
\end{equation}
So the end limit sits at the critical point $0$, $s=6$, $\pi=(3,4,5,6,4,2)$ and $\cC=\{1,4\}$. This marking is expansive with $N=1$. By \cref{thm:realization-new} it is realised in the natural family of every real entire map with this model; \cref{fig:cobweb-fat} shows one realisation.

\begin{proposition}[A full labelled family that misses an expansive marking]\label{prop:fat-tail}
For $0<r<1$ let
\[
 h_r(x)=\frac{r+\beta x^2}{(1+\gamma x^2)^{3/2}},
 \qquad \sqrt{1+\gamma}\,(2-\gamma)=2r,\quad \beta=\frac{3\gamma r}{2-\gamma},
\]
where $\gamma\in(0,2)$ is unique because the left side decreases from $2$ to $0$ there.
\begin{enumerate}[label=(\roman*)]
\item Every $h_r$ is even, positive and holomorphic on the common strip $\abs{\Im z}<1/\sqrt2$. Its only real critical points are $0$ and $\pm1$, with $h_r(0)=r$ and $h_r(\pm1)=1$, and $h_r(x)\sim\beta\gamma^{-3/2}\abs{x}^{-1}$.
\item The maps $c+\lambda h_r$ ($c\in\R$, $\lambda>0$, $0<r<1$) form a full labelled family. Each vector of singular values occurs exactly once, and the family is closed under post-composition with increasing affine maps.
\item No map in this family realises \eqref{eq:fat-marking}. Hence its Thurston map has no fixed point, and the family violates \eqref{eq:axiom-bound} for every $\Phi$ with $m\Phi(m)^2\to0$.
\end{enumerate}
\end{proposition}

\begin{proof}
(i) The only singularities of $h_r$ are the branch points $\pm i/\sqrt\gamma$, and $1/\sqrt\gamma>1/\sqrt2$. On the strip, $1+\gamma z^2$ avoids $(-\infty,0]$, so the principal branch of $(1+\gamma z^2)^{3/2}$ is holomorphic there. Moreover
\[
 h_r'(x)=\beta\gamma\,x(1-x^2)(1+\gamma x^2)^{-5/2},
\]
and $h_r(\pm1)=1$ because $r+\beta=2r(1+\gamma)/(2-\gamma)=(1+\gamma)^{3/2}$.

(ii) A vector of singular values consists of an end limit $v_1$, a minimum value $v_2$ and a maximum value $v_3$, with $v_1<v_2<v_3$. It is realised exactly by $c=v_1$, $\lambda=v_3-v_1$ and $r=(v_2-v_1)/(v_3-v_1)$. All these maps have the same lap structure, so $\psi_v$ can be defined lap by lap, and (a)--(c) follow from (i) and the continuity of $\gamma$ and $\beta$ in $r$.

(iii) First, $h_r(1/r)>r$ for every $r\in(0,1)$. Put $\tau=\sqrt{1+\gamma}\in(1,\sqrt3)$. Then $r=\tau(3-\tau^2)/2$, $\gamma=\tau^2-1$ and $\beta=\tfrac32\tau(\tau^2-1)$. Since $h_r(1/r)=r(r^3+\beta)/(r^2+\gamma)^{3/2}$, the claim is $r^3+\beta>(r^2+\gamma)^{3/2}$, and
\[
 (r^3+\beta)^2-(r^2+\gamma)^3=\tfrac1{16}(2-\tau)(\tau-1)^2(\tau+1)^2(\tau+2)\,Q(\tau^2),
\]
where $Q(1+w)=3w^2(w^2-w+2)+12w+4>0$ for $w\ge0$.

Now suppose that a map of the family realises \eqref{eq:fat-marking}. Its end limit is the critical point $0$, so it is $f=\lambda h_r$. Then $f(0)=\lambda r=\xi_2<1<\xi_4=\lambda=f(1)$, so $1<\lambda<1/r$. As $h_r$ decreases on $[1,\infty)$,
\[
 f^2(1)=\lambda h_r(\lambda)>\lambda h_r(1/r)>\lambda r=f(0),
\]
whereas \eqref{eq:fat-marking} requires $f^2(1)=\xi_1<\xi_2=f(0)$. A fixed point $x^*$ of $T$ would give such a realisation, namely $M_{x^*}\circ f_{v(x^*)}$, which is $F$ conjugated by $M_{x^*}$ and lies in the family by (ii). So $T$ has no fixed point. Since $N=1$, \cref{thm:axiomatic} shows that \eqref{eq:axiom-bound} fails for every $\Phi$ with $m\Phi(m)^2\to0$.
\end{proof}

So \cref{prop:finite-type} needs its covering hypothesis: a common strip of holomorphy for a full labelled family does not by itself guarantee realisation. Each $h_r$ is holomorphic on a neighbourhood of the closed strip $\abs{\Im z}\le1/\sqrt2$, so its values on the boundary lines are asymptotic values, and $h_r$ is not of finite type on the strip.

% Cobweb of the entire realisation of (E.2) (eps = 0.1169723509, lambda = 5.355886551):
% xi1 = 0.31323, xi2 = 0.62649, xi3 = 1.75555, xi4 = 2.21481, z = 1.96362.
\FloatBarrier


\begin{thebibliography}{99}

\bibitem{Ahl79}
Lars V. Ahlfors,
\emph{Complex Analysis}, third edition,
McGraw-Hill, New York, 1979.

\bibitem{BFR}
Walter Bergweiler, N\'uria Fagella, and Lasse Rempe-Gillen,
\emph{Hyperbolic entire functions with bounded Fatou components},
Commentarii Mathematici Helvetici \textbf{90} (2015), no.~4, 799--829.
\href{https://doi.org/10.4171/CMH/371}{doi:10.4171/CMH/371}.

\bibitem{BMS}
Araceli Bonifant, John Milnor, and Scott Sutherland,
\emph{The W. Thurston algorithm applied to real polynomial maps},
Conformal Geometry and Dynamics \textbf{25} (2021), 179--199.
\url{https://arxiv.org/abs/2005.07800}.

\bibitem{dFdMSV26}
Edson de Faria, Welington de Melo, Pedro A. S. Salom\~ao, and Edson Vargas,
\emph{Blaschke-type models for multimodal circle maps},
preprint, 2026, arXiv:2605.05823.
\url{https://arxiv.org/abs/2605.05823}.

\bibitem{dMvS93}
Welington de Melo and Sebastian van Strien,
\emph{One-Dimensional Dynamics},
Ergebnisse der Mathematik und ihrer Grenzgebiete, third series, vol.~25,
Springer-Verlag, Berlin, 1993.
\href{https://doi.org/10.1007/978-3-642-78043-1}
{doi:10.1007/978-3-642-78043-1}.

\bibitem{DH93}
Adrien Douady and John H. Hubbard,
\emph{A proof of Thurston's topological characterization of rational functions},
Acta Mathematica \textbf{171} (1993), no.~2, 263--297.

\bibitem{Duren}
Peter L. Duren,
\emph{Univalent Functions},
Grundlehren der mathematischen Wissenschaften, vol.~259,
Springer-Verlag, New York, 1983.

\bibitem{EL92}
Alexandre E. Eremenko and Mikhail Yu. Lyubich,
\emph{Dynamical properties of some classes of entire functions},
Annales de l'Institut Fourier \textbf{42} (1992), no.~4, 989--1020.

\bibitem{FvS25}
Gustavo R. Ferreira and Sebastian van Strien,
\emph{Holomorphic motions, natural families of entire maps, and
multiplier-like objects for wandering domains},
Mathematische Annalen \textbf{392} (2025), 701--732.
\href{https://doi.org/10.1007/s00208-025-03107-8}
{doi:10.1007/s00208-025-03107-8}.

\bibitem{Gar87}
Frederick P. Gardiner,
\emph{Teichm\"uller Theory and Quadratic Differentials},
Pure and Applied Mathematics, Wiley-Interscience, New York, 1987.

\bibitem{HS94}
John H. Hubbard and Dierk Schleicher,
\emph{The spider algorithm},
in \emph{Complex Dynamical Systems} (Cincinnati, OH, 1994),
Proceedings of Symposia in Applied Mathematics, vol.~49,
American Mathematical Society, Providence, RI, 1994, 155--180.

\bibitem{HSS09}
John Hubbard, Dierk Schleicher, and Mitsuhiro Shishikura,
\emph{Exponential Thurston maps and limits of quadratic differentials},
Journal of the American Mathematical Society \textbf{22} (2009), no.~1,
77--117.
\href{https://doi.org/10.1090/S0894-0347-08-00609-7}
{doi:10.1090/S0894-0347-08-00609-7}.

\bibitem{LV73}
Olli Lehto and Kaarlo I. Virtanen,
\emph{Quasiconformal Mappings in the Plane}, second edition,
Grundlehren der mathematischen Wissenschaften, vol.~126,
Springer-Verlag, New York--Heidelberg, 1973.

\bibitem{MB} 
Helena, Mihaljević-Brandt, \emph{  A landing theorem for dynamic rays of geometrically finite entire functions}  J. Lond. Math. Soc. (2) \textbf{81} (2010), no. 3, 696--714.

\bibitem{MT88}
John Milnor and William Thurston,
\emph{On iterated maps of the interval},
in \emph{Dynamical Systems} (College Park, MD, 1986--87),
Lecture Notes in Mathematics, vol.~1342,
Springer-Verlag, Berlin, 1988, 465--563.

\bibitem{Pro24}
Nikolai Prochorov,
\emph{Finite and infinite degree Thurston maps with a small postsingular set},
preprint, 2024, arXiv:2410.01146.
\url{https://arxiv.org/abs/2410.01146}.

\bibitem{RW06}
Carl Edward Rasmussen and Christopher K. I. Williams,
\emph{Gaussian Processes for Machine Learning},
MIT Press, Cambridge, MA, 2006.

\bibitem{Smania05}
Daniel Smania,
\emph{Phase space universality for multimodal maps},
Bulletin of the Brazilian Mathematical Society, New Series
\textbf{36} (2005), no.~2, 225--274.
\href{https://doi.org/10.1007/s00574-005-0038-y}
{doi:10.1007/s00574-005-0038-y}.
The fixed-point criterion is Proposition~8.1, p.~269.

 \bibitem{vS-complex}
 Sebastian van Strien,
 \emph{The Thurston algorithm for general entire maps on the complex plane:
 existence and uniqueness},
 preprint, 2026, arXiv:2611.27925.
 \url{https://arxiv.org/pdf/2611.27925}.

\end{thebibliography}
\end{document}